\documentclass[12pt]{amsart}

\usepackage[utf8]{inputenc}
\usepackage[english]{babel}
\usepackage{polski}
\usepackage{amsthm,amsfonts,amsmath, amssymb, mathtools}
\usepackage{csquotes}
\usepackage{indentfirst}
\usepackage{comment}
\usepackage{enumerate}
\usepackage{cancel}
\usepackage{caption}
\usepackage{subcaption}
\usepackage{graphicx}
\usepackage{tikz}
\usepackage{quiver}
\usepackage[all,cmtip]{xy}
\usepackage{url}
\usepackage{hyperref}
\hypersetup{
    colorlinks=true,
    linkcolor=black, % default = red
    filecolor=black,      
    urlcolor=black,
    citecolor=black, % default = green 
    pdftitle={Quadratic Hubbard trees, automata, and ray connections}, % change!
    pdfpagemode=FullScreen,
    }
\usepackage[style = alphabetic, backend = biber, doi = false, isbn = false, url = false, eprint = false, date = year, maxnames = 99]{biblatex}
\renewbibmacro{in:}{} % Removes "In:" before journals in the citation
\usepackage[margin=1in]{geometry}
\usepackage[colorinlistoftodos, tickmarkheight=0.2cm,
textwidth=\marginparwidth-.25cm]{todonotes}

\newcommand{\CC}{\mathbb{C}}     
\newcommand{\DD}{\mathbb{D}}

\newcommand{\NN}{\mathbb{N}}

\newcommand{\QQ}{\mathbb{Q}}     
\newcommand{\RR}{\mathbb{R}}     
\DeclareMathOperator{\R}{\mathcal{R}}
\DeclareMathOperator{\Char}{\mathrm{Char}}

\newtheorem{thm}{Theorem}[section]
\newtheorem{lemma}[thm]{Lemma}
\newtheorem{cor}[thm]{Corollary}
\newtheorem{prop}[thm]{Proposition}
\theoremstyle{definition}
\newtheorem{exm}[thm]{Example}
\newtheorem{definition}[thm]{Definition}
\newtheorem{remark}[thm]{Remark}
\newtheorem{question}{Question}

\title{Quadratic Hubbard trees, automata, and ray connections}
\author{Eduardo Ventilari Sodré}
\address{Department of Mathematics, Brown University}
\email{\url{eduardo_ventilari_sodre@brown.edu}}
\date{\today}

\begin{document}

\begin{abstract}
We describe an algorithm that, given a postcritically finite quadratic polynomial, outputs a finite state machine whose regular language is the set of angles whose dynamical rays land on the Hubbard tree, generalizing an example of Milnor. We do this by extracting orbit portraits from angled internal addresses, and constructing a suitable Markov partition of angle intervals. A main application is computing direct ray connections between Hubbard trees in matings of quadratic polynomials, and understanding when ray connections between the trees first appear when moving along veins of the Mandelbrot set. If one of the mating parameters is fixed, we prove that first ray connections are picked up at parabolic parameters or tips of veins. We also prove that the locus of parameter pairs in two given veins for which there are no ray connections between the trees is open. Additionally, we show that there are no ray connections between small Julia sets in matings of locally connected rabbit tunings. 
\end{abstract}
% Abstract must precede \maketitle with amsart

\maketitle

\setcounter{tocdepth}{1}
\tableofcontents

\section{Introduction}

Questions concering the mating construction for polynomial maps are now ubiquitous in one-dimensional complex dynamics. In one instance, Milnor explored several aspects of the mating $f_{c_{1/4}} \coprod f_{c_{1/4}}$, where $c_{1/4}$ is the landing point of the external parameter ray of the Mandelbrot set of with angle $1/4$, and $f_{c_{1/4}}$ is the corresponding postcritically finite polynomial \cite{Milnor2004}. In particular, he identified the ray connections passing through the Hubbard trees of the two Julia sets: all such ray connections are preperiodic, and eventually map to the ray connection at angle $0$ connecting the $\beta$ fixed points, or to a period $4$ cycle of rays corresponding to the orbit $\{1/5, 2/5, 4/5, 3/5\}$ under doubling. He was then able to uniformly bound the size and diameter of all ray connections, providing another proof of the existence of the topological mating without appealing to Tan Lei's criterion for quadratic polynomials \cite{Jung2017, Lei1992}.

Milnor does this by labeling the sides of the edges of the Hubbard tree of $f_{c_{1/4}}$, corresponding to external angle intervals, and encoding the dynamics of mappings between these intervals as a Markov partition. By assinging a bit value to each labeled edge-side according to whether it is above the spine $[\beta, -\beta]$ or not, the bit sequences represented by allowable infinite words in the resulting automaton correspond exactly to the binary expansions of angles whose external rays land on the Hubbard tree. This is essentially a ``double cover'' of the usual Markov partition of the Hubbard tree, for which one computes the core entropy of the polynomial \cite{Dudko2020, Tiozzo2015}; each edge of the Hubbard tree has two sides.

In this paper, we generalize this construction and describe an algorithmic implementation of it for all quadratic postcritically finite polynomials. The key idea is to compute the ``forbidden region'' of the Hubbard tree $H_c$, which corresponds to how the branches of the preimage Hubbard tree $f_c^{-1}(H_c)$ are attached to $H_c$, and which angle sectors envelop these branches. Knowing the forbidden region, we build a Markov partition and a finite state automaton, whose regular language is exactly the set of angles whose external rays land on $H_c$.

In order to find the forbidden region, we first find the collection of all satellite orbit portraits for $f_{c}$ such that $q \geq 3$ rays land at each periodic point. Geometrically, this corresponds to finding all the instances in which, on the ``path'' to the parameter $c$ within the Mandelbrot set $M$ starting from the main cardioid, we undergo a satellite bifurcation of rotation number $p/q$, where $q \geq 3$. The main combinatorial tool for this is the notion of (angled) internal addresses, as described by Schleicher \cite{Schleicher2016}.

The collection of all satellite orbit portraits for which $q \geq 3$ tells us what are the external arguments of all periodic branch points of the Hubbard tree $H$. By pulling these back along preimages, we find the ``attaching points'' of the preimage Hubbard tree, when dynamical sectors formed by adjacent angles in the circular ordering stop containing points of the postcritical set. This is what culminates in a description of the forbidden region of $H$, and allows us to construct the aforementioned automaton computationally.

As the main application of this algorithm, by ``mating'' two of these automata, we are able to compute the set of direct ray connections between the two Hubbard trees for all PCF quadratic polynomials. Several interesting phenomena arise from these computations. For instance, in the formal mating $f_{1/6} \sqcup f_{3/16}$, there are no direct ray connections between the two Hubbard trees, and this can be shown to imply there are also no indirect ones. Therefore the Hubbard trees ``survive'' in the geometric mating as disjoint invariant trees containing the postcritical set, even though the Julia set is the entire sphere. 

We also utilize the combinatorial nature of the set of direct ray connections between the Hubbard trees to describe when they first appear as we move along veins of $M$. This culminates in the following result:

\begin{thm} \label{theorem: first ray connection}
Let $c_1$ and $c_2$ be two postcritically finite parameters in the Mandelbrot set. Then exactly one of three cases below holds for the formal mating $f_{c_1} \sqcup f_{c_2}$.

\begin{enumerate}[1.]
    \item There are no ray connections between the Hubbard trees $H_{c_1}$ and $H_{c_2}$.

    \item The parameter $c_1$ is the tip of a vein, and all direct ray connections between the Hubbard trees $H_{c_1}$ and $H_{c_2}$ in the formal mating $f_{c_1} \sqcup f_{c_2}$ are of the form $\R_{f_{c_1} \sqcup f_{c_2}}(\phi)$, where $\phi$ is in the grand orbit of the external argument of $c_1$. For all parameters $c \prec c_1$, there are no ray connections between the generalized Hubbard tree $T_c$ and $H_{c_2}$ in the formal mating $f_c \sqcup f_{c_2}$.

    \item There exists a unique parabolic parameter $p_1 \preceq c_1$ such that, for all parameters $c \prec p_1$, there are no ray connections between the generalized Hubbard tree $T_{c}$ and $H_{c_2}$ in the formal mating $f_c \sqcup f_{c_2}$. If $P_M(\theta_-, \theta_+)$ is the parameter ray pair landing at $p_1$, then $\R_{f_{c_1} \sqcup f_{c_2}}(\theta_-)$ or $\R_{f_{c_1} \sqcup f_{c_2}}(\theta_+)$ is a direct ray connection between the Hubbard trees.
\end{enumerate}
\end{thm}

The notation $\R_{f_{c_1} \sqcup f_{c_2}}(\phi)$ means the union of the external ray at angle $\phi$ for $f_{c_1}$ with the external ray at angle $-\phi$ for $f_{c_2}$, along with their common point on the glued circles at infinity in the formal mating (see Section \ref{section: mating}). The vein of a PCF parameter in the Mandelbrot set is a canonical topological arc connecting it to the main cardioid, and the relation $c' \preceq c$ means that $c'$ is on the vein to $c$. The generalized Hubbard tree is the regular hull of the (possibly infinite) postcritical set. For parameters along veins, this will still be a compact finite tree; see the discussion in Section \ref{section: laminations}. Furthermore, a parabolic quadratic polynomial still admits a canonically defined Hubbard tree, analogously to PCF parameters.

In other words, by fixing one of the parameters and varying the other along its vein, either there are no ray connections; or the first one is picked up at the tip; or it is picked up a parabolic parameter. Monotonocity of the characteristic lamination along veins (Proposition \ref{proposition: monotonicity along veins}) guarantees that Theorem \ref{theorem: first ray connection} also holds for $c_1$ not postcritically finite, and we get uniform boundedness on the length of any ray connection in the formal mating. Compare these results with Proposition 4.1, 4.4 and Theorem 4.2 in \cite{Jung2017}, which also use forbidden regions of angles to assert boundedness of ray connections, and produce large subsets of $M \times M$ for which the topological mating exists. The following consequence is immediate and known to experts, but important to state nonetheless:

\begin{cor} \label{corollary: existence of topological mating}
Suppose $c_1$ and $c_2$ are postcritically finite parameters for which there are no direct ray connections between the Hubbard trees. Then, for all $c_1'$ and $c_2'$ in the veins of $c_1$ and $c_2$ respectively, the topological mating $f_{c_1'} \coprod f_{c_2'}$ exists, and the generalized Hubbard trees $T_{c_1'}$ and $T_{c_2'}$ are disjointly embedded in the topological $2$-sphere as invariant finite trees.
\end{cor}

Understanding what happens when both parameters are allowed to vary along their veins is more delicate (see Example \ref{example: 1/4 with 3/8} and Question \ref{question: varying both parameters}). Given PCF veins $V_1$ and $V_2$, let $D$ be the set of parameter pairs $(c_1, c_2) \in V_1 \times V_2$ such that there are no ray connections between the generalized trees $T_{c_1}$ and $T_{c_2}$ in the formal mating $f_{c_1} \sqcup f_{c_2}$. We call it the \textit{disjointness locus} for the veins $V_1$ and $V_2$.

\begin{thm} \label{theorem: disjointness locus open}
The disjointness locus $D$ is an open subset of $V_1 \times V_2$.
\end{thm}

Motivated by questions pertaining to disjointness of embeddings in topological and geometric matings, we also describe a subset of $M \times M$ of pairs of  renormalizable parameters for which there are no ray connections between the small Julia sets. Given a rational number $p/q \in S^1$, let $M_{p/q}$ be the small Mandelbrot set corresponding to the $p/q$-sublimb of the main cardoid. Alternatively, $M_{p/q}$ is the set of parameters $c$ that are tunings of the $p/q$-rabbit polynomial; under this interpretation one could denote this set by $r_{p/q} * M$, where $r_{p/q}$ is the parameter of the critically periodic $p/q$-rabbit polynomial.

\begin{thm} \label{theorem: tuned rabbits}
If $c_1 \in M_{p_1/q_1}$ and $c_2 \in M_{p_2/q_2}$ have locally connected Julia sets and are not in conjugate limbs, then there are no ray connections between the small Julia sets in the mating $f_{c_1} \sqcup f_{c_2}$.
\end{thm}

This result answers Question 13 in \cite{Buff2013} for this specific loci of renormalizable parameters in a quadratic mating. More precisely, in this situation, the only possible intersections of boundaries of immediate basins in the geometric mating (when it exists) come from the individual polynomials. Evidently this also implies disjointness of the Hubbard trees.

\subsection{Outline of the paper}

Section \ref{section: background} provides some of the necessary background in complex dynamics. We recall relevant facts about Hubbard trees to define the forbidden region, and Milnor's theory of orbit portraits. It contains a brief introduction to the notion of internal addresses, and descriptions of some well known algorithms for the combinatorics of quadratic polynomials. Specifically, Lavaurs' algorithm for computing QML, a faster way of computing the angle of colanding to a given periodic angle under doubling, the relationship between tuning and external parameter rays, and how to compute the external arguments of the points in the postcritical orbit. The reader familiar with these topics may skip this recollection. 

Section \ref{section: finding forbidden region} details the actual procedure of computing the forbidden region of the Hubbard tree. To do so in an algorithm guaranteed to terminate, we describe bounds on the number of entries of the internal address we have to compute, in the case of a Misiurewicz parameter, whose internal address is infinite. As a small digression, we also explore the periodicity of increments in the internal address for these preperiodic kneading sequences. Proposition \ref{proposition: periodic approximation lands at root} allows us to recover the external arguments of a satellite orbit portrait with valence $v \geq 3$ from just the angled internal address, something not immediate merely from the combinatorics.

In Section \ref{section: hubbard tree automaton} we first give some general definitions and results on automata. Knowing the forbidden region of $H_c$, we construct the Hubbard tree automaton $\mathcal{H}(c)$, show that its regular language coincides with the binary expansions of external arguments of points in $H_c$, and prove several of its properties. We exhibit a collection of examples of Hubbard tree automata in Section \ref{section: examples of hubbard tree automata}.

With the intent of applying these notions to quadratic matings, Section \ref{section: mating} describes the relevant background and the procedure of mating Hubbard tree automata. Section \ref{section: laminations} recovers some standard definitions and results in the theory of quadratic invariant laminations, with the purpose of reproving the known result of monotonicity of characteristic leaves with respect to veins of the Mandelbrot set. We opt to describe veins as topological arcs instead of their classical combinatorial descriptions, which lends itself to a more intuitive visualization of the examples and computations in Section \ref{section: ray connections between hubbard trees}.

Section \ref{section: characteristic states automaton} contains necessary discussions about the nature of characteristic points in ray connections between Hubbard trees, and what this translates to in the language of mating automata. With this we are able to prove Theorem \ref{theorem: first ray connection} in Section \ref{section: algorithm for first ray connection}: the proof describes an algorithm in the mating automaton that finds a ``minimal characteristic arc path'', corresponding to the minimal characteristic landing point of a direct ray connection between the Hubbard trees on $[\alpha_{c_1}, c_1]$. 

In section \ref{section: disjointness locus} we discuss a few examples of what happens when we vary both parameters along their veins in terms of existence of ray connections between the trees, and we prove Theorem \ref{theorem: disjointness locus open}. The proof of Theorem \ref{theorem: tuned rabbits} in Section \ref{section: proof of tuned rabbits} uses only the topology and rotational dynamics of the periodic small Julia sets, and can be understood independently of the previous results on automata. Finally, Section \ref{section: future} points to future directions in understanding ray connections between Hubbard trees.

\subsection{Implementation}
A Python implementation of the algorithm is freely available at \url{https://github.com/evsmath/Hubbard-tree-automata}, with MIT license. It currently contains functions that compute the forbidden region, the Hubbard tree automaton, and the product mating automaton. Development for more functions, such as finding the first periodic ray connection along a vein, is ongoing.

\subsection*{Acknowledgements}

The author wishes to thank Richard Birkett, Caroline Davis, John Hubbard, Wolf Jung, Lasse Rempe, Willie Rush Lim and Dierk Schleicher for their insights into this work, and specially Jeremy Kahn. Rush Lim pointed out that a previous version of Theorem \ref{theorem: tuned rabbits}, stated in terms of Hubbard trees, is valid in cosiderable more generality, applicable to all parameters having locally connected Julia sets and implying disjointness of the small Julia sets. All labeled figures in this paper were obtained in Mandel (\url{https://mndynamics.com/indexp.html}), and the directed graphs were drawn on quiver (\url{https://q.uiver.app/}), and the drawings were made in \url{https://mathcha.io}.

\subsection*{AI disclosure statement}

No artificial intelligence models and software (such as LLMs) were used in the research, writing, reference collecting, and editing of this paper, nor in the coding of the algorithms here described. Minor usage was employed for proofreading and bug catching within the Python implementation. 

\section{Background} \label{section: background}

Throughout this paper, we mainly restrict ourselves to quadratic polynomials. We recall basic facts about their dynamics, and defer to standard texts in the subject for a more thorough exposition \cite{Carleson1993, Douady1984}. Other aspects of the theory are mentioned, sometimes without proof, as needed.

Every complex quadratic polynomial is affinely conjugate to a unique monic polynomial of the form $f_c(z) \coloneqq z^2 + c$, where $c \in \CC$ is arbitrary. The polynomial $f_c$ has a unique finite, simple critical point at $c_0 = 0$, and we label its postcritical orbit as $P_c = \{c_1, c_2, \ldots\}$, where $c_i \coloneqq f^{i}(c_0)$. 

\begin{remark}
For the rest of this paper, when defining an object depending on the polynomial $f_c$, we may use $c$ as the indexing object instead, and may also omit the dependence on $c$ whenever clear from the context.
\end{remark}

The \textit{filled Julia set} $K_c \coloneqq K_{f_c}$ is defined to be set of points $z \in \CC$ such that $|f_c^ n(z)|$ remains bounded under iteration, and the Julia set $J_c$ is the topological boundary of $K_c$. The \textit{Fatou set} $F_c$ is the complement $\CC \setminus J_c$, and there is a unique unbounded component of the Fatou set. The \textit{Mandelbrot set} $M$ is the set of parameters $c \in \CC$ for which the corresponding postcritical set $P_c$ is contained in a compact subset of $\CC$, or equivalently, such that $K_c$ is connected. If $c \in M$, the filled Julia set $K_c$ is compact, full, and there is a unique conformal uniformizing map $\varphi_c : \CC \setminus \overline{\DD} \to \CC \setminus K_c$ tangent to the identity at infinity that conjugates the dynamics of $f_c$ outside of $K_c$ with that of $z^2$ outside of $\overline{\DD}$. This allows us to define \textit{external rays} $\R_c(\theta) \coloneqq \{\varphi_c(re^{i2\pi \theta}): r > 1\}$ of the polynomial $f_c$.

\begin{remark}
Henceforth, the circle $S^1$ is always measured in turns, so that an angle $\phi \in S^1$ varies from $0$ to $1$, which is identified with $0$.
\end{remark}

The polynomial is said to be \textit{postcritically finite} (PCF) if $P$ is finite, that is, the critical point is preperiodic. It is well known that this implies that $K_c$ is locally connected, that every external ray has a well defined landing point $\lim_{r \to 1^+} \varphi_c(r e^{i2\pi \theta})$ on $\partial K_c$ which is continuous with respect to $\theta$, and $K_c$ is arcwise connected. Given $z \in \partial K_c$, its \textit{external arguments} are the angles whose dynamical rays \textit{land on} $z$, that is, the $\theta$ for which $\lim_{r \to 1^+}\varphi_c(re^{2z\pi i \theta}) = z$.

If the critical point is periodic, there is a canonical uniformization of the periodic components of the interior of $K_c$. Since every such component $U$ (which are exactly the bounded Fatou components) is preperiodic, there is a well defined notion of its internal rays. A \textit{regular arc} in $K_c$ is one that traverses interior components only at internal rays. These uniformizations extend homeomorphically to the boundary, and the image of $1 \in \overline{\DD}$ in the closure $\overline{U}$ is called the \textit{root} of the Fatou component. For periodic Fatou components, this is exactly the fixed point on $\partial U$ of the first return map of that component to itself.

\begin{definition} \label{definition: hubbard tree}
The \textit{Hubbard tree} $H_c$ of a postcritically finite quadratic polynomial $f_c$ is the regular hull of the finite postcritical set $P_c$ within $K_c$; That is, it is the smallest set in $K_c$ containing the postcritical orbit $P_c$ and for which any two points are connected by a regular arc. It is a finite topological tree, compact, and forward invariant under the dynamics of $f_c$.
\end{definition}

If $f_c$ is PCF, we will also say that the parameter $c$ is PCF. More generally, $f_c$ is \textit{geometrically finite} if the critical point falls into a repelling cycle, or is attracted to an attracting or parabolic cycle. If the critical point is attracted to an attracting cycle (in particular if the critical point is periodic), then $f_c$ is \textit{hyperbolic}.

The Mandelbrot set is compact, connected, full, and there is a unique uniformizing map $\Phi: \CC \setminus \overline{\DD} \to \CC \setminus M$ that parametrizes the complement of $M$ and is tangent to the identity at infinity. Every hyperbolic quadratic polynomial $f_c$ has a naturally associated component $W$ of the interior of $M$ such that $c \in W$. These are called the \textit{hyperbolic components} of $M$, and they are canonically uniformized by the multiplier $\mu(c)$ of the attracting periodic cycle of $f_c$ for $c \in W$. This uniformization extends homeomorphically to the boundary, and the parameter for which $\mu(c) = 1$ is called the \textit{root} of the hyperbolic component. The parameter for which $\mu(c) = 0$ is called the \textit{center} of $W$, and the corresponding polynomial is poscritically finite with the critical point periodic. A root is contained in the closure of at most two hyperbolic components; if it is contained in the closure of only one, the component that $c$ is the root of is said to be \textit{primitive}, and if not, it is \textit{satellite}. Every parameter for which $\mu(c) = e^{2\pi i p/q}$, where $p/q \in \QQ$, is the root of another hyperbolic component.

The parametrization of the exterior of $M$ allows us to define \textit{parameter rays} $\R_M(\theta) \coloneqq \{\Phi(re^{i2\pi \theta}): r > 1 \}$ of the Mandelbrot set. It is known that if $\theta$ is rational, then $\R_M(\theta)$ lands at a unique parameter in $\partial M$. If $\theta \in S^1$ is periodic under doubling (that is, it has odd denominator), that parameter is the root of a unique hyperbolic component, and it is the landing point of exactly two external rays, which are rational and have the same period under doubling. We denote by $c_{\theta}$ the center of that hyperbolic component, and so $f_{\theta} \coloneqq f_{c_{\theta}}$ is postcritically finite and the critical point is periodic of the same period as $\theta$. If $\theta$ is strictly preperiodic under doubling, then $\R_M(\theta)$ lands at a parameter such that finitely many rays land on it, all of which are preperiodic with same preperiod and period. These are called \textit{Misiurewicz parameters}. We let this parameter be $c_{\theta}$, so that $f_{\theta}$ is postcritically finite, and its critical value $c = c_{\theta}$ is preperiodic of same preperiod as $\theta$, but may have period strictly dividing that of $\theta$.

\begin{figure}[ht]
    \centering
    \begin{subfigure}[b]{0.4\textwidth}
        \centering
        \includegraphics[width=\textwidth]{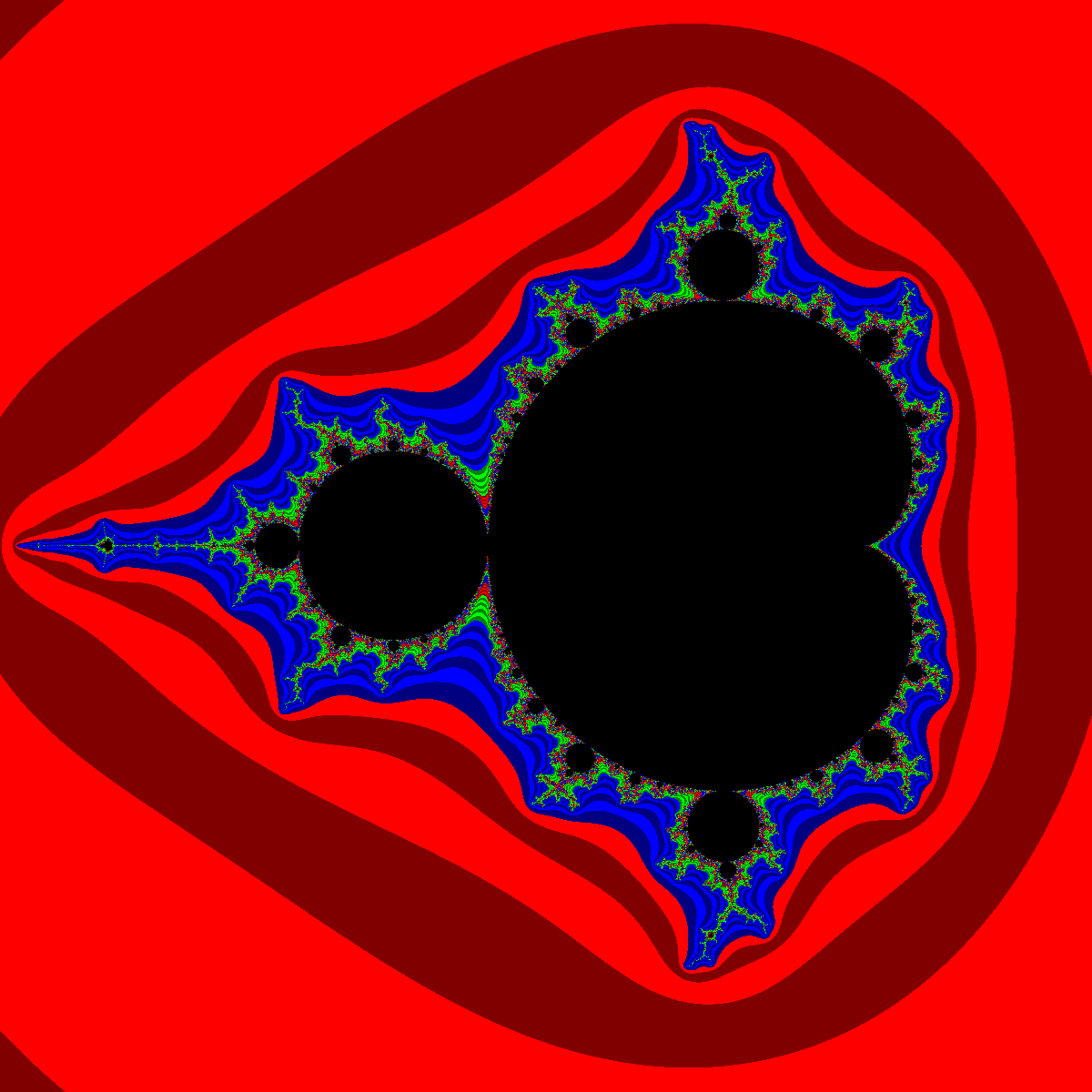}
        \caption{The Mandelbrot set.}
    \end{subfigure}
    \begin{subfigure}[b]{0.4\textwidth}
        \centering
        \includegraphics[width=\textwidth]{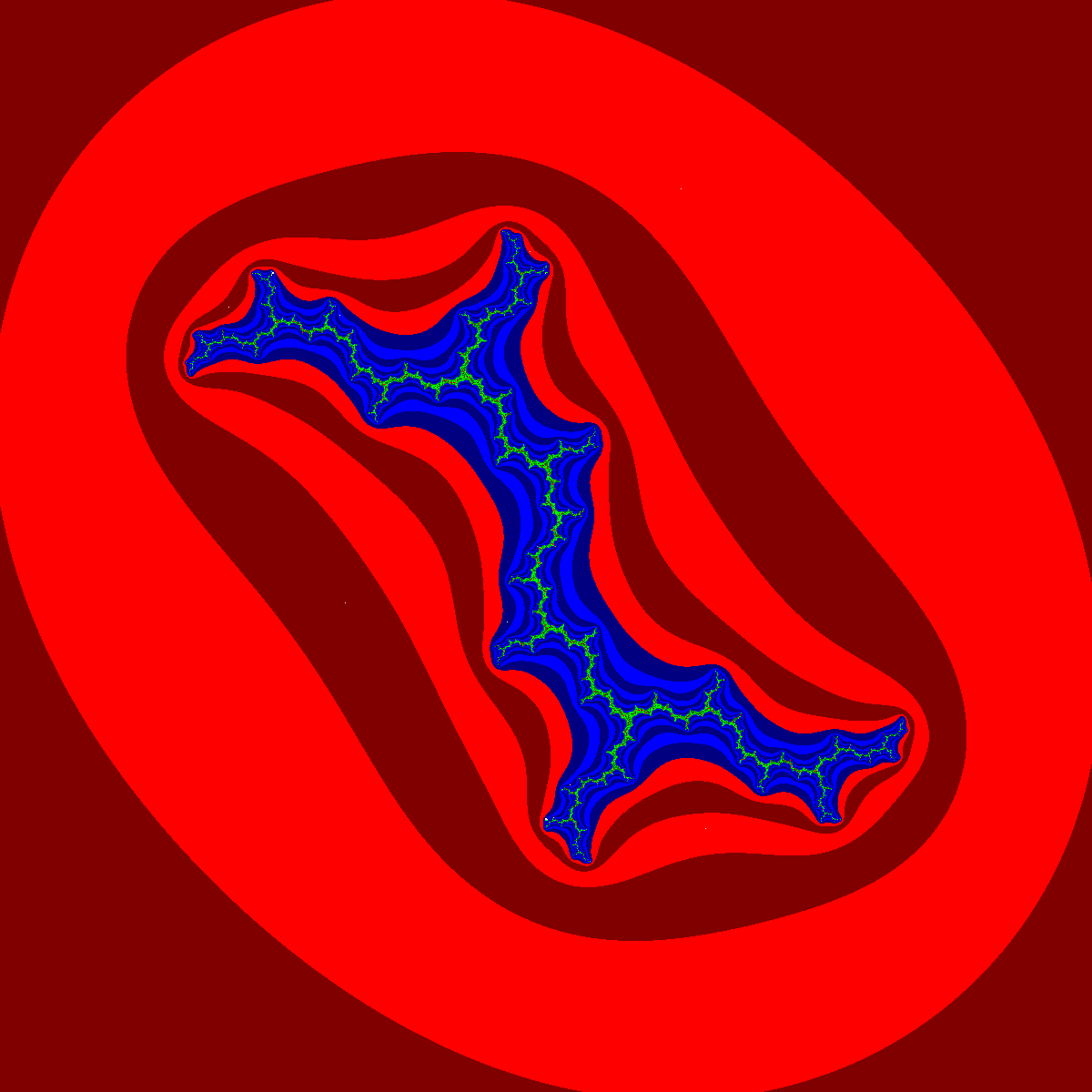}
        \caption{The Julia set for $z^2 + i$.}
    \end{subfigure}
    \caption{Parameter and dynamical planes.}
    \label{parameter dynamical}
\end{figure}

\subsection{Orbit Portraits}

The landing data of rational rays implies many dynamical-parameter correspondences \cite{Lavaurs1986, Schleicher2000}, and one of major use to us is the theory of orbit portraits \cite{Milnor2000}. Proofs of the statements below are deferred to these sources.

Let $f_c$ be be a quadratic polynomial and $z \in J_c$ a periodic point, with orbit $\mathcal{O} = \{z_1, \ldots, z_m\}$ parabolic or repelling. In the language of Milnor, the collection $\mathcal{P} = \{A_1, \ldots, A_m\}$ of the external arguments of each $z_i$ forms an \textit{orbit portrait}. The sets $A_i$ are pairwise unlinked in $S^1$, every angle in all of the $A_i$ is periodic under doubling and has the same period $n$ (the \textit{ray period}), there are exactly $v$ angles in each set $A_i$ (the \textit{valence}), and each $A_i$ is mapped bijectively under doubling to $A_{i+1}$ preserving cyclic order. The portrait is said to be non-trivial if $v \geq 2$ or it is equal to the zero portrait $\mathcal{P} = \{\{0\}\}$, and trivial otherwise.

Non-trivial portraits fall into two possibilities. Either all the angles in each $A_i$ are permutted transitively by the first return map from $A_i$ to itself, and so all angles form a single periodic orbit with some combinatorial rotation number $p/q$ where $q = v$ and $\mathrm{gcd}(p,q) = 1$; or for some collection $A_i$ the rays are not permutted transitively, and so $v = 2$, there are exactly two periodic orbits of angles in the portrait, and each angle is fixed under the first return map. The former is said to be a \textit{satellite} orbit portrait, and the latter a \textit{primitive} orbit portrait.

\begin{prop} \label{proposition: definition of characteristic arc}
Every non-trivial orbit portrait has a unique angle interval $(\theta, \theta')$ of shortest length such that $\theta$ and $\theta'$ are in the same $A_i$. It contains no other angles of the portrait, and it determines the portrait uniquely.
\end{prop}

This arc $(\theta, \theta')$ of the circle is called the \textit{characteristic arc} of the portrait. In the case of the zero portrait, by an abuse of notation we consider the arc $(0,1)$ to be the characteristic one.

\begin{prop} \label{proposition: parameter dynamical correspondence of portraits}
If $(\theta, \theta')$ is the characteristic arc of an orbit portrait, then the parameter rays $\R_M(\theta)$ and $\R_M(\theta')$ land together at a parabolic parameter $c$, which is the root of a hyperbolic component having the same period as that of $\theta$ and $\theta'$. This hyperbolic component is satellite if and only if the portrait is satellite, and primitive if and only if the portrait is primitive. Conversely, every parabolic parameter in $M$ is the landing point of exactly two parameter rays $\R_M(\theta)$ and $\R_M(\theta')$, forming the characteristic arc of an orbit portrait.
\end{prop}

The parameter ray pair, together with the common landing point, cuts out a region of the parameter plane containing that hyperbolic component, called the \textit{wake} $\mathcal{W}_{\mathcal{P}}$ of the portrait.

\begin{prop} \label{proposition: wake dependence}
For every parameter $c$ within the wake of the portrait $\mathcal{P} = \{A_1, \ldots, A_m\}$, there exists a unique periodic cycle $\{z_1, \ldots, z_m\}$ for $f_c$ such that the external arguments of $z_i$ coincides with $A_i$. The periodic points move holomorphically with respect to the parameter $c$ within the wake.
\end{prop}

We say that $f_c$ \textit{realizes} the portrait $\mathcal{P}$ if it has a non-attracting periodic orbit whose external arguments coincide with $\mathcal{P}$. From the above, $f_c$ realizes $\mathcal{P}$ if and only if $c \in \mathcal{W}_{\mathcal{P}}$.

\subsection{Hubbard trees and the forbidden region}

This subsection defines the forbidden region of a Hubbard tree, which characterizes which external rays land on the Hubbard tree or not. The facts below are classical, and included for completeness (See for instance \cite{Li2007}, Section 3.1). Recall from Definition \ref{definition: hubbard tree} that the Hubbard tree $H_c$ of a postcritically finite quadratic polynomial $f_c$ is the regular hull of its postcritical set within the filled Julia set.

For $z \in H$, let $b(z)$ be the number of branches of $H$ at $z$; note that this is in general distinct from the number of components of $K_c \setminus \{z\}$, that is, how many ``branches'' of the filled Julia set are joined at $z$. For $z_1, z_2 \in H$, denote by $[z_1, z_2]$ (respectively $(z_1, z_2)$) the unique closed (open) regular arc in $K_c$ which has $z_1$ and $z_2$ as endpoints. Half-open and half-closed arcs are analogously defined. A point $J_c$ is called \textit{biaccessible} if at least two distinct external parameter rays land on it.

The following is classical \cite[Proposition 4.4]{Douady1984}:

\begin{prop} \label{proposition: hubbard trees}
For all $z \in H$, $b(f(z)) \geq b(z)$, except in the case of $z = c_0$, $b(c_0) = 2$ and $b(f(c_0)) = b(c_1) = 1$.

\begin{itemize}
    \item[(a)] If $f_c$ is hyperbolic and $c_0 = 0$ is periodic of period $n \geq 2$, then there exists an $r$ with $2 \leq r \leq n$ such that, for $1 \leq i \leq r$, $b(c_i) = 1$, and for $r < i$, $b(i) = 2$. The internal angles of the Fatou components corresponding to the branches at the points $c_i$ are $0$, for $b(c_i) = 1$, and $\{0, 1/2\}$, for $b(c_i) = 2$.

    \item[(b)] If $f_c$ is subhyperbolic and $c_0 = 0$ is strictly preperiodic of preperiod $m$ and period $n$, then $b(c_0) = 2$, $b(c_1) = b(c_2) = 1$ and
    $$
    1 = b(c_1) \leq b(c_2) \leq \cdots \leq b(c_m) = b(c_{m+1}) =  \cdots = b(c_{m+n}).
    $$
\end{itemize}
\end{prop}

In the periodic case, note that $b(c_0) = 1$ exactly when $f_c$ is a $p/q$-rabbit, and $2$ otherwise, and in the preperiodic case, $b(c_0)  = 2$ always. This Proposition implies that branch points of $H$ are always in the Julia set.

\begin{lemma} \label{lemma: biaccessible}
If $f_c$ is postcritically finite and $z \in J_c$ is biaccessible, then for some $n \geq 0$ we have $f^n(z) \in H_c$.
\end{lemma}

\begin{proof}
See \cite{Lei1987}, Chapter 3.
\end{proof}

Let $T^0 = H_c$ and $T^1 = f^{-1}(T^0)$ be the preimage of the Hubbard tree. If $P$ is the (finite) forward orbit of the critical point, then $T^{1}$ is the regular hull of $f^{-1}(P)$. Since $T^0$ maps over itself surjectively, we have $T^0 \subseteq T^1$, and the only case where $T^1 = T^0$ is when the Hubbard tree is the entire Julia set, which occurs only for the parameter $c = -2$. In all other cases, we have a disjoint union $T^1 \setminus T^0 = \bigsqcup_{j = 1}^N C_j$ into connected components, where each closure $\overline{C_j}$ is a subtree in $T^1$ whose intersection with $T^0$ is a single point.

\begin{definition}
The \textit{attaching points} of the preimage Hubbard tree $T^1$ are the intersections of the closures $\overline{C_j}$ with the Hubbard tree $T^0$.
\end{definition}

Note that different components may have the same attaching point. Moreover, if $z \in H$ is not the critical point, it is an attaching point if and only if $b(f(z)) > b(z)$, since there will be branches of $T^1 \setminus T^0$ that are then mapped injectively into $T^0$ in the image.

\begin{lemma} \label{lemma: attaching points}
Let $w$ be an attaching point of the preimage Hubbard tree. Then $w$ is preperiodic.

\begin{itemize}
    \item[(a)] If $w$ is not an extremity of the Hubbard tree, then $w$ is strictly preperiodic, $w \in J_{f}$, and at least three external rays land on $w$. Moreover, $w$ is not the critical point.

    \item[(b)] If $w$ is an extremity of the Hubbard tree and $c_0$ is periodic, then $w$ is periodic, and it is the unique attaching point which is an extremity of the Hubbard tree.

    \item[(c)] If $w$ is an extremity of the Hubbard tree and $c_0$ is strictly preperiodic, then $w$ is strictly preperiodic, and it is the unique attaching point which is an extremity of the Hubbard tree.
\end{itemize}
\end{lemma}

\begin{remark}
If the polynomial $f_c$ is hyperbolic, exactly one of the extremities of the Hubbard tree is an attaching point. If $f_c$ is a rabbit, then $w = c_0$, otherwise $w = c_i$ for the unique $c_i$ for which there is a jump in the number of branches:
$$
1 = b(c_i) < b(c_{i+1}) = 2.
$$
Similarly, if $f_c$ is subhyperbolic, the postcritical attaching points are exactly those $c_i$ for which there is a jump in the number of branches to $b(c_{i+1})$.
\end{remark}

\begin{proof}[Proof of Lemma \ref{lemma: attaching points}]
Recall that we have the strict inequality $b(f(w)) > b(w)$ except if $w$ is the critical point $c_0$. In the latter case, it is clear that $w$ is preperiodic. In the former, either $b(w) = 1$, and so $w \in P$ and is preperiodic, or $b(w) \geq 2$, so that $b(f(w)) \geq 3$ and $f(w)$ is a branch point. That branch points of $H$ are preperiodic follows from the No Wandering Triangles Theorem \cite{Thurston2009}.

To prove (a), if $w$ is not an extremity, then either $w = c_0$ or $f(w)$ is a branch point, as above. If $w = c_0$, then $b(w) = 2$ and $b(f(w)) = b(c_1) = 1$, which is not possible if $w$ were to contain extra branches of $T^{1}$. Hence $f(w)$ is a branch point. Since the number of branches strictly increases from $w$ to $f(w)$, the point $w$ cannot be periodic. Moreover, $w$ cannot be in the Fatou set, since no point in a periodic Fatou component can be a branch point because the internal rays at the center corresponding to branches of the Hubbard tree are only $0$ or $\{0, 1/2\}$. For (b), that $w$ is periodic is immediate, and uniqueness follows from the remark. Item (c) is also immediate.
\end{proof}

To each closure $\overline{C_j}$ we associate an interval of angles $I_j$ as follows. If $w$ is the attaching point and $w \in J_c$, then $I_j = (\theta_j^-, \theta_j^+)$, where $\theta_j^-, \theta_j^+$ are the angles in counterclockwise order from $0$ to $1$ adjacent to the branch of $T^1$ at $w$ following into $C_j$. If $w$ is in the Fatou set, so that $w = c_i$ and it is an extremity of the Hubbard tree, then let $w^*$ be the unique periodic point on the boundary of the Fatou component $U_i$ containing $c_i$ that has the same period as $c_i$. We say that $w^*$ is the root of $U_i$, and $ w \in T_0 \cap J_c \cap \partial U_i$. Let $I_j = (\theta_j^-, \theta_j^+)$, where $\theta_j^-$ and $\theta_j^+$ are the angles adjacent to $U_i$ landing at $w^*$. In all cases, these intervals are pairwise disjoint.

\begin{definition} \label{definition: forbidden region}
The union $F = \bigsqcup_{j = 1}^N I_j$ of the intervals constructed above is the \textit{forbidden region} of the Hubbard tree $T^0$.  
\end{definition}

\begin{figure}[ht]
    \centering
    \includegraphics[width=0.7\linewidth]{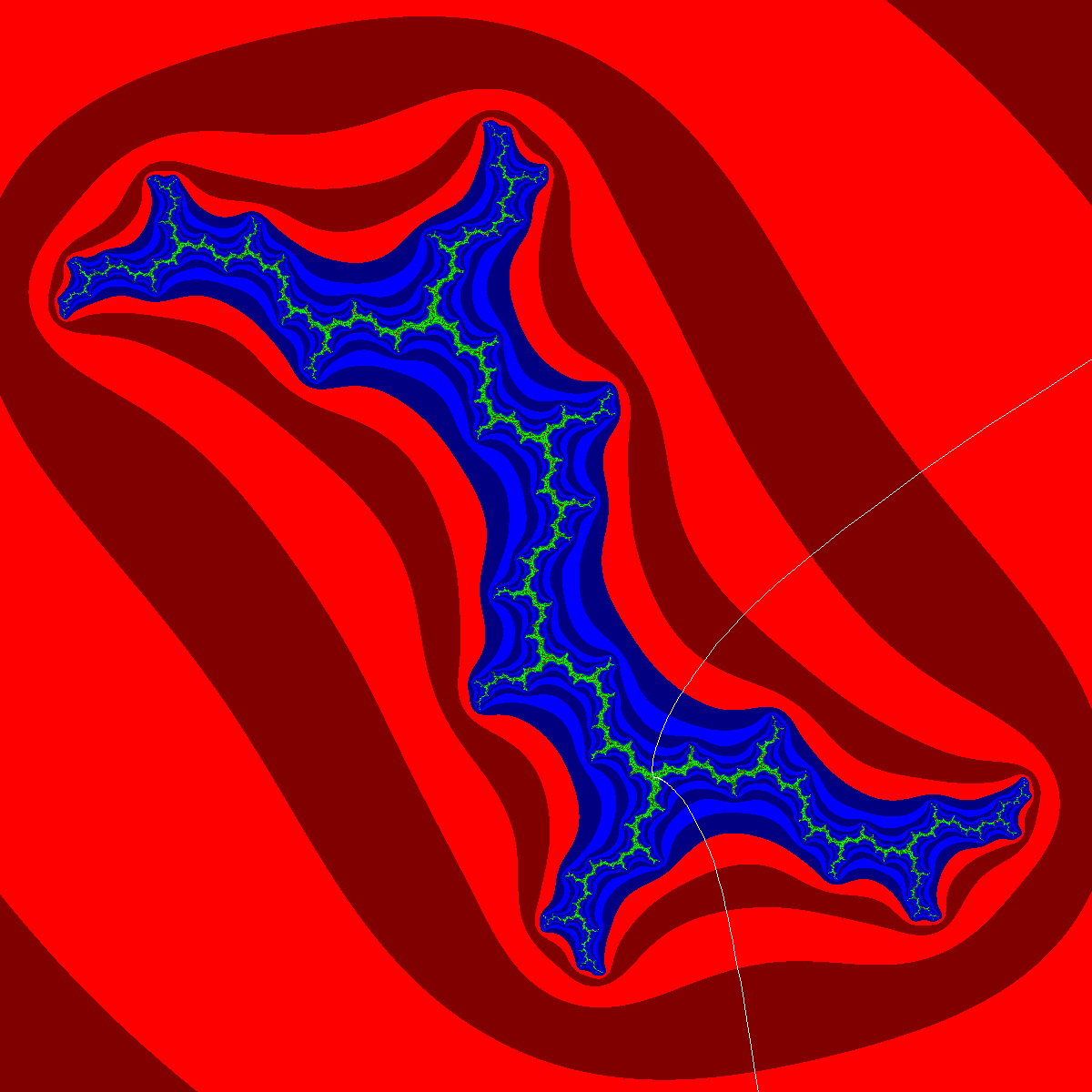}
    \caption{Julia set for $f_{1/6}$ with the external rays $\R_c(1/14)$ and $\R_c(11/14)$ depicted; the forbidden region is the angle interval $(11/14, 1/14) \subset S^1$.}
    \label{fig: 1/6 forbidden region}
\end{figure}

The name is justified by the following result, which appears as Lemma 3.10 in \cite{Li2007}:

\begin{prop} \label{proposition: forbidden region and external arguments}
An external ray at angle $\phi$ lands at a point on the Hubbard tree if and only if the forward orbit of $\phi$ under doubling never enters the forbidden region.
\end{prop}

\begin{proof}
If $\R(\phi)$ lands on $T^0$, for all $k$ the external ray $\R(2^k \phi)$ also lands on $T^0$, and so cannot be in the forbidden region since each forbidden arc encloses components of $T^1 \setminus T^0$. Conversely, suppose that $\R(\phi)$ lands at a point $z \notin T^0$. The regular hull of $P \cup \{z\}$ is a tree that contains $T^0$; let $[z^* , z]$ be the arc in this tree connecting $z$ to $T^0$, so that $z^*$ is in $T^0$. As before, $z^*$ is either an extremity of $T^0$, or a branch point of this larger tree, and therefore preperiodic by the No Wandering Triangles Theorem \cite{Thurston2009}.

\textit{Case 1: $z^*$ is the center of a Fatou component $U$}. Necessarily the arc $(z^*, z]$ contains a unique point $x$ in the boundary of $U$ such that $x \notin T^0$. Either $x$ is an extremity of $K_c$, and so $z = x$ and the only external argument of $x$ is $\phi$, or $x$ is biaccessible, having two external arguments $\alpha^-$ and $\alpha^+$ adjacent to the Fatou component. In the former case, the forward orbit of $z$ will eventually fall into the boundary of the periodic Fatou components, and for some forward image, $z \in \partial U_i$ where $c_i$ is the unique periodic attaching point in the Fatou set. Because $z$ has only one external argument and its orbit is contained in the Julia set, this forward image will also have only one external argument, and so $z$ cannot be the point $w^*$ for the attaching point $w = c_i$. In particular, $\phi \in I_i$. 

In the latter case, $\phi \in [\alpha^-, \alpha^+]$. Because $x$ is biaccessible, it eventually iterates into the Hubbard tree, and at the boundary of a Fatou component. If $k \geq 1$ is the smallest integer such that $f^k(x) \in T^0$, then $f^{k-1}(x) \in T^1 \setminus T^0$, and so the $(k-1)$-th forward image of $[\alpha^- , \alpha^+]$ under doubling is strictly contained in the component of the forbidden region associated to the attaching point $w = c_i$, and so is $\phi$.

\textit{Case 2: $z^*$ is in the Julia set}. If the arc $(z^*, z)$ is contained in the Fatou set, then there is a smallest $k \geq 1$ such that $f^k((z^*, z)) \cap T^0 \neq \emptyset$, since this open arc is inside of a preperiodic Fatou component, whose center eventually falls into the periodic critical cycle. If the arc is not fully contained in the Fatou set, there is a point $z' \in (z^*, z) \cap J_c$, which is biaccessible, and eventually iterates into the Hubbard tree; again we conclude that there exists a smallest $k$ such that $f^k((z^*, z)) \cap T^0 \neq \emptyset$. Therefore $f^{k-1}((z^*, z))$ is an open arc not contained in $T^0$, it contains some point in $T^1 \setminus T^0$, and $f^{k-1}(z^*) \in T^0$, because the forward orbit of $z^*$ is always in the Hubbard tree.

We claim that this implies that $f^{k-1}(z^*)$ is an attaching point of $T^1$. Note that if $y \in f^{k-1}((z^*, z)) \cap ( T^1 \setminus T^0)$, then the arc $(f^{k-1}(z^*), y]$ is contained in $T^1 \setminus T^0$, being the regular hull of these two points in $T^1$, and the closure intersects $T^0$ only at $f^{k-1}(z^*)$. Hence $f^{k-1}(z^*)$ is an attaching point. Moreover, the external arguments $\gamma^-, \gamma^+$ adjacent to the arc $[z^*, z]$ form an interval $(\gamma^-, \gamma^+)$ which contains $\phi$ and maps into the forbidden region, so that $\phi$ eventually maps into the forbidden region.
\end{proof}

\begin{lemma} $0$ is contained in the forbidden region if and only if $\theta$ is not a dyadic rational.
\end{lemma}

\begin{proof}
The $\beta$ fixed point, which is the landing point of the ray $\R_c(0)$, is an extremity of $J_c$, and is in the Hubbard tree if and only if it is a postcritical point, which is equivalent to $\theta$ being dyadic. Therefore if $\theta$ is not dyadic, $0 \in \bigcup_{i \geq 0} f^{-i}(F)$. But as $0$ is fixed under doubling, $0 \in F$.
\end{proof}

\subsection{Internal Addresses} \label{subsection: internal addresses}

Given a parameter $c \in M$, one way to describe its position within the Mandelbrot set is through the notion of internal addresses \cite{Schleicher2016}. This roughly corresponds to looking at the path needed to reach $c$ from the main cardioid, and keeping track of the hyperbolic components of lowest period found along the way. To complement our notations, we say that $P_M(\theta, \theta')$ is a \textit{parameter ray pair} if the parameter rays $\R_M(\theta)$ and $\R_M(\theta')$ land together at a parameter $c \in \partial M$. Recall that if $\theta$ is rational and periodic under doubling, then $\theta$ and $\theta'$ share the same ray period, and the rays land at the root of a hyperbolic component. In this case, we may call it a periodic ray pair.

\begin{definition} \label{internal address}
The \textit{internal address} of a parameter $c \in M$ is the sequence
$$
S_0 \to S_1 \to \cdots \to S_k \to \cdots
$$
defined recursively as follows. We start with $S_0 = 1$ and the ray pair $P_M(0,1)$ (where, by a slight abuse of notation, we treat the angles $0$ and $1$ as distinct but having parameter rays landing together, corresponding to the zero orbit portrait). Given the initial segment $S_0 \to \cdots \to S_k $ and the ray pair $P_M(\theta_k, \theta_k')$, let $P_M(\theta_{k+1}, \theta_{k+1}')$ be the periodic ray pair of lowest period that separates, in the parameter plane, the parameter $c$ from $P_M(\theta_k, \theta_{k}')$, or for which $c \in P_M(\theta_{k+1}, \theta_{k+1}')$. Define $S_{k+1}$ to be the common period of $\theta_{k+1}$ and $\theta_{k+1}'$ under doubling. If there is no such ray pair, then the sequence terminates at $S_k$. If it exists, the next ray pair is unique by Lavaurs' lemma \cite{Lavaurs1986}.
\end{definition}

Note that we do not keep track of the ray pairs themselves, only their periods, or equivalently the periods of the hyperbolic components that the rays land at the roots of.

We can augment the internal address with extra information, which a posteriori uniquely specifies the combinatorial class of the parameter $c$. For each $k$, let $W_k$ be the hyperbolic component corresponding to the entry $S_k$ in the internal address. If $S_k$ is not the last entry, then $c$ is contained in a \textit{$p_k/q_k$-subwake} of $W_k$. That is, if $c'$ is the parameter for which $\mu_{W_k}(c') = e^{2\pi i p_k/q_k}$ and $P_M(\psi, \psi')$ is the periodic parameter ray pair landing on $c'$, then $c$ is in the wake $\mathcal{W}$ determined by that ray pair. The intersection $M \cap \mathcal{W}$ is called the \textit{$p_k/q_k$-sublimb} of $W_k$.

\begin{definition}
The \textit{angled internal address} of a parameter $c \in M$ is the sequence
$$
(S_0)_{p_0/q_0} \to (S_1)_{p_1/q_1} \to \cdots \to (S_k)_{p_k/q_k} \to \cdots,
$$
where $c$ is contained in the $p_k/q_k$-subwake of each hyperbolic component $W_k$. If the internal address of $c$ terminates at $S_k$, then the angled internal address is
$$
(S_0)_{p_0/q_0} \to (S_1)_{p_1/q_1} \to \cdots \to (S_{k-1})_{p_{k-1}/q_{k-1}} \to S_k.
$$
\end{definition}

If $\theta$ is rational and $c_{\theta}$ is the associated postcritically finite parameter, there is an alternative combinatorial way to describe its internal address, through kneading sequences.

For all $\theta \neq 0$ in $S^1$, it has two preimages $\theta/2$ and $(1+\theta)/2$ under doubling. We label the arc between $\theta/2$ and $(1+\theta)/2$ in the counterclockwise order with the symbol $\texttt{1}$, the complementary open arc with $\texttt{0}$, and the common endpoints with $*$. The \textit{kneading sequence} $\nu(\theta)$ of $\theta$ is the itinerary of $\theta$ under doubling with respect to this partition of $S^1$, so that $\nu = \nu_1\nu_2\cdots$ is an infinite sequence in the set of symbols $\{\texttt{1}, \texttt{0}, *\}$, starting with $\texttt{1}$. The internal address of a kneading sequence $\nu$ can be understood as the periods of the ``best periodic approximations'' to the actual sequence. Starting with $S_0 = 1$ and $\nu^ 0 = \overline{\texttt{1}}$, we define $S_{k+1}$ recursively as the position of the first difference between $\nu$ and $\nu^k$, and let $\nu^{k+1}$ be the periodic continuation of the first $S_{k+1}$ entries of $\nu$. If $\nu$ is periodic, then the internal address terminates with final entry $S_k$ the period of $\nu$.

We will mainly restrict ourselves to the case where $\theta$ is rational. If $\theta$ is periodic under doubling, then $\nu(\theta)$ is a periodic sequence having symbol $*$ and same period as $\theta$. In line with the notation in \cite{Schleicher2016}, the \textit{upper kneading sequence} $\mathcal{A}(\nu)$ is the sequence obtained by replacing every instance of $*$ consistently with $\texttt{1}$ or $\texttt{0}$, guaranteeing that the period of $\theta$ appears in the internal address of $\mathcal{A}(\nu)$. Similarly, the \textit{lower kneading sequence} $\overline{\mathcal{A}}(\nu)$ is the one obtained by replacing every instance of $*$ by the other symbol. Geometrically, if $P_M(\theta, \theta')$ is the ray pair landing at the root of $W$, the kneading sequence $\mathcal{A}(\nu)$ corresponds to the initial segment of the kneading sequence of every angle whose external ray lands at $\partial W$, except at the root. The internal address of $c_{\theta}$ (the center of the hyperbolic component that $\R_M(\theta)$ lands at the root of) is therefore the one obtained from the kneading sequence $\mathcal{A}(\nu(\theta))$.

On the other hand, if $\theta$ is strictly preperiodic under doubling, then $\nu(\theta)$ is strictly preperiodic, having same preperiod but period dividing that of $\theta$. The internal address of $c_{\theta}$ is simply the one corresponding to $\nu(\theta)$.

\subsection{The partner angle algorithm}

Lavaurs algorithmically described the landing pattern of periodic rays under doubling on the boundary of the Mandelbrot set through a lamination of the unit disk \cite{Lavaurs1986}. Recall that the corresponding parameter rays land in pairs, defining an involution on this set of angles, where by convention we assume that the rays $\R_M(0)$ and $\R_M(1)$ are distinct and land together at the cusp of the main cardioid.

The description is recursive, producing at each step a lamination $L_k$, meaning a collection of leaves on the closed unit disk $\overline{\DD}$, where two points on the boundary are connected by a leaf if they correspond to a parameter ray pair. Given $L_{k}$, let $A_{k+1}$ be the set of angles of period $k+1$ seen as points on $S^1$, arranged in a linear order assuming $0 < 1$. Starting from the first, we connect it to the next lowest angle in the set $A_{k+1}$ with a leaf, such that the resulting leaf does not cross a leaf of $L_k$. Next we consider the new lowest angle in $A_{k+1}$ not an endpoint of a leaf, and repeat this process. The union of these leaves with $L_k$ will produce $L_{k+1}$, and the union $\bigcup_{k \geq 1} L_k$ represents how angles of the same period are paired up.

\begin{definition} \label{definition: QML}
The closure $\overline{\bigcup_{k \geq 1} L_k}$ is the \textit{quadratic minor lamination}, abbreviated QML.
\end{definition}

The quadratic minor lamiantion was first introduced by Thurston \cite{Thurston2022} in his 1984 manuscript; see Figure II.9 in \cite{Thurston2009}. It is the set of minor leaves of all invariant quadratic laminations, and provides a combinatorial model for the Mandelbrot set (see Section \ref{section: laminations}).

Throughout our algorithm, we will repeatedly have to find \textit{partner angles}. That is, if $\theta$ is rational and periodic under doubling, we will have to find the angle $\theta'$ such that $P_M(\theta, \theta')$ forms a parameter ray pair, or equivalently $\{\theta, \theta'\}$ specifies a leaf in QML. One way to do is is by computing QML up to period $n$, being the period of $\theta$, and check which angle ends up being paired up with $\theta$. This is in general computationally expensive, since we have to compute all laminations of smaller period beforehand, whose number of leaves roughly double in quantity from one period to the next.

A more direct way to do this, linear in the period of the angle $\theta$, is as follows; the author thanks Jeremy Kahn and Lasse Rempe for their explanations of this algorithm, which is also present in \cite{Bruin2011} (pages 168-172, under the notion of ``conjugate angles''). Similarly to the construction of kneading sequences, let $\theta \neq 0$ be rational with odd denominator, and $\theta/2$, $(1+\theta)/2$ its two preimages. On the closed unit disk, connect these two preimages by a diameter, and consider this as a partition of $S^1$. Let $\theta_1 = \theta$ and its orbit be $\{\theta_1, \theta_2, \ldots, \theta_n\}$, where $n$ is the period of $\theta$ and $\theta_n$ is the periodic preimage of $\theta$. The idea is to ``pullback'' the strictly preperiodic preimage of $\theta$ along the orbit of $\theta$ for $n$ steps, and the closest candidate angle will be the partner.

Let $\alpha_n$ be the strictly preperiodic preimage of $\theta$. Suppose we are given $\alpha_k$, for $1 < k \leq n$; we define $\alpha_{k-1}$ to be the preimage of $\alpha_k$ that is on the same side of $S^1$ with respect to the partition. Therefore $\alpha_1$ is a strictly preperiodic angle whose itinerary coincides with that of $\theta$. 

\begin{prop} \label{proposition: partner angle algorithm}
If $\theta'$ is the partner angle to $\theta$, the first $n$ digits in the binary expansion of $\alpha_1$ agree with those of $\theta'$.
\end{prop}

\begin{proof}
Let $J$ be the open arc on $S^1$ whose endpoints are $\alpha_n$ and $2^{n-1}\theta'$. It has length less than $1/2$, and so maps bijectively onto $(\theta, \theta')$. The orbit of $\theta$ never enters $(\theta, \theta')$ by minimality of the characteristic leaf of the portrait, and so it also never enters $J$ and its iterated pullbacks under doubling. Hence the $(n-1)$-th pullback of $J$ consists of angles all sharing the same first $n$ entries in their intineraries.

Moreover, since $0$ is not contained in $(\theta, \theta')$, $0$ and $1/2$ are not contained in $J$ and in any on their pullbacks. This guarantees that the angles also share the same first $n$ digits in their binary expansion.
\end{proof}

With this, we compute $\theta'$ by finding the first $n$ digits in the binary expansion of $\alpha_1$, and repeating them as a periodic block.

\subsection{Angle tuning} \label{subsection: angle tuning}

Each hyperbolic component $W$ of the Mandelbrot set gives rise to a homeomorphic copy $M_W$ of $M$ \cite{Douady1986, Milnor2000}; if $c_0$ is the center of $W$, the image $c' \in M_W$ of $c \in M$ by this homeomorphism is called the \textit{tuning} of $c_0$ by $c$. In particular, every component $W$ has its own $p/q$-sublimbs, where the attracting periodic cycle associated to it undergoes a bifurcation at $\mu(c) = e^{2\pi i p/q}$, and each point of the cycle splits into $q$ points. By the first return map, they are permuted with rotation number $p/q$. These and other facts can be formalized through the notions of renormalization and quadratic-like mappings \cite{Douady1985}.

Of interest to use is how parameter rays landing on $M$ correspond to rays landing on $M_W$:

\begin{thm}[Theorem 2 in\cite{Douady1986}] \label{theorem: angle tuning}
Let $P_M(\theta_-, \theta_+)$ be the parameter ray pair landing at the root of $W$, and $\R_M(\phi)$ a parameter ray landing at $c \in \partial M$. If $c_0$ is the center of $W$ and $c^* = c_0 * c$ is the tuning of $c_0$ by $c$, then $c^*$ has an external argument $\phi'$ associated to $\phi$ as follows. If
$$
\begin{cases}
\theta_{-} & = .\overline{b_1^{0}b_2^{0}\ldots b_k^{0}}, \\
\theta_{+} & = .\overline{b_1^{1}b_2^{1}\ldots b_k^{1}}, \\
\phi & = s_1s_2\ldots s_m \ldots
\end{cases}
$$
are the binary expansions of $\theta_{\pm}$ and $\phi$, then the binary expansion of $\phi'$ is
$$
\phi' = .(b_1^{s_1} b_2^{s_1}\ldots b_k^{s_1})(b_1^{s_2} b_2^{s_2}\ldots b_k^{s_2})\ldots
$$
That is, every instance of $0$ in the binary expansion of $\phi$ is replaced by the binary period of $\theta_-$, and every instance of $1$ is replaced by the binary period of $\theta_{+}$.
\end{thm}

We will employ tuning to find the angles landing at the roots of the $p/q$-sublimbs of a hyperbolic component $W$. For this, we also need to know what these angles are with respect to the main cardioid of $M$, that is, the unique hyperbolic component whose attracting cycle has period $1$. The formula below gives these angles:

\begin{prop}[Theorem 1 in \cite{Douady1986}] \label{proposition: cardioid angles}
Let $t = p/q$ be a fraction in $(0,1)$ in lowest terms, and let $\theta_{\pm} = a_{\pm}/(2^q - 1)$ be the angles of the parameter rays landing at the root of the $p/q$-sublimb of the main cardioid of the Mandelbrot set. Then
\begin{equation} \label{equation: series}
\theta_- = \sum_{0 < s < t} \dfrac{1}{2^{q(s)} - 1} = \sum_{0 < a/b < t}\dfrac{1}{2^b},
\end{equation}
where the first sum ranges over $s$ rational and $q(s)$ is the denominator of $s$ in lowest terms, and the second sum ranges over all rationals $a/b$, possibly not in lowest terms.
\end{prop}

Because we seek a computational implementation of tuning, we convert this infinite series into a finite sum which gives us the desired limit $\theta_-$, knowing that it must be rational with denominator $2^q - 1$. The right-most expression in (\ref{equation: series}) may be rewritten as
\begin{gather*}
\theta_- = \sum_{b = 1}^{\infty} \sum_{0 < a/b < t} \dfrac{1}{2^b} = \sum_{b = 1}^{\infty} \sum_{0 < a < bt}\dfrac{1}{2^b} = \sum_{b = 1}^{\infty} \dfrac{\lceil bt \rceil -1}{2^b},
\end{gather*}
and we estimate the tail end of the series:
$$
\sum_{b > b_0} \dfrac{\lceil b t \rceil - 1}{2^b} < \sum_{b > b_0} \dfrac{b}{2^b},
$$
where, for $x \geq 2$, $f(x) = x/2^x$ is decreasing. Hence, for $b_0 \geq 2$,
$$
\sum_{b > b_0} \dfrac{b}{2^b} \leq \int_{b_0}^{\infty} \dfrac{x}{2^x}dx = \dfrac{1}{(\log 2)2^{b_0}}\left[b_0 - \dfrac{1}{\log 2} \right] < \dfrac{2b_0}{2^{b_0}},
$$
and as soon as the tail of the series is less than $1/(2^q - 1)$, the limit $\theta_-$ is the lowest fraction of the form $a_-/(2^q - 1)$ larger than the respective partial sum. As soon as $b_0 - \log_2 b_0 \geq q+1$, we find
$$
\dfrac{2b_0}{2^{b_0}} < \dfrac{1}{2^{q}} < \dfrac{1}{2^q - 1},
$$
and since $\log_2 b_0 < b_0/2$ for $b_0 \geq 5$, as soon as $b_0 - \frac{1}{2}b_0 > q + 1$ for $b_0 \geq 5$, we find our desired partial sum. Concisely, since $q \geq 2$ always, we get
\begin{equation} \label{equation: theta minus}
\theta_- = \left \lceil(2^q - 1)\sum_{b = 1}^{2q+2}\dfrac{\lceil bt \rceil - 1}{2^b}\right\rceil/ (2^q - 1).
\end{equation}

\subsection{Angles landing at the critical orbit}

Given a rational $\theta$, in order to find the forbidden region of the Hubbard tree $H_{\theta} \coloneqq H_{c_{\theta}}$, we also need the following:

\begin{itemize}
    \item If $\theta$ is periodic under doubling, to find what are all the periodic rays landing at the roots of the periodic Fatou components, and moreover what are the adjacent rays to each component.

    \item If $\theta$ is strictly preperiodic, to find all the rays landing at each point of the postcritical set.
\end{itemize}

The case of $\theta$ periodic is immediate from the partner angle algorithm, which gives us the characteristic leaf of the orbit portrait, and from it all the ray pairs landing at the roots. Assume now that $\theta = a/2^m (2^n - 1)$ is preperiodic of preperiod $m \geq 1$ and period $n$ under doubling, so that its kneading sequence
$$
\nu = \nu_1\ldots \nu_m \overline{\eta_1 \ldots \eta_{l}} = A \overline{B}
$$
is preperiodic of preperiod $m$ and period $l$, where $l$ divides $n$. Explicitly, $n$ is the ray period of $\R_{c_{\theta}}(\theta)$ in the dynamical plane of $c_{\theta}$, and $l$ is the period of the preperiodic point $c_1$ at which the ray lands \cite{Schleicher2000}. This implies that $n/l$ is the number of rays in the orbit of $\theta$ landing at each point of the orbit. If $n/l > 1$, these are all the rays that land on the postcritical set $P$, and $c_1$ eventually iterates into a satellite orbit portrait. We know how these rays are grouped up in terms of landing since the $n$ periodic angles in the orbit are partitioned into groups of $n/l$, and one may pullback these collections along the orbit of $\theta$ to find which rays land together at each $c_i$. 

The pullback works as follows: if $A = \{\phi_1, \ldots, \phi_{k}\}$ is any collection of all the external arguments of a point $z \in J_{\theta}$ which is not the critical value, and $\{\theta/2, (1+\theta)/2\}$ are the two preimages of $\theta$, the $2n/l$ preimages of the angles in $A$ are divided into two sets of $n/l$ angles based on the partition of $S^1$ with respect to the diameter connecting the preimages of $\theta$. Dynamically, this partition represents the common landing point $c_0 = 0$ of the dynamical ray pair $\{\R(\theta/2), \R((1+\theta)/2)\}$ in the preperiodic case, or the critical Fatou component $U_0$ in the periodic case. Later, we will use this pullback procedure for branch points of the Hubbard tree as well.

If $n/l= 1$, there are two possibilities; either only one ray lands at each postcritical point, therefore forming a trivial portrait, and so all points are extremities of the Hubbard tree $H_{\theta}$, or exactly two rays of different ray orbits and same ray period land at each point of the orbit, and therefore the angles correspond to a primitive orbit portrait.

If $\theta_1 = \theta$, let $\{\theta_{m+1}, \ldots, \theta_{m+n}\}$ be the periodic orbit of $\theta$. Each of these angles has a corresponding partner angle that lands with it in the Mandelbrot set, say $\{\theta_{m+1}', \ldots, \theta_{m+n}'\}$, computed by the partner angle algorithm. Each ray pair $(\theta_j, \theta'_j)$ lands at the root of a satellite or primitive hyperbolic component, that is, they form the characteristic leaf of a satellite or primitive portrait, according to whether they are or aren't in the same orbit under doubling, respectively. All of the corresponding angle sectors are either equal, nested or disjoint. The satellite portraits are not realized if $n/l = 1$, and the primitive ones that are realized are exactly those for which $\theta$ lies in the interior of the corresponding angle sector. This allows us to decide if the portrait is trivial or primitive, and in the latter case, what are the rays landing at the points of $P$. In fact, because the characteristic leaf of an orbit portrait determines it uniquely (Proposition \ref{proposition: definition of characteristic arc}), we have:

\begin{lemma} \label{lemma: uniqueness of non-trivial portrait}
For any preperiodic $\theta$, there exists at most one ray pair $(\theta_j, \theta_j')$ for which $\theta \in (\theta_j, \theta_j')$.
\end{lemma}

Therefore, by iterated pullback along the postcritical orbit, we obtain the set of all external arguments of each $c_i \in P$, in all cases.

\section{Finding the Forbidden Region} \label{section: finding forbidden region}

Let $\theta \in \QQ$ be a rational angle, and $c_{\theta}$ the corresponding postcritically finite parameter in $M$ (either a Misiurewicz parameter or the center of a hyperbolic component). Our goal is to algorithmically find the forbidden region $F_{\theta}$ for the Hubbard tree $H_{\theta}$. As outlined in the introduction, the first steps of the algorithm are:

\begin{enumerate}[1)]
    \item Compute the angled internal address of $c_{\theta}$, up to a point where all entries $(S_k)_{p_k/q_k}$ with $q_k \geq 3$ are found.

    \item Find all periodic branch points of $H_{\theta}$ and the rays landing on them.

    \item Find the attaching points of the preimage Hubbard tree $T^1 \supseteq T^0$ and the forbidden angle arcs. 
\end{enumerate}

In general, the angled internal address of a parameter $c \in M$ may contain infinitely many entries for which $q_k \geq 3$. However, for postcritically finite $c = c_{\theta}$, there are only finitely many such entries:

\begin{lemma}
If $f_c$ is poscritically finite, it realizes only finitely many orbit portraits with valence $v \geq 3$.
\end{lemma}

\begin{proof}
Every periodic point $z$ in an orbit portrait of valence $v \geq 3$ is in the Hubbard tree $H_c$, by Lemma \ref{lemma: biaccessible}. Every pair of adjacent external arguments of $z$ determines a branch of $H_c$ at $z$, since the first return map is locally a homeomorphism, the arguments are permuted transitively by it, and the Hubbard tree is forward invariant. Therefore $z$ is a branch point, and $H_c$ has finitely many branch points.
\end{proof}

Computation of the angled internal address in step 1) is done algorithmically through the kneading sequence $\nu(\theta)$ of $\theta$. If $\theta$ is preperiodic, the internal address is infinite, but there are explicit bounds on how much of the internal address we need to compute in order to find all the finitely many entries with $q_k \geq 3$. Step 2) follows from step 1) once we compute the associated ray pairs $(\theta_k, \theta_k')$, through the procedure detailed in Subsection \ref{subsection: satellite orbit portraits}. This is due to the following facts, the first of which is Proposition 2.6 in \cite{Schleicher2016}, and the second is a consequence of Lemma \ref{lemma: biaccessible} and the results in \cite{Milnor2000}.

\begin{prop} \label{proposition: denominator not 1/2}
If $c \in M$ is contained in the subwake of a hyperbolic component $W$ at internal angle other than $1/2$, then $W$ occurs in the internal address of $c$. 
\end{prop}

\begin{prop}
Let $z$ be a periodic branch point of $H_c$, $\mathcal{P}$ be the portrait corresponding to its orbit, and $(\theta, \theta')$ be its characteristic arc. Then the parameter ray pair $P_M(\theta, \theta')$ lands at the root of a satellite hyperbolic component $W_{p/q}$ which is on the $p/q$-sublimb of a hyperbolic component $W$, where $q = b(z)$.
\end{prop}

In summary, every orbit of periodic branch points of $H_c$ corresponds to a unique entry $(S_k)_{p_k/q_k}$ with $q_k \geq 3$ in the angled internal address, and vice-versa.

If $\theta$ is periodic under doubling, step 1) is just computing the internal address from the kneading sequence. However, $\nu(\theta)$ is $*$-periodic, and we have to figure out what is the upper kneading sequence $\mathcal{A}(\nu(\theta))$, that is, whether we have to subsitute $*$ for a $0$ or a $1$, before computing the internal address. A simple way to do this is just try each substitution, and check which one contains the period of $\theta$ as an entry in the corresponding internal address. Another way to do so is via tuning:

\begin{lemma} \label{lemma: upper kneading sequence}
Let $\theta$ be periodic of period $n$ under doubling, and $W$ be the hyperbolic component that $\R_M(\theta)$ lands at the root of. Then $\mathcal{A}(\nu(\theta))$ is the periodic kneading of period $n' \leq n$ whose first $n$ entries coincide with the $*$-periodic kneading sequence of $\nu(\gamma)$, where $\gamma$ lands at the root of $W_{1/2}$.  
\end{lemma}

The proof of this fact is immediate, as there are no intermediate components of period diving $n$ between $W$ and $W_{1/2}$. The tuning and kneading algorithms allow us to compute $\gamma$, and therefore $\nu(\gamma)$.

If, however, $\theta$ is strictly preperiod under doubling, the internal address is infinite. As mentioned before, we only need to calculate the address up until we are sure that we have found all entries $(S_k)_{p_k/q_k}$ such that $q_k \geq 3$, as this allows us to find all periodic branch points of $H_{\theta}$.

\subsection{Preperiodic internal addresses} \label{subsection: preperiodic internal addresses}

Assume that $\theta$ is preperiodic under doubling, with $\nu(\theta)$ having preperiod $m$ and period $l$. The Hubbard tree $H$ of $f_{\theta}$ has at most $m + l$ extremities, and a straightforward computation asserts that if $b$ is the number of branch points, then $b \leq m + l - 2$. By realizability of satellite portraits \cite{Milnor2000}, for every entry $(S_k)_{p_k/q_k}$ of the angled internal address of $\nu(\theta)$ there is a unique non-attracting periodic orbit of period $S_k$ and valence $q_k$ on the Hubbard tree, realized within the $p_k/q_k$-subwake of the hyperbolic component associated to $S_k$.

If $q_k \geq 3$, then that periodic orbit consists of branch points of the Hubbard tree, and so
\begin{equation} \label{equation: branch bounds}
S_k \leq b \leq m + l - 2.
\end{equation}
This guarantees that every entry $S_k$ in the internal address greater than $m + l - 2$ necessarily has denominator $q_k = 2$, and so no more satellite periodic orbit portraits will be found for that entry or later ones. In other words, to find all satellite orbit portraits with valence $v \geq 3$, we need only compute the internal address up to $S_k \leq m + l - 2$, giving us an explicit upper bound. It is in general not sharp, since every valence $\geq 3$ orbit portrait found before that entry in the internal address already accounts for some of the branch points of the Hubbard tree.

As a digression, we further analyse the behavior of internal addresses coming from preperiodic angles. Because the kneading sequences are preperiodic, there is essentially only a finite amount of information that it provides to compute the internal address. In fact, we will find that it eventually stabilizes into a periodic sequence of increments from one entry to the next, determined by the ``offsets'' of the approximating periodic sequences.

\begin{lemma} \label{lemma: preperiodic kneading bounds}
Let $\nu(\theta) = A\overline{B}$ be a preperiodic kneading sequence of preperiod $m$ and period $l$, and $S_0 \to S_1 \to \cdots \to S_k  \to \cdots$ the corresponding internal address. 
\begin{itemize}
    \item[(a)] If $S_k < m$, then $S_{k+1} \leq m + \mathrm{lcm}(S_k, l) < m(l+1)$.

    \item[(b)] If $S_k \geq m$, then $S_{k+1} - S_k \leq m + l$.
\end{itemize}
\end{lemma}

\begin{proof}
(a) If $S_k < m$, then $\nu^k = \overline{A'}$, where $A'$ is some initial segment of the preperiod $A$ and $\nu^k$ is the $k$-th approximating periodic kneading sequence of $\nu = \nu(\theta)$ as in Subsection \ref{subsection: internal addresses}. If $S_{k+1} > m + \mathrm{lcm}(S_k, l)$, that is, $\nu^k$ and $\nu$ agree up to position $m + \mathrm{lcm}(S_k, l)$, there is some cyclic permutation $A^*$ of the block $A'$ such that the blocks $(A^*)^{\mathrm{lcm}(S_k, l)/S_k}$ and $(B)^{\mathrm{lcm}(S_k, l)/l}$ are equal, implying that $\nu^k = \nu$. This in turn implies that the internal address terminates at $S_k$, but this only occurs for $\nu$ periodic, whereas $\nu$ is strictly preperiodic.

(b) If $S_k \geq m$, then the approximating periodic kneading sequence $\nu^k$ is of the form
\begin{equation} \label{periodic nu}
\nu^k = \overline{AB\dots BB'} = \overline{A(B)^{u}B' }, 
\end{equation}
where $u \geq 0$ and $B'$ is a (possibly empty) strict initial segment of the block $B$ (that is, $B' \neq B)$. Suppose that $B''$ is the complementary block, that is, $B = B' B'' $, and consider the cyclic permutation $B^* = B''B'$ of $B$. Then, comparing $\nu^k$ with $\nu$, we have
$$
\begin{array}{lcl}
    \nu & = & ABB\ldots BB'|B^*B^*\ldots B^{**}B^{**}B^{**} \ldots\\
    \nu^k & = & ABB\ldots BB'|\ldots A\ldots B \ \ B \ \  B \ \ \ldots
\end{array}
$$
where $B^{**}$ is some other cyclic shift of $B$. If the $\nu$ and $\nu^k$ agreed up to position $S_k + m + l$, then $B^{**} = B$, and because the block $B$ is not periodic of period less than $l$, necessarily $A$ is of the form $(B^*)^v$ for some $v \geq 1$. But this also implies that $\nu = \nu^k$, which is not possible for $\nu$ strictly preperiodic.
\end{proof}

For each $k$ such that $S_k \geq m$, the periodic block for $\nu^k$ will terminate in some initial segment $B' _k$ of $B$. From the comparison between $\nu$ and $\nu^k$, we see that this initial segment determines the shifted block $B^*_k$, and if $S_k \geq m + l$, also the shifted block $B^{**}_k$. This in turn determines the increment $S_{k+1} - S_k$, without necessarily knowing $k$ or $S_k$. In fact, this implies that $B'_k$ determines the next initial segment $B'_{k+1}$ of $B$ that is at the end of the period of $\nu^{k+1}$. Therefore the sequence of initial segments is eventually periodic and cycles, determining the consecutive differences $S_{k+1} - S_k$.

Lemma 2.2 in \cite{Schleicher2016} gives the following formula for the denominators in the angled internal address:
\begin{equation} \label{equation: denominators}
q_k = \begin{cases}
\dfrac{S_{k+1} - r}{S_k} + 1, & \text{ if } S_k \in \mathrm{orb}_{\rho_{\nu}}(r), \\
\dfrac{S_{k+1} - r}{S_k} + 2, & \text{ if } S_k \notin \mathrm{orb}_{\rho_{\nu}}(r),
\end{cases}
\end{equation}
where $r \in \{1, \ldots, S_k\}$ is congruent to $S_{k+1}$ modulo $S_k$, $\rho_{\nu}(r) \coloneqq \min\{ k > r : \nu_{k} \neq \nu_{k-r}\}$, and the $\rho_{\nu}$-orbit of $r$ is denoted by $\mathrm{orb}_{\rho_{\nu}}(r)$. This orbit can be understood as follows: we start with the initial segment of $\nu$ of length $r$ and repeat it to form a periodic sequence $\eta^0$. The position of the first difference between $\eta^0$ and $\nu$ is equal to $\rho_{\nu}(r)$, and we now take the initial segment of $\nu$ of length $\rho_{\nu}(r)$ to be the next periodic block that forms $\eta^1$. Iterating this construction, we compute the $\rho_{\nu}$-orbit or $r$. Note that the internal address is, in this language, just the $\rho_{\nu}$-orbit of $r = 1$. Hence the first formula for $q_k$ applies if the $\rho_{\nu}$-orbits if $1$ and $r$ ``merge'' before or at the entry $S_k$, and the second formula applies if not.

By (\ref{equation: branch bounds}), we have $q_k = 2$. Then, in formula (\ref{equation: denominators}), necessarily $(S_{k+1} - r)/S_k = 1$, that is, $r = S_{k+1} - S_k$, so that the increment from one entry of the internal address to the next is also determined by just the initial segment $B'_k$. Note however that different initial segments may determine the same increment.

Let $\{I_1, \ldots, I_t\}$ be the periodic sequence of initial segments of $B$ obtained after the entries of the internal address are larger than or equal to $m+l$, so that if $I_j = B'_k$ is the initial segment of $B$ corresponding to $\nu^k$ as in (\ref{periodic nu}), then $I_{j+1}$ is the corresponding initial segment of $B$ for $\nu^{k+1}$, modulo $t$. Also let $r_j$ be the increments determined by the initial segments $I_j$. From the formula for $q_k = 2$, we necessarily have $S_k \in \mathrm{orb}_{\rho_{\nu}}(r_j)$, showing that the $\rho_{\nu}$-orbits of $r_j$ and $1$ have all merged prior to or at the entry $S_k$. Moreover, this concludes that the sequence of increments $r_k = S_{k+1} - S_k$ follows a periodic pattern, determined by the period of the initial segments $B'_k$.

\subsection{Satellite orbit portraits} \label{subsection: satellite orbit portraits}

We now describe how to convert the information of an entry $(S_k)_{p_k/q_k}$ with $q_k \geq 3$ in the angled internal address into the corresponding satellite orbit portrait, by computing the corresponding angle pair $(\theta_k, \theta_k')$ separating $c$ from the main cardioid. At first, the solution seems to be through the quadratic minor lamination, that is, computing it to a sufficiently high period and relating its nesting structure for $\theta$ with respect to the known angled internal address. In order to avoid extraneous calculations, we could even compute the QML only partially, ``near'' $\theta$. Fortunately, in the case of satellite bifurcations with $q_k \geq 3$, there is a simpler way of figuring out what the corresponding angle pair is:

\begin{prop} \label{proposition: periodic approximation lands at root}
Let $\theta$ be an angle having binary expansion $0.b_1b_2b_3\ldots$ and $c_{\theta} \in M$ any parameter in the impression of $\theta$, having angled internal address $(S_0)_{p_0/q_0} \to \cdots \to (S_k)_{p_k/q_k} \to \cdots$. Suppose that $k \geq 0$ is such that $q_k \geq 3$. If $\varphi_k = .\overline{b_1b_2\ldots b_{S_k}}$, then $\varphi_k$ lands at the root of the hyperbolic component corresponding to the entry $S_k$ in the internal address of $c_{\theta}$.
\end{prop}

In other words, one of the angles in the angle pair $(\theta_k, \theta_k')$ is the ``naive'' approximation of $\theta$ by an $S_k$-periodic angle. We say naive because $\varphi_k$ may not be the $S_k$-periodic angle that is indeed closest to $\theta$, but a posteriori it will be the closest to $\theta$ from one side.

We state a couple of preliminary facts about widths of hyperbolic components \cite{Schleicher2016}.

\begin{lemma} \label{lemma: size of hyperbolic components}
For a given hyperbolic component $W$ of $M$ such that the parameter ray pair $(\theta, \theta')$ lands at its root, let $|W|$ be its width $|\theta - \theta'|$. Then:

\begin{itemize}
    \item If $W_{p/q}$ is the hyperbolic component corresponding to the $p/q$-bifurcation from $W$, then
    $$
    |W_{p/q}| = |W| \dfrac{(2^n - 1)^2}{2^{qn} - 1}.
    $$

    \item If $W$ has period $n$, then
    $$
    \dfrac{|W| - |W_{1/2}|}{2} < |W| \dfrac{1}{2^n}.
    $$
\end{itemize}
\end{lemma}

\begin{proof}[Proof of \ref{proposition: periodic approximation lands at root}]
Let $W$ be the hyperbolic component corresponding to the entry $S_k$, $(\theta_k, \theta'_k)$ be the parameter ray pair landing at the root of $W$, and $(\alpha_k, \alpha'_k)$ be the parameter ray pair landing at the root of $W_{1/2}$. Either $\theta \in (\theta_k, \alpha_k)$, or $\theta \in (\alpha'_k, \theta'_k)$; assume the former without loss of generality.

The statement is clear if $k = 0$; if $k > 0$, then $|W| < 1/3$. Since $|\theta_k - \alpha_k| < |W|/2^{S_k} \leq 1/(3 \cdot 2^{S_k})$ by Lemma \ref{lemma: size of hyperbolic components}, there can be no other parameter rays of period dividing $S_k$ in the interval $(\theta_k, \alpha_k)$, and in particular $\theta_k$ is the closest angle of period dividing $S_k$ under doubling to $\theta$ from the left (here, by abuse of notation, we may assume the circular order as $0 < 1$).

We show that $\varphi_k = \theta_k$. Note that $|\theta - \varphi_k| \leq .0\ldots 0\overline{1} = 1/2^{S_k}$, and that if $b_1 = b_{S_{k}+1}$, then this difference is less than $1/2^{S_k + 1}$. If this is the case, then $\theta_k = \varphi_k$, because the distance of $\theta$ to any angle of period dividing $S_k$ to its right would be greater than $2/(3 \cdot 2^{S_k})$.

Consider the case where $W$ is not a real hyperbolic component. We then claim that $b_1 = b_{S_k + 1}$; this is because the angle $\theta$, within the angle sector of $W$, has its parameter ray in the first renormalization locus of $W$, since it is not in the $1/2$-sublimb $W_{1/2}$. Hence $2^{S_k}\theta$ is in $(\theta_k, \theta'_k)$, and therefore shares the same bit value as $\theta_k$, and as $\theta$ itself. Now suppose that $W$ is a real hyperbolic component. If $b_1 \neq b_{S_{k}+1}$, then because we are assuming that $\theta \in (\theta_k, \alpha_k)$, it is of the form
$$
\theta = .0b_2 \ldots b_{S_k}1 b_{S_k + 2}\ldots,
$$
and
$$
\varphi_k = .0b_2\ldots b_{S_k}0\ldots 
$$
But this implies that $\varphi_k < \theta$, and so necessarily $\varphi_k = \theta_k$. The other case, where $\theta \in (\alpha'_k, \theta'_k)$, is analogous.
\end{proof}

Knowing the parameter ray pairs $(\theta_k, \theta_k')$, and through the tuning algorithms described previously, we find the ray pair landing at the root of the hyperbolic component $W_{p_k/q_k}$. This ray pair forms the characteristic leaf of the satellite orbit portrait, and so we complete step 2).

\subsection{Finding the attaching points} \label{subsection: finding the attaching points}

Knowing all the satellite orbit portraits realized by $f_{\theta}$ with valence $v \geq 3$, or equivalently, the collections of rays landing at each periodic branch point of the Hubbard tree, we use this information to find the attaching points of the preiamge Hubbard tree $T^1$, and subsequently find the forbidden region $F$.

We separately consider the cases of extremities of $H_{\theta}$ being attaching points. Recall that if $\theta$ is periodic, exactly one extremity $c_i$ is an attaching point, and it is the unique postcritical point $c_i$ such that $1 = \nu(c_i) < \nu(c_{i+1}) = 2$, or $c_0$ in the case of a $p/q$-rabbit; and if $\theta$ is preperiodic, every $c_i$ such that $\nu(c_i) < \nu(c_{i+1})$ is an attaching point.

For $\theta$ preperiodic, let $A_i = \{\phi^i_1, \ldots, \phi^i_{t}\}$ be the collection of external arguments of $c_i$ for $i \geq 1$, ordered counterclockwise. We may check if each consecutive angle interval $(\phi^i_j, \phi^ i_{j+1})$ contains an element of the orbit of $\theta$, excluding those that are external arguments of $c_i$. If it does, then to this angle interval there corresponds a branch of the Hubbard tree at $c_i$, since the corresponding dynamical ray pair $\{\R_{c_{\theta}}(\phi^i_j), \R_{c_{\theta}}(\phi^i_{j+1})\}$ separates the postcritical set in the dynamical plane. Therefore if $(\phi^i_j, \phi^i_{j+1})$ does not contain elements of the orbit of $\theta$ but the image interval under doubling does, excluding the external arguments of $c_{i+1}$, the interval $(\phi^i_j, \phi^ i_{j+1})$ corresponds to a branch of the preimage Hubbard tree attached at $c_i$, and therefore an open arc of the forbidden region. A similar check is done in the case of $\theta$ periodic, but now the angles in the orbit of $\theta$ land at the roots of the Fatou components $U_i$, and separation of the postcritial set by angle intervals is checked analogously.

As for attaching points that are not extremities, we need only pullback the periodic branch points of the Hubbard tree, found on step 2). Let $z$ be a periodic, non-postcritical branch point of $H$, and $C_z = \{\phi_1, \ldots, \phi_t\}$ be the collection of its external arguments. We find the external arguments $C_w$ of the strictly preperiodic preimage $w$ of $z$ (the periodic preimage is just another periodic branch point, and therefore another ``seed'' for this pullback algorithm), check which angle intervals $(\phi_j, \phi_{j+1})$ are such that their preimages stop separating the postcritical set, and add those to the list of open arcs of the forbidden region.

There are three possibilities:

\begin{itemize}
    \item If at least three consecutive angle intervals for the preimage $w$ contain angles in the orbit of $\theta$, then $w$ is still a branch point, and we continue this procedure with its two preperiodic preimages. 

    \item If exactly two angle intervals of $w$ contain angles in the orbit of $\theta$, then $w$ sits inside of an edge of the $H$, and we need not continue the algorithm for $w$.

    \item If exactly one angle interval contains angles in the orbit of $\theta$, then either this preimage is an extremity of the Hubbard tree, which cannot be the case since we have considered postcritical periodic points beforehand, or not in the Hubbard tree at all. In that case, $w$ is contained in a component of $T^1 \setminus T^0$, and so all of its angle intervals will be already be contained in a larger angle interval of the forbidden region; in other words, it cannot be an attaching point.
\end{itemize}

As there are only finitely many branch points of the Hubbard tree (and in fact, we have the upper bound in (\ref{equation: branch bounds})), the algorithm eventually terminates, with an upper bound on the number of times a pullback gives a branch point.

\section{Automata and The Hubbard Tree Automaton} \label{section: hubbard tree automaton}

\subsection{General automata}

We define automata as convenient for our purposes. A more general treatment of these objects is found in \cite{Hopcroft2001}, containing other classical definitions and properties of those objects.

\begin{definition} \label{definition: automaton}
We define an \textit{($\omega$-)automaton}, or \textit{finite state state machine}, as a pair $\mathcal{A} = (\mathcal{G}, \chi)$, where:
\begin{itemize}
    \item $\mathcal{G} = (\mathcal{S}, \mathcal{E})$ is a finite directed graph, with set of vertices $\mathcal{S}$ and edges $\mathcal{E}$;
    \item $\chi : \mathcal{S} \to \Sigma$ is a function on the set of vertices with values in a finite set $\Sigma$.
\end{itemize}
The vertices of $\mathcal{G}$ are called the \textit{states} of $\mathcal{A}$, and the directed edges are the \textit{transitions}. If $A \to B$ is a directed edge, then it is a transition \textit{from} $A$ and \textit{to} $B$.
\end{definition}

We will assume that, for all automata under consideration, there is a unique transition from a state to another, and we emphasize that loops (a transition from a state to itself) are allowed. In what follows, ``infinite'' will always mean ``one-sided infinite''.

\begin{definition} \label{definition: regular language}
A \textit{word} in the alphabet $\mathcal{S}$ is a finite or infinite sequence of states. A \textit{path} in $\mathcal{A}$ is a word $w = w_0w_1w_2\ldots$, finite or infinite, such that there exists a transition from $w_i$ to $w_{i+1}$ for all $i \geq 0$. The function $\chi$ extends to sets of words and paths, with values in $\Sigma^{\leq \omega}$. An \textit{allowable ($\Sigma$-)sequence} in $\Sigma^{\NN}$ is the image under $\chi$ of an infinite path in $\mathcal{A}$. The \textit{regular language} $\mathcal{L} = \mathcal{L}(\mathcal{A})$ is the colection of all allowable sequences for $\mathcal{A}$. We say that $\mathcal{A}$ \textit{recognizes} a sequence $s \in \Sigma^{\NN}$ if it is in the regular language.
\end{definition}

There is a natural shift funtion $\sigma$ defined on infinite words, infinite paths and infinite sequences on $\Sigma$, which forgets the first position. By a slight abuse of notation, we will denote all of these functions with the same symbol.

\begin{prop} \label{proposition: regular language is shift invariant and compact}
Suppose that $\Sigma^{\NN}$ is equipped with the Tychonoff topology. Then $\mathcal{L}(\mathcal{A})$ is $\sigma$-invariant and compact.
\end{prop}

\begin{proof}
If $s \in \Sigma^{\NN}$ is recognized by $\mathcal{A}$, then $s = \chi(w)$ for some infinite path $w$ in $\mathcal{A}$. But then $\sigma(s) = \chi(\sigma(w))$, where $\sigma(w)$ is also an infinite path in $\mathcal{A}$. We may also equip the set of infinite words $\mathcal{S}^{\NN}$ with the Tychonoff topology, so that it is compact. The set of infinite paths is then a subshift of finite type (see \cite{Devaney2021}, §1.13), so it is a closed subset and compact. The function $\chi : \mathcal{S}^{\NN} \to \Sigma^{\NN}$ is continuous, and maps compact sets to compact sets. The image of the subshift of finite type is exactly $\mathcal{L}(\mathcal{A})$.
\end{proof}

A \textit{(pre)periodic word} in $\mathcal{S}^{\NN}$ is a word whose sequence of states is (pre)periodic. Analogously we may define \textit{(pre)periodic paths} and \textit{(pre)periodic allowable sequences}.

\begin{lemma} \label{lemma: nonempty regular language implies existence of periodic path}
If $\mathcal{L}(\mathcal{A}) \neq \emptyset$, then there exists a periodic infinite path in $\mathcal{A}$.    
\end{lemma}

\begin{proof}
If $\mathcal{L}(\mathcal{A}) \neq \emptyset$, there exists an infinite path in $\mathcal{A}$. Since the set of states is finite, some state must be repeated in the path, that is, $w_k = w_{k + l}$ for $k \geq 0$ and $l \geq 0$. The word $\overline{w_k w_{k+1}\ldots w_{k+l-1}}$ is therefore an infinite periodic path in $\mathcal{A}$.
\end{proof}

As we will mainly be concerned with infinite paths in an automaton, we describe an operation that reduces into to a simpler one, while preserving the regular language.

\begin{definition} \label{definition: trimming}
Let $\mathcal{A} = (\mathcal{G}, \chi)$ be an automaton. \textit{Trimming} the automaton $\mathcal{A}$ is an iterative procedure that produces another automaton $\mathcal{A}_{\mathrm{trim}}$ as follows. From $\mathcal{A}$, we erase all states that have no transitions from them, and all transitions into those states. (More precisely, we take the corresponding subgraph of $\mathcal{G}$.) If the remaining collection of states $\mathcal{S}_1 \subseteq \mathcal{S}_0 \coloneqq \mathcal{S}$ contains states that have no transitions from them, we erase those as well and transitions to them, giving us a subset $\mathcal{S}_2 \subseteq \mathcal{S}_1$. This creates a sequence
$$
\mathcal{S}_0 \supseteq \mathcal{S}_1 \supseteq \mathcal{S}_2 \supseteq \cdots.
$$
As $\mathcal{S}$ is finite, this process must eventually stabilize, so that there is a final set of states $\mathcal{S}_{\mathrm{trim}}$, possibly empty, such that all states $A \in \mathcal{S}_{\mathrm{trim}}$ have a transition to some other state in $\mathcal{S}_{\mathrm{trim}}$. We call $\mathcal{A}_{\mathrm{trim}} = (\mathcal{G}_{\mathrm{trim}}, \chi|_{\mathcal{S}_{\mathrm{trim}}})$ the \textit{trimmed automaton}, where $\mathcal{G}_{\mathrm{trim}}$ is the subgraph of $\mathcal{G}$ restricted to the states in $\mathcal{S}_{\mathrm{trim}}$ and transitions between them.
\end{definition}

By construction, the trimmed automaton $\mathcal{A}_{\mathrm{trim}}$ has the property that every state has a transition from it. In particular, every finite path in $\mathcal{A}_{\mathrm{trim}}$ can be continued to an infinite path. We say that an automaton is \textit{trimmed} if $\mathcal{A}_{\mathrm{trim}} = \mathcal{A}$.

\begin{prop} \label{proposition: trimmed automaton has same regular language}
$\mathcal{L}(\mathcal{A_{\mathrm{trim}}}) = \mathcal{L}(\mathcal{A})$.
\end{prop}

\begin{proof}
Since $\mathcal{G}_{\mathrm{trim}}$ is a subgraph of $\mathcal{G}$, evidently $\mathcal{L}(\mathcal{A}_{\mathrm{trim}}) \subseteq \mathcal{L}(\mathcal{A})$. Suppose that $s \in \Sigma^{\NN}$ is recognized by $\mathcal{A}$, being equal to $\chi(w)$ for some infinite path $w$ in $\mathcal{A}$. If $w$ were not an infinite path in $\mathcal{A}_{\mathrm{trim}}$, it is because one of its states is contained in $\mathcal{S}_i \setminus \mathcal{S}_{i+1}$ for some $i \geq 0$. Let $i_0$ be smallest $i$ such that a state of $w$ is in $\mathcal{S}_i \setminus \mathcal{S}_{i+1}$. If $w_k = A$ was erased at this step, it is because it had no transitions from it; but clearly it transitions to $w_{k+1}$, which by assumption has not been erased at an earlier step.
\end{proof}

The following definitions and results will be useful in Sections \ref{section: characteristic states automaton} and \ref{section: algorithm for first ray connection}.

\begin{definition} \label{definition: cross section}
A collection of states $\mathcal{C} \subseteq \mathcal{S}$ of an automaton $\mathcal{A}$ is called a \textit{cross section} if every infinite path in $w$ must eventually pass through a state in $\mathcal{C}$; that is, for some $k \geq 0$, $w_k \in \mathcal{C}$.
\end{definition}

Cross sections partially describe the shape of an automaton:

\begin{prop} \label{proposition: cross section shape of automaton}
Suppose that $\mathcal{C} = \{C_1, \ldots, C_k\}$ is a cross section of a trimmed automaton $\mathcal{A}$. Then, for each $C_i$ and $C_j$ in $\mathcal{C}$, there exists a finite collection of finite paths $\mathcal{P}_{C_iC_j}$ such that, for all $w \in \mathcal{P}_{C_iC_j}$, the first state of $w$ is $C_i$, the last state of $w$ is $C_j$, no other state of $w$ is in $\mathcal{C}$, and every path in $\mathcal{A}$ starting at $C_i$ must follow one of the paths in $\bigcup_{j} \mathcal{P}_{C_i C_j}$.

Similarly, for every $A \notin \mathcal{C}$ and $C_j \in \mathcal{C}$, there exists a finite collection of finite paths $\mathcal{P}_{AC_j}$ such that for $w \in \mathcal{P}_{AC_j}$, the first state of $w$ is $A$, the last state is $C_j$, no other state of $w$ is in $\{A\} \cup \mathcal{C}$, no state of $w$ repeats, and every path in $\mathcal{A}$ starting at $A$ must follow one of the paths in $\bigcup_j \mathcal{P}_{A C_j}$.
\end{prop}

Intuitively, the automaton takes the form below, illustrated for $|\mathcal{C}| = 2$, where every arrow represents a distinct path of finite length through states not in $\mathcal{C}$, and $\mathcal{T} = \mathcal{S \setminus \mathcal{C}}$:
% https://q.uiver.app/#q=WzAsMyxbMCwxLCJDXzEiXSxbMiwxLCJDXzIiXSxbMSwwLCJcXG1hdGhjYWx7VH0iXSxbMCwwLCIiLDEseyJvZmZzZXQiOjEsImFuZ2xlIjotOTB9XSxbMCwwLCIiLDEseyJyYWRpdXMiOjEsImFuZ2xlIjotOTB9XSxbMiwwLCIiLDEseyJvZmZzZXQiOjEsImN1cnZlIjoxfV0sWzIsMCwiIiwxLHsib2Zmc2V0IjotMSwiY3VydmUiOjF9XSxbMiwxLCIiLDEseyJvZmZzZXQiOjEsImN1cnZlIjotMX1dLFsyLDEsIiIsMSx7Im9mZnNldCI6LTEsImN1cnZlIjotMX1dLFsxLDEsIiIsMSx7InJhZGl1cyI6MSwiYW5nbGUiOjkwfV0sWzEsMSwiIiwxLHsib2Zmc2V0IjoxLCJhbmdsZSI6OTB9XSxbMCwxLCIiLDEseyJjdXJ2ZSI6LTF9XSxbMCwxLCIiLDEseyJvZmZzZXQiOi0yLCJjdXJ2ZSI6LTF9XSxbMSwwLCIiLDEseyJvZmZzZXQiOi0yLCJjdXJ2ZSI6LTF9XSxbMSwwLCIiLDEseyJjdXJ2ZSI6LTF9XV0=
\[\begin{tikzcd}
	& {\mathcal{T}} & \\
	{C_1} && {C_2}
	\arrow[shift right, curve={height=6pt}, from=1-2, to=2-1]
	\arrow[shift left, curve={height=6pt}, from=1-2, to=2-1]
	\arrow[shift right, curve={height=-6pt}, from=1-2, to=2-3]
	\arrow[shift left, curve={height=-6pt}, from=1-2, to=2-3]
	\arrow[shift right, from=2-1, to=2-1, loop, in=145, out=215, distance=10mm]
	\arrow[from=2-1, to=2-1, loop, in=150, out=210, distance=5mm]
	\arrow[curve={height=-6pt}, from=2-1, to=2-3]
	\arrow[shift left=2, curve={height=-6pt}, from=2-1, to=2-3]
	\arrow[shift left=2, curve={height=-6pt}, from=2-3, to=2-1]
	\arrow[curve={height=-6pt}, from=2-3, to=2-1]
	\arrow[from=2-3, to=2-3, loop, in=330, out=30, distance=5mm]
	\arrow[shift right, from=2-3, to=2-3, loop, in=325, out=35, distance=10mm]
\end{tikzcd}\]
This illustration should not be taken literally as an automaton equivalent to $\mathcal{A}$, but rather how all the infinite paths in $\mathcal{A}$ arise from a given starting state, through concatenation of a finite set of choices of first returns from $\mathcal{C}$ to $\mathcal{C}$ (and possibly first transitions from $\mathcal{T}$ to $\mathcal{C}$).

\begin{proof}
If a path $w$ in $\mathcal{A}$ starting at $C_i$ did not pass through $\mathcal{C}$ again, then $\sigma(w)$ would avoid $\mathcal{C}$, which is not possible. Moreover, any path starting at $C_i$ must return to $\mathcal{C}$ in an uniformly bounded number of transitions. Otherwise, if there were paths that took arbitrarily long to return to $\mathcal{C}$, since the number of states is finite, the paths of length greater than a fixed number would have to repeat some state. But this produces an infinite periodic path which avoids $\mathcal{C}$. We may then partition the set of paths starting at $C_i$ with respect to their first return to $\mathcal{C}$ at some state $C_j$, giving us the finite union $\bigcup_{j} \mathcal{P}_{C_i C_j}$. Similarly, for $A \notin \mathcal{C}$, a path starting at $A$ cannot take arbitrarily long to pass through a state in $\mathcal{C}$, otherwise some infinite periodic path would avoid $\mathcal{C}$.
\end{proof}

\subsection{The Hubbard tree automaton} \label{subsection: the hubbard tree automaton}

Let $F = \bigsqcup_{i = 1}^N I_j$ be the forbidden region of the Hubbard tree $H_c$ of $f_c$. For each of the angle intervals $I_j = (\phi^j_{-}, \phi^j_{+})$ associated to attaching points of the preimage tree, let $S$ be the set consisting of all angles $\phi^j_{\pm}$, their forward images under doubling, and the angles $\{0, 1/2\}$. We call $S$ the set of \textit{marked angles}. As will be discussed later, it may be enlarged arbitrarily by including finitely many other preperiodic angles and their forward images. By construction, $S$ is forward invariant.

Let $\mathcal{S}_{\mathrm{ang}}$ be a set of distinct symbols, each associated to a different marked angle $\phi \in S$, and $\mathcal{S}_{\mathrm{arc}}$ be another set of symbols, each associated to a distinct open arc $(\phi_1, \phi_2)$ of adjacent angles $\phi_1, \phi_2 \in S$, oriented counterclockwise starting from $0$. Let $\mathcal{S} = \mathcal{S_{\mathrm{ang}}} \cup \mathcal{S}_{\mathrm{arc}}$.

\begin{remark}
As a slight abuse of notation, we will from now on interchangeably refer to elements of $\mathcal{S}$ as the angles of $S$ or open arcs with endpoints in $S$ that they correspond to, and vice-versa.
\end{remark}

We view the elements of $\mathcal{S}$ as vertices of a directed graph $\mathcal{G}$ as follows. Let $A, B \in \mathcal{S}$.

\begin{itemize}
    \item If $A \in \mathcal{S}_{\mathrm{arc}}$, then we add a unique directed edge $A \to B$ if and only if the image of $A$ under doubling contains $B$.

    \item If $A \in \mathcal{S}_{\mathrm{ang}}$, we add a unique directed edge $A \to B$ if and only if the image of $A$ is the marked angle $B$.
\end{itemize}

Each $A \in \mathcal{S}$ maps over a union of elements of $\mathcal{S}$ homeomorphically, since $0$ and $1/2$ are marked angles. There is a function $\chi: \mathcal{S} \to \{0,1\}$ that associates to each element $A \in \mathcal{S}$ the bit value $0$ or $1$ according to whether $A$ is contained in $[0,1/2)$ or $[1/2, 1)$, respectively. This is well defined since $0$ and $1/2$ are marked angles, and we note that $\chi(0) = 0$, $\chi(1/2) = 1$. The pair $(\mathcal{G}, \chi)$ therefore describes a trimmed automaton with set of states $\mathcal{S} = \mathcal{S}_{\mathrm{arc}} \cup \mathcal{S}_{\mathrm{ang}}$, which essentially just encodes the dynamics of doubling on $S^1$ from the partition given by the marked angles $S$ as endpoints of arcs.

There is a distinguished subset $\mathcal{F} \subseteq \mathcal{S}$ of states that are contained in the forbidden region $F$ of $H_c$. We call these the \textit{forbidden states}. We then erase from $\mathcal{G}$ the vertices in $\mathcal{F}$, and trim the resulting automaton (Definition \ref{definition: trimming}).

\begin{definition}
The automaton $\mathcal{H}_S(c)$ obtained by erasing the states $\mathcal{F}$ and trimming is called the \textit{Hubbard tree automaton} of $f_c$. 
\end{definition}

To solidify notation, we denote the set of states of $\mathcal{H}_S(c)$ by $\mathcal{S}_{\mathrm{good}} \subseteq \mathcal{S}$ (the \textit{good states}), its directed graph by $\mathcal{G}_{\mathrm{good}}$, and identify $\chi$ with its restriction to $\mathcal{S}_{\mathrm{good}}$. We also call the allowable sequences in $\mathcal{L}(c) \subseteq\{0,1\}^{\NN}$ \textit{allowable bit sequences}. If $\phi$ is the angle given by the bit sequence of a path $w$ in $\mathcal{H}_S(c)$ seen binary expansion, we say that $\phi$ is the \textit{angle associated to $w$}. Note that if $\phi$ is associated to $w$, then $2\phi$ is associated to $\sigma(w)$.

We include the dependence of the automaton on $S$ as part of the notation for the following reason. Even though $\mathcal{H}_S(c)$ can be canonically defined based on the uniqueness of the forbidden region, and just establishing the set of marked angles to be the forward images of its endpoints along with $\{0, 1/2\}$, the construction of the automaton may be applied to any enlargement $S' \supseteq S$ of the set of marked angles by including finitely many preperiodic angles and their forward orbits. The resulting automata $\mathcal{H}_{S'}(c)$ will be distinct, but share the same regular language: see Corollary \ref{corollary: regular language is the same} below.

Intuitively, an infinite path $w = A_1A_2A_3 \ldots$ in $\mathcal{H}_S(c)$ is a infinite sequence of arcs and marked angles, mapping one over the next, that never falls into the forbidden region, and its corresponding bit sequence is the itinerary with respect to the partition of $S^1$ by $\{0, 1/2\}$.

\begin{definition} \label{definition: arc words and arc paths}
We say that a word $w \in \mathcal{S}^{\leq \omega}$ is an \textit{arc word} if all of its states are arc states. Analogously, a path in $\mathcal{H}(c)$ is an \textit{arc path} if it passes through only arc states.
\end{definition}

Let $f: S^1 \to S^1$ be the doubling map, $A = A_0$ a state in $\mathcal{S}$, and $A_0 \to A_1 \to \cdots \to A_k$ a finite path in graph $\mathcal{G}$. We define a subset of $S^1$ corresponding to this itinerary:
\begin{equation} \label{equation: definition of J arcs}
J_{A_0A_1\ldots A_k} \coloneqq A_0 \cap f^{-1}(A_1) \cap \cdots \cap f^{-k}(A_k),
\end{equation}
such that $J_{A_0A_1\ldots A_k}$ consists of all $\phi \in S^1$ that follow the itinerary $A_0 \to A_1 \to \cdots \to A_k$ under doubling. As we include $\{0,1/2\}$ in the marked angle set, for $k > 0$ and $0 \leq i \leq k$ we have
\begin{equation}
f^i(J_{A_0A_1\ldots A_k}) = J_{A_i \ldots A_k},
\end{equation}
where the maps are homeomorphisms onto their images. If any of the $A_0, \ldots, A_k$ is an angle state, then $J_{A_0A_1\ldots A_k}$ is a single angle, and a preimage of a marked angle. If all $A_i$ are arc states, then $J_{A_0A_1\ldots A_k}$ is an open arc contained in $A$ whose endpoints are iterated preimages of marked angles. We will call the sets of this form \textit{itinerary arcs}, or \textit{$J$-arcs}.

If $A_0$ is an arc state and $k \geq 1$, then it is partitioned into a union $\bigsqcup J_{A_0 A_{i_1}\ldots A_{i_k}}$ over all possible finite paths $A_0 \to A_{i_1} \to \cdots \to A_{i_k}$ of length $k+1$. In general some of these subsets may not be good, that is, they iterate into a forbidden state for $\mathcal{H}_S(c)$.

The following is a basic lemma from real analysis:

\begin{lemma} \label{lemma: nested intersection of open intervals}
Let $\bigcap_{n \geq 0} (a_n, b_n)$ be a nested intersection of non-empty, open bounded intervals of $\RR$, whose diameters tend to $0$. Then the intersection is empty if and only if there exists $N > 0$ such that either for all $n \geq N$ we have $a = a_n$, or for all $n \geq N$, $b = b_n$.
\end{lemma}

\begin{proof}
If $a_n = a$ for all large $n$, then $\bigcap_{n \geq N} (a, b_n) = \emptyset$, since $b_n \to a$. The analogous fact holds if $b b_n$ for all $n$ large. Conversely, by passing to a subsequence we may assume that the endpoints are ordered as below:
$$
a_0 < a_1 < a_2 < \cdots < a_n < \cdots < b_n < \cdots < b_2 < b_1 < b_0,
$$
and because the diameters tend $0$, $l = \lim_n a_n = \lim b_n$ is the unique element in the intersection.
\end{proof}

\begin{prop} \label{proposition: regular language lands on tree}

Let $S'$ be an enlargement of the set of marked angles as above, and $\mathcal{H}_{S'}(c)$ the resulting automaton with regular language $\mathcal{L}_{S' }(c)$. An external ray $\R_c(\phi)$ lands at a point on the Hubbard tree $H_c$ if and only a binary expansion of $\phi$ is in $\mathcal{L}_{S'}(c)$.
\end{prop}

\begin{comment}
\begin{remark}
If $\phi$ above is a dyadic rational, we assume by convention that its binary expansion is the one with an infinite sequence of zeros.
\end{remark}
\end{comment}

\begin{proof}
Suppose that $\R_c(\phi)$ lands at $z \in H_c$. By Proposition \ref{proposition: forbidden region and external arguments}, the forward orbit of $\phi$ never enters the forbidden region. Since the states $\mathcal{S}'$ partition the circle, there is a unique state $A_i$ such that $2^i\phi \in A_i$. By construction of $S'_{\mathrm{good}}$, the word $A_0A_1A_2\ldots$ forms an infinite path in $\mathcal{G}_{\mathrm{good}}$. Moreover, its bit sequence is a binary expansion of $\phi$, since it coincides with its itinerary under doubling with respect to the partition $S^1 = [0,1/2) \cup [1/2, 1)$. (If $\phi$ is dyadic, it is the preimage of a marked angle, and the binary expansion associated to its bit sequence is the one with trailing zeros since $\chi(0) = 0$ and $\chi(1/2) = 1$.) 

Conversely, suppose that a binary expansion of $\phi$ is the image bit sequence of some infinite path $A_0A_1A_2\ldots$ in $\mathcal{H}_{S'}(c)$. We divide into two cases:

\textit{Case 1: for some $A_k$, we have $A_k \in \mathcal{S}'_{\mathrm{ang}}$.} Let $A_{k_0}$ be the first state in the path which is in $\mathcal{S}'_{\mathrm{ang}}$, that is, which corresponds to a marked angle. The subset $J_{A_0A_1\ldots A_{k_0}}$ is therefore an iterated preimage of a marked angle $\psi$, the word corresponds to the itinerary of $\psi$, and its orbit never passes through the forbidden region. Since its binary expansion coincides with that of $\phi$, we have $\psi = \phi$, and $\R_c(\phi)$ lands on $H_{c}$.

\textit{Case 2: $w$ is an arc path.} The itinerary arcs $J_{A_0A_1\ldots A_k}$ for $k \geq 0$ form a nested intersection
\begin{equation} \label{equation: nested arcs}
A_0 = J_{A_0} \supseteq J_{A_0 A_1} \supseteq \cdots \supseteq J_{A_0 A_1 \ldots A_k} \supseteq \cdots,
\end{equation}
which is possibly empty. All elements in the intersection share the same bit sequence and therefore the same binary expansion. Hence if the intersection is non-empty, it must consist of a single point, and the diameters of the arcs tend to $0$ in either case.

If $\psi \in \bigcap_{k \geq  0} J_{A_0 \ldots A_k}$, the same arguments in case 1 conclude that $\psi = \phi$. If the intersection is empty, then by Lemma \ref{lemma: nested intersection of open intervals} one of the endpoints of the $J_{A_0\ldots A_k}$ stabilizes to an angle $\eta$. We choose angles $\eta_k \in J_{A_0\ldots A_k}$ so that $\eta_k \to \eta$. The angles $\eta_k$ have binary expansions that agree with a given one for $\phi$ up to the $k$-th entry, and so $\eta_k \to \phi = \eta$. Moreover, each $\eta_k$ is contained in the set $C_k \coloneqq S^1 \setminus \bigcup_{i = 0}^{k}f^{-i}(F)$, which form a nested collection of closed sets. Then the limit $\phi$ is contained in the nested intersection $\bigcap_{k \geq 0} C_k = S^1 \setminus \bigcup_{i = 0}^{\infty}f^{-i}(F)$, meaning that it never iterates into the forbidden region.
\end{proof}

The proof of \ref{proposition: regular language lands on tree} actually reveals more information. If an arc path $w$ in $\mathcal{H}(c)$ is such that the intersection in (\ref{equation: nested arcs}) is empty, the angle $\phi$ associated to $w$ is necessarily an iterated preimage of a marked angle, and there is another path in $\mathcal{H}(c)$ whose bit sequence is a binary expansion of $\phi$ (namely, one passing through an angle state). This is because if an endpoint of the $J_{A_0 \ldots A_k}$ stabilizes, the infinite itinerary on good arcs corresponds to an angle whose ray lands on $H_c$ ``infinitesimally close'' to that endpoint. This is the case if and only if $\phi$ can be approximated by other angles whose rays land on $H_c$: see Remark \ref{remark: distinction} and Example \ref{example: rabbit tuned by basilica}. In this case, we say that the word $w = A_0 A_1 \ldots$ and the nested itinierary arcs $J_{A_0\ldots A_k}$ \textit{accumulate} on $\phi$.

Moreover, every angle $\phi$ whose binary expansion is recognized by $\mathcal{H}(c)$ either has a unique path $w$ in $\mathcal{H}(c)$ whose associated angle is $\phi$, or $\phi$ is an iterated preimage of a marked angle, and there are at most three paths whose bit sequences give $\phi$. One is always the path passing through an angle state, and the other two possible paths are the arc paths in $\mathcal{H}(c)$ accumulating on $\phi$ from either of its sides in $S^1$. If $\phi$ is an iterated preimage of a marked angle and it is not accumulated on by an infinite arc path, we say that it is an \textit{isolated angle of (or landing on) $H_c$}. Equivalently, it is an isolated angle in the set of all the angles whose rays land on $H_c$.

\begin{lemma} \label{lemma: dyadic angles accumulated on either side}
If $\phi$ is dyadic and $\R_c(\phi)$ lands on $H_c$, then $\phi$ is accumulated on either side by angles of rays also landing on $H_c$.
\end{lemma}

\begin{proof}
From the hypothesis, $\beta \in H_c$, since it is in the forward orbit of the landing point $z$ of $\R_c(\phi)$. Since $\beta$ has only one external argument, it is an extremity of $H_c$, which shows that the critical point is strictly preperiodic and the parameter $c$ is Misiurewicz of $\beta$ type. In particular, $H_c \subseteq J_c$.

Suppose $z$ is an extremity of $H_c$, and in particular postcritical, so that for some neighborhood $N$ of it we have $N \cap H_c = (z_*, z]$. We may take a sequence of biaccessible parameters $z_n \to z$ on $(z_*, z]$ such that, if $(\phi_n, \phi_n')$ are the external arguments of $z_n$ adjacent-most to $(z_n, z]$, we have $\phi_n, \phi_n'\to \phi$. This is because, by possibly passing to a subsequence, the $\phi_n, \phi_n'$ have to convergence to an external argument of $z$ by continuity, and there is a unique one.

If $z$ is not an extremity of $H_c$, then it is the critical point or precritical, and for some neighborhood $N$ of $z$ we have $H_c \cap N = (z_1, z_2)$. The same arguments above apply, where we need to approximate $\phi$ landing on $z$ only on the side of $(z_1, z_2)$ that it lands on.
\end{proof}

The lemma above shows that if $\phi$ is dyadic and a binary expansion of $\phi$ is in $\mathcal{L}_{S'}(c)$, then both are. In fact, both binary expansions are realized by arc paths in $\mathcal{H}(c)$, accumulating on $\phi$ from above and below. This also shows that we could have attributed the bit values $\chi(0) = 1$ and $\chi(1/2) = 0$ and still gotten the same regular language.

\begin{cor} \label{corollary: regular language is the same}
The regular language $\mathcal{H}_{S'}(c)$ does not depend on the enlargement $S'$ of $S$.
\end{cor}

\begin{proof}
The only ambiguity between angles landing on $H_c$ and their binary expansions on $\mathcal{L}_{S'}(\phi)$ is with respect to dyadic rationals; this is resolved by the dicussion following Lemma \ref{lemma: dyadic angles accumulated on either side}.
\end{proof}

Based on this corollary, we may denote the regular language of the automaton simply by $\mathcal{L}(c)$. The expression $\mathcal{H}(c)$ is taken to mean \textit{a} Hubbard tree automaton for $f_c$ for some choice of set of marked angles, when that choice is immaterial. Moreover, we may use the notation $\mathcal{H}(\theta)$ for the Hubbard tree automaton if $c = c_{\theta}$.

\begin{remark} \label{remark: distinction}
It is important to be careful about the distiction between a bit sequence being recognized by the automaton, and a path being in $\mathcal{H}(c)$. It is possible that a word $w \in \mathcal{S}^{\NN}$ has an associated bit sequence $\chi(w)$ that is recognized by the automaton, but $w$ itself is not a path in $\mathcal{H}(c)$. This is the case if the angle $\phi$ given by $\chi(w)$ is a marked angle that is isolated from one side. Examples \ref{example: rabbit tuned by basilica} and \ref{example: 43/127} exhibit these isolated angles.
\end{remark}

The following is a useful property that, up to an enlargement of the set of marked angles, we can assume of good arc states.

\begin{lemma} \label{lemma: enlargement for endpoint landing on tree}
Given any enlargement $S'$ of the set of marked angles, there exists a further enlargement $S'' \supseteq S'$ with the property that, for any good arc state $A = (\phi_1, \phi_2)$, the rays $\R_c(\phi_1)$ and $\R_c(\phi_2)$ land on $H_c$.
\end{lemma}

\begin{proof}
If any $\phi \in S$ is such that the landing point $z$ of $\R_c(\phi)$ is not in $H_c$, consider the regular arc $[z, z^*]$ connects $z$ to $H_c$, as in the proof of Proposition \ref{proposition: forbidden region and external arguments}. The point $z^* \in H_c$ is either an extremity of $H_c$ or is a branch point of $J_c$, so that it has at least three external arguments. In either case, $z_*$ is preperiodic. We adjoin to $S'$ all the external arguments of $z^*$, and we do this for every other $\phi' \in S$ whose ray does not land on $H_c$.

Now let $A = (\phi_1, \phi_2) \in \mathcal{S}''_{\mathrm{arc}} \cap \mathcal{S}''_{\mathrm{good}}$, and suppose by way of contradiction that $\R_c(\phi_1)$ lands at $z \notin H_c$. Letting $[z, z^*]$ be the regular connecting $z$ to $H_c$ as above, there are two external arguments $\eta_1$ and $\eta_2$ of $z^*$ whose rays are adjacent-most to $[z, z^*]$. The arc $A$ is contained in $(\eta_1, \eta_2)$, and a suitable application of the arguments in the proof of Proposition \ref{proposition: forbidden region and external arguments} guarantee that $A$ eventually iterates into the forbidden region. This is a contradiction with the assumption that $A \in \mathcal{S}_{\mathrm{good}}$.
\end{proof}

An important enlargement of the set of marked angles, which we will assume throughout the rest of this section and in Sections \ref{section: characteristic states automaton} and \ref{section: algorithm for first ray connection}, is described below.

\begin{itemize}
    \item We first adjoin all the external arguments of all branch points and extremities of $H_c$. By the No Wandering Triangle Theorem \cite{Thurston2009}, they are all preperiodic. These can be found a pullback procedure analogous to the one described in Section \ref{subsection: finding the attaching points}. If $f_c$ is critically periodic, instead of the external arguments of the extremities we include the external arguments of the roots if the periodic Fatou components.

    \item If $c$ is in the $1/2$-wake of the Mandelbrot set, so that the $\alpha$ fixed point has external arguments $1/3$ and $2/3$, we also include them.

    \item If $S$ includes the external argument of some preperiodic point in the interior of an edge of $H_c$, we adjoin all other external arguments of that point.

    \item If $c$ is a Misiurewicz parameter, we adjoin the external arguments of the critical point. If $f_c$ is critically periodic, we adjoin the two external arguments of the root of the periodic Fatou component containing $0$, and the two external arguments of the preimage of this root under the first return map.

    \item Finally, we consider the enlargement given by Lemma \ref{lemma: enlargement for endpoint landing on tree}, so that the endpoints of any arc $A \in \mathcal{S}_{\mathrm{arc}} \cap \mathcal{S}_{\mathrm{good}}$ have rays that land on $H_c$.
\end{itemize}

In order to constructively find the angles for the enlargement described by Lemma \ref{lemma: enlargement for endpoint landing on tree}, we may take $\phi \in S$ whose ray $\R_c(\phi)$ doesn't land on $H_c$, iterate it until it is first in a component $I_j$ forbidden region $F$, and pullback $I_j$ along the orbit of $\phi$. The endpoints of the arc containing $\phi$ will be new marked angles, along with the other external arguments of the common landing point of the rays.

\begin{lemma} \label{lemma: for good arc states, landing points of endpoints are not equal}
Assume that $S$ is sufficiently enlarged as above. Let $A = (\phi_1, \phi_2)$ be a state in $\mathcal{S}_{\mathrm{arc}} \cap \mathcal{S}_{\mathrm{good}}$, and let $z_1 = \R_c(\phi_1)$, $z_2 = \R_c(\phi_2)$. Then $z_1 \neq z_2$.
\end{lemma}

\begin{proof}
Suppose by way of contradiction that $z = z_1 = z_2$. Then $\phi_1$ and $\phi_2$ are external arguments of the same point $z$, and the Julia set must intersect the two open regions cut out by the ray pair $\R_c(\phi_1) \cup \{z\} \cup  \R_c(\phi_2)$. Otherwise, all angles on one of the arcs $(\phi_1, \phi_2)$ or $(\phi_2, \phi_1)$ would land at $z$, which is not possible. Let $V$ be the region corresponding to the angle sector $A = (\phi_1, \phi_2)$. As there are no other external arguments of $z$ in $A$ by hypothesis on $S$, the intersection $J_c \cap V$ is connected (that is, there is a single ``branch'' of the Julia set following into $V$). A suitable adaptatition of the proof of Proposition \ref{proposition: forbidden region and external arguments} asserts that $z$ is an iterated preimage of an attaching point of the preimage Hubbard tree, and $A$ eventually maps inside of the forbidden region, which contradicts $A \in \mathcal{S}_{\mathrm{good}}$.
\end{proof}

\begin{definition} \label{definition: regular arc associated to good arc state}
Let $A = (\phi_1, \phi_2)$ be a good arc state, and $z_1$, $z_2$ be the landing points of $\R_c(\phi_1)$, $\R_c(\phi_2)$, respectively. By Lemmas \ref{lemma: enlargement for endpoint landing on tree} and Lemma \ref{lemma: for good arc states, landing points of endpoints are not equal}, the regular arc $[z_1, z_2]$ is not a single point and is contained in the Hubbard tree. We call it the \textit{regular arc on the Hubbard tree associated to $A$}.
\end{definition}

\begin{lemma} \label{lemma: landing points of good arc states and transitions}
Suppose that $A$ and $B$ are good arc states, and $A \to B$ is a transition. Let $B_1 = (\eta_1, \eta_2)$ be the preimage of $B$ that is contained in $A = (\phi_1, \phi_2)$. If $u_1$, $u_2$ are the landing points of $\R_c(\eta_1)$, $\R_c(\eta_2)$ and $[z_1, z_2]$ is the regular arc associated to $A$, then $u_1 \neq u_2$, $[u_1, u_2] \subseteq [z_1, z_2]$, and there are no external arguments of $u_1$ and $u_2$ on $B_1$.
\end{lemma}

\begin{proof}
We have $u_1 \neq u_2$ since $f_c(u_1) \neq f_c(u_2)$ by Lemma \ref{lemma: for good arc states, landing points of endpoints are not equal}. For the second conclusion, it is sufficient to show that $u_1$ and $u_2$ are in the Hubbard tree. This is because there can be no branch points of $H_c$ in $[z_1, z_2]$ for which a branch leads to $u_j$, otherwise some external argument of it is a marked angle that is contained in $A$.

If $u_1$ did not land on $H_c$, there is a regular arc $[u, z^*]$ connecting it to $H_c$. Moreover, $z^* \in (z_1, z_2)$; if $z^* \notin [z_1,z_2]$, then $\R_c(\eta_1)$ would cross one of the rays $\R_c(\phi_1)$ or $\R_c(\phi_2)$, and if $z^* = z_1$ or $z_2$, there is another part of the Julia set ``branching off of'' $z_1$ or $z_2$ and contained in the region of the plane bounded by $\R_c(\phi_1) \cup [z_1, z_2] \cup \R_c(\phi_2)$, producing another external argument of it in $A$, which is not possible.

The image $f_c(u_1)$ is on the Hubbard tree, being the landing point of the ray $\R_c(2\eta_1)$ where $B = (2\eta_1, 2\eta_2) \in \mathcal{S}_{\mathrm{arc}} \cap \mathcal{S}_{\mathrm{good}}$, and our assumption on the set of marked angles being enlarged according to Lemma \ref{lemma: enlargement for endpoint landing on tree}. Hence $u_1$ is in the preimage Hubbard tree, and the arc $[u_1, z^*]$ as well; this implies that $z^*$ is an attaching point, and its external arguments are marked angles as well. This would likewise imply that some marked angle is in $A$. Finally, if there were an external argument of $u_1$ on $B_1$, the image would be a marked angle in $B$, a contradiction with our assumptions. The same arguments apply for $u_2$.
\end{proof}

The arguments in the proof of Lemma \ref{lemma: landing points of good arc states and transitions}, can be extended inductively to all finite paths $A_0 \to A_1 \to \cdots \to A_k$ in $\mathcal{H}(c)$ passing only arc states:

\begin{cor} \label{corollary: landing points of J arcs}
Let $A_0 \to A_1 \to \ldots \to A_k$ be an arc path in $\mathcal{H}(c)$, $A_0 = (\phi_1, \phi_2)$, and $J_{A_0A_1\ldots A_k} = (\psi_1, \psi_2)$. If $\R_c(\phi_j)$ lands at $z_j$ and $\R_c(\psi_j)$ lands at $v_j$ for $j = 1, 2$, then $v_1 \neq v_2$, $[v_1, v_2] \subseteq [z_1, z_2]$, and there are no external arguments of $v_1$ and $v_2$ on $J_{A_0\ldots A_k}$.
\end{cor}

In particular, Definition \ref{definition: regular arc associated to good arc state} can be extended to all itinerary arcs $J_{A_0 \ldots A_k}$ for which $A_0, \ldots, A_k \in \mathcal{S}_{\mathrm{arc}} \cap \mathcal{S}_{\mathrm{good}}$, so that $[v_1, v_2]$ is its \textit{associated regular arc} in the Hubbard tree. They in fact form a nested sequence with respect to the nesting of the $J_{A_0 \ldots A_k}$. For an arc path $w = A_0 A_1 \ldots$ in $\mathcal{H}(c)$, the nested intersection of the associated regular arcs gives the landing point on $H_c$ of $\R_c(\psi)$, where $\psi \in \bigcap_{k \geq 0} J_{A_0 \ldots A_k}$, or the landing point of the iterated preimage of a marked angle that the word accumulates to.

\subsection{An involution on arc states} \label{subsection: an involution on arc states}

Let $A = (\phi_1, \phi_2)$ be an arc state such that $\R_c(\phi_1)$ and $\R_c(\phi_2)$ have distinct landing points $z_1$ and $z_2$ on $H_c$ (for example, if $A$ is a good arc state). Also let $E_1$ and $E_2$ be the sets of external arguments of $z_1$ and $z_2$, respectively. These are unlinked sets of angles on $S^1$ and contained in $S$. There are exactly two open arcs on $S^1$ for which one endpoint is in $E_1$, the other is in $E_2$, and such that there are no elements of $E_1 \cup E_2$ in the interior. The arc $A = (\phi_1, \phi_2)$ is one of them, and let $A' = (\phi_2', \phi_1')$ be the other arc, where $0 < \phi_2'< \phi_1'< 1$, $\phi_2'\in E_2$ and $\phi_1'\in E_1$. Also by assumption on $S$, there are no marked angles in $A'$. If there were, either $\phi' \in A'\cap S$ is another external argument of $z_1$ or $z_2$, or the regular arc connecting the landing point of $\R_c(\phi')$ to $[z_1, z_2]$ produces a branch point $z_*$ of $H_c$ in $(z_1, z_2)$. The former is not possible by our hypotheses on $\phi_2'$ and $\phi_1'$, and the latter would imply that $S$ is missing an external argument of $z_*$ accessing it through the $A$ arc. Hence $A'$ is a state in the angle doubling automaton, that is, its endpoints are two adjacent marked angles.

\begin{lemma} \label{lemma: A' is good arc state}
$A$ is a good arc state if and only if $A'$ also is.
\end{lemma}

\begin{proof}
It is sufficient to prove one direction of the equivalence. If $A' \notin \mathcal{S}_{\mathrm{good}}$, it is because every path starting at $A'$ in graph $\mathcal{G}$ eventually passes through a forbidden state. This implies that every angle in $A'$ does not land on the Hubbard tree; this can only be possible if $(z_1, z_2) \subseteq U$ for some Fatou component $U$. But this implies also that every angle in $A$ eventually iterates into the forbidden region, and so $A \notin \mathcal{S}_{\mathrm{good}}$. More precisely, $z_1$ and $z_2$ eventually iterate to the root of the periodic Fatou component that contains the unique attaching point in the Fatou set, and $A$ iterates into the arc determined by the the ray pair landing at this root adjacent-most to the component.
\end{proof}

\begin{definition} \label{definition: reflection and involution}
We call $A' \in \mathcal{S}_{\mathrm{arc}} \cap \mathcal{S}_{\mathrm{good}}$ the \textit{reflection} of $A$. This defines a function $\iota : \mathcal{S}_{\mathrm{arc}} \cap \mathcal{S}_{\mathrm{good}} \to \mathcal{S}_{\mathrm{arc}} \cap \mathcal{S}_{\mathrm{good}}$ mapping $A \mapsto A'$ as above, which we call the \textit{involution} or \textit{reflection on good arc states}. It has a natural extension to words in $\mathcal{S}_{\mathrm{arc}} \cap \mathcal{S}_{\mathrm{good}}$ given by
\begin{equation}
\iota(w) = \iota(A_0A_1A_2\ldots) \coloneq \iota(A_0)\iota(A_1)\iota(A_2)\ldots.
\end{equation}
\end{definition}

It is clear from the definition of $A'$ that $\iota$ is indeed an involution, and that the regular arc of $H_c$ associated to $A' = \iota(A)$ is the same arc associated to $A$. Intuitively, the reflection maps one side of an edge of the Hubbard tree, represented by an arc state, to the side opposite to it on the same edge.

Due to Corollary \ref{corollary: landing points of J arcs} and its following discussion, we may also define the \textit{reflection} $\iota(J_{A_0 \ldots A_k})$ of an itinerary arc $J_{A_0 \ldots A_k}$, when $A_0 \to \cdots \to A_k$ is an arc path on $\mathcal{H}(c)$.

\begin{prop} \label{proposition: J arcs respect involution}
For an arc path $A_0 \to \cdots \to A_k$ in $\mathcal{H}(c)$ passing only through states in $\mathcal{S}_{\mathrm{arc}} \cap \mathcal{S}_{\mathrm{good}}$, we have
$$
\iota(J_{A_0\ldots A_k}) = J_{\iota(A_0\ldots A_k)}.
$$
\end{prop}

\begin{proof}
We prove this by induction on the length of such a path. For $k = 0$ there is nothing to prove. Assume it holds for all such paths of length $k$, and consider $A_0 \to \cdots \to A_k \to A_{k+1}$. Let $J_{A_0\ldots A_k A_{k+1}} = (\eta_1, \eta_2)$, $\iota(J_{A_0\ldots A_k A_{k+1}}) = (\eta_1', \eta_2')$, and $J_{\iota(A_0 \ldots A_k A_{k+1})} = (\eta_1'', \eta_2'')$. Also let $[x_1,x_2]$ be the regular arc associated to $J_{A_0 \ldots A_k A_{k+1}}$ and $\iota(J_{A_0 \ldots A_k A_{k+1}})$, and $[y_1, y_2]$ be the regular arc associated to $J_{\iota(A_0 \ldots A_k A_{k+1}})$.

We claim that $x_1 = y_1$ and $x_2 = y_2$. This is because $J_{A_1\ldots A_{k+1}} = (2\eta_1, 2\eta_2)$ and
$$
\iota(J_{A_1\ldots A_{k+1}}) = J_{\iota(A_1\ldots A_k A_{k+1})} = (2\eta_1'', 2\eta_2''),
$$
so that, for $j = 1, 2$, the rays $\R_c(2\eta_j)$, $\R_c(2\eta_j'')$ land on $f_c(x_j)$. But $\R_c(2\eta_j'')$ lands on $f_c(y_j)$, so that $f_c(x_j) = f_c(y_j)$. All of the points $x_j$ and $y_j$ are contained in the regular arc $[z_1, z_2]$ determined by $A_0$ as in Corollary \ref{corollary: landing points of J arcs}, and since the external arguments of the critical point are marked, $[z_1, z_2]$ maps injectively onto its image. Therefore $x_j = y_j$ for $j_2$. Finally, for any arc $J_{B_0 \ldots B_i} = (\psi_1, \psi_2)$, we are able to recover it and $\iota(J_{B_0 \ldots B_i})$ from just regular arc $[v_1, v_2]$ associated to it. The arcs will be the angle sectors bounded by the external arguments of $v_1$ and $v_2$ that are adjacent-most to the external arguments of the other point. This concludes that $\iota(J_{A_0 \ldots A_k A_{k+1}}) = J_{\iota(A_0 \ldots A_k A_{k+1})}$.
\end{proof}

This in particular guarantees that $J_{A_0\ldots A_k}$ and $J_{\iota(A_0\ldots A_k)}$ have the same associated regular arc.

\begin{cor} \label{corollary: involution preserves arc paths}
If $A \to B$ is a transition between states in $\mathcal{S}_{\mathrm{arc}} \cap \mathcal{S}_{\mathrm{good}}$, then $\iota(A) \to \iota(B)$ is also a transition. In particular, if $w$ is a finite or infinite word on $\mathcal{S}_{\mathrm{arc}} \cap \mathcal{S}_{\mathrm{good}}$, Then $w$ is a path in $\mathcal{H}(c)$ if and only if $\iota(w)$ also is.
\end{cor}

\section{Examples of Hubbard Tree Automata} \label{section: examples of hubbard tree automata}

\begin{comment}
Consider $\theta = 1/5 = .\overline{0011}$. It is a periodic angle of period $4$, and its parameter ray lands at the root of a hyperbolic component with center $c_{\theta}$, the Kokopelli polynomial/ The kneading sequence of $\theta$ is $\nu(\theta) = \overline{\texttt{110*}}$, and its upper kneading sequence is $\mathcal{A}(\nu(\theta)) = \overline{\texttt{1100}}$, which has internal address $1 \to 3 \to 4$. Note that the other choice of subsutition for $*$, $\overline{\texttt{1101}}$, has infinite internal address, which skips the entry $S_k = 4$:
$$
1 \to 3 \to 6 \to 7 \to 10 \to 11 \to 14 \to 15 \to 18 \to \cdots
$$

The component associated to $S_1 = 3$ is contained in the $1/3$-subwake of the component associated to $S_0 = 1$, Corresponding to the rabbit bifurcation off of the main cardioid. We only need to see in which $p/q$-subwake of $S_1 = 3$ the component for $S_2 = 4$ is contained in. According to the formula \ref{denominators}, we have $r = 1$ in this case, and so the $\rho_{\nu}$-orbit is already the internal address. Hence the angled internal address of $c_{\theta}$ is
$$
(1)_{1/3} \to (3)_{1/2} \to 4.
$$

The only satellite periodic orbit with branch number $q \geq 3$ is that corresponding to the $\alpha$ fixed point, with external arguments $\{1/7, 2/7, 4/7\}$. 
\end{comment}

\begin{exm}
Let $\theta = 3/7 = .\overline{011}$, whose postcritically finite parameter $c_{3/7}$ corresponds to the airplane polynomial, and whose filled Julia set is depicted in Figure \ref{fig: airplane}. It's kneading sequence is $\nu(\theta) = \texttt{10*}$, being $*$-periodic of period 3. The upper kneading sequence is $\mathcal{A}(\nu) = \texttt{100}$, which can be see from the fact that the internal address for $\texttt{101}$ starts with $1 \to 2 \to 4 \to \cdots$ and skips $3$.

\begin{figure}[ht]
    \centering
    \includegraphics[width=0.7\linewidth]{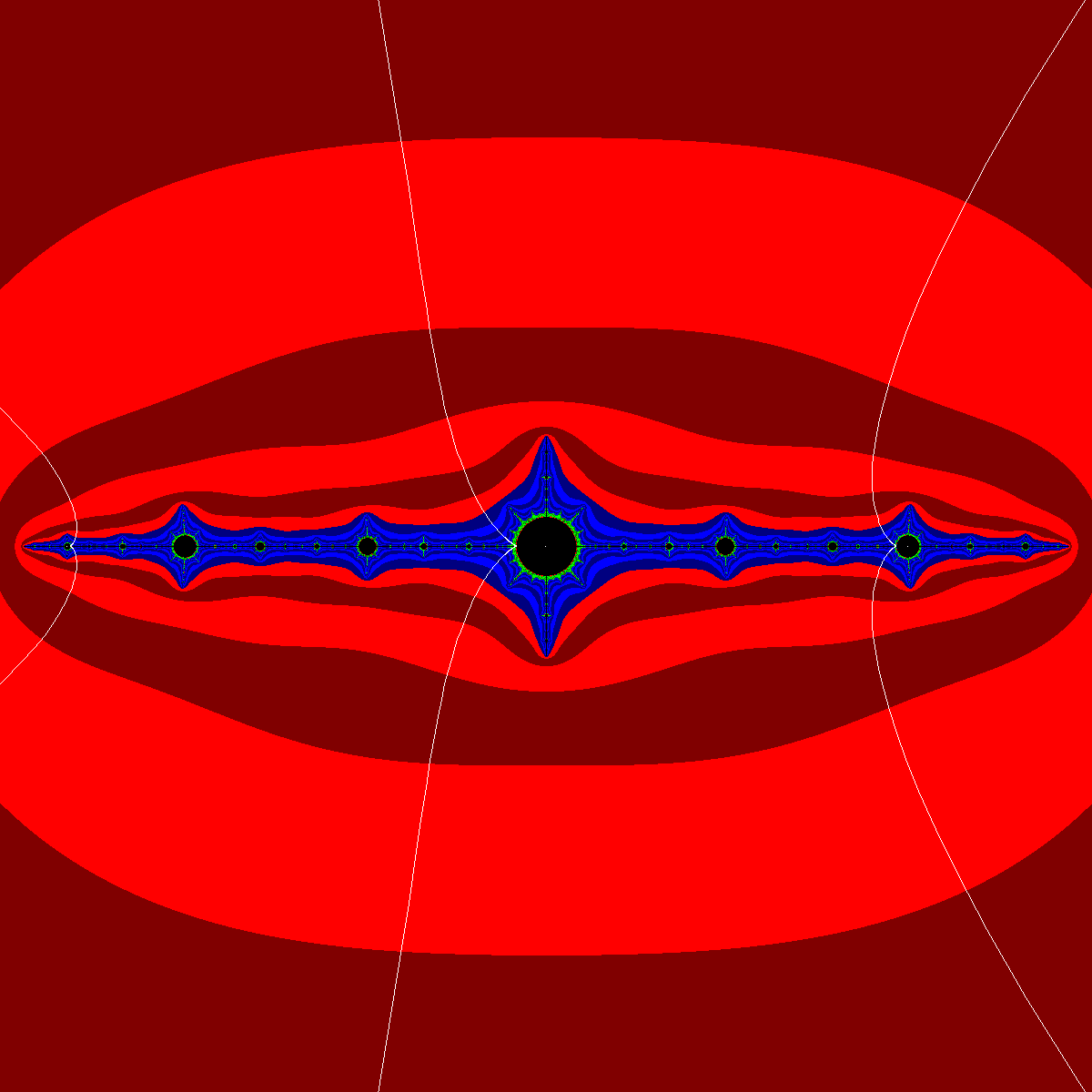}
    \caption{The airplane Julia set, with the orbits of rays landing at the roots of the periodic Fatou components marked. The forbidden region is $(6/7, 1/7)$.}
    \label{fig: airplane}
\end{figure}

Because the airplane polynomial is real, the Hubbard tree is an interval contained in $\RR$. The only attaching point of the preimage Hubbard tree is the center of the Fatou component which the angles $\{1/7, 6/7\}$ land at the root of, so that the forbidden region is $(6/7, 1/7)$. The set of marked angles is
$$
S = \left\{ \dfrac{1}{7}, \dfrac{2}{7}, \dfrac{3}{7}, \dfrac{4}{7}, \dfrac{5}{7}, \dfrac{6}{7}, \dfrac{1}{2}, \dfrac{0}{1}\right\},
$$
partitioning the circle into $8$ open arcs and $8$ common endpoints of adjacent arcs. The forbidden states for the Hubbard tree automaton are $0$, $(0,1/7)$, and $(6/7, 0)$. We may erase the state $1/2$, since it maps to the forbidden $0$, and similarly we erase the states $(3/7, 1/2)$ and $(1/2, 4/7)$, since they map to forbidden states. The automaton $\mathcal{H}_S(3/7)$, in the shape of a labeled directed graph, is depicted below:
% https://q.uiver.app/#q=WzAsMTAsWzIsMCwiMCJdLFsyLDEsIjAiXSxbMywxLCIxIl0sWzMsMCwiMSJdLFsxLDEsIjEiXSxbMSwwLCIwIl0sWzAsMSwiMSJdLFs0LDAsIjAiXSxbNCwxLCIxIl0sWzUsMCwiMCJdLFswLDFdLFsxLDJdLFsyLDNdLFszLDBdLFsxLDRdLFswLDVdLFszLDddLFsyLDhdLFs5LDddLFs3LDhdLFs4LDldLFs0LDVdLFs1LDZdLFs2LDRdLFsxLDMsIiIsMSx7ImN1cnZlIjotMX1dLFszLDEsIiIsMSx7ImN1cnZlIjotMX1dXQ==
\[\begin{tikzcd}
	& 0 & 0 & 1 & 0 & 0 \\
	1 & 1 & 0 & 1 & 1
	\arrow[from=1-2, to=2-1]
	\arrow[from=1-3, to=1-2]
	\arrow[from=1-3, to=2-3]
	\arrow[from=1-4, to=1-3]
	\arrow[from=1-4, to=1-5]
	\arrow[curve={height=-6pt}, from=1-4, to=2-3]
	\arrow[from=1-5, to=2-5]
	\arrow[from=1-6, to=1-5]
	\arrow[from=2-1, to=2-2]
	\arrow[from=2-2, to=1-2]
	\arrow[curve={height=-6pt}, from=2-3, to=1-4]
	\arrow[from=2-3, to=2-2]
	\arrow[from=2-3, to=2-4]
	\arrow[from=2-4, to=1-4]
	\arrow[from=2-4, to=2-5]
	\arrow[from=2-5, to=1-6]
\end{tikzcd}\]
The center square represents the mappings over the arc states, while the two period $3$ cycles that they may ``escape'' to are the cycles of the marked angles $\{1/7, 2/7, 4/7\}$ and $\{3/7, 6/7, 5/7\}$. Since those marked angles are not isolated landing on $H_{3/7}$, and therefore are represented by bit sequences over the arc states, we may erase the angle states and get a simpler automaton having the same regular language:
% https://q.uiver.app/#q=WzAsNCxbMCwwLCIwIl0sWzAsMSwiMCJdLFsxLDEsIjEiXSxbMSwwLCIxIl0sWzIsM10sWzAsMV0sWzMsMF0sWzEsMl0sWzEsMywiIiwxLHsiY3VydmUiOi0xfV0sWzMsMSwiIiwxLHsiY3VydmUiOi0xfV1d
\[\begin{tikzcd}
	0 & 1 \\
	0 & 1
	\arrow[from=1-1, to=2-1]
	\arrow[from=1-2, to=1-1]
	\arrow[curve={height=-6pt}, from=1-2, to=2-1]
	\arrow[curve={height=-6pt}, from=2-1, to=1-2]
	\arrow[from=2-1, to=2-2]
	\arrow[from=2-2, to=1-2]
\end{tikzcd}\]
One may show that a bit sequence is recognized by this automaton if and only if there are no more than two consecutive $0$'s and two consecutive $1$'s anywhere in the sequence.
\end{exm}

\begin{exm} \label{example: 1/6}
Let $c = i$, which is the landing point of the single parameter ray $\R_M(1/6)$, and let $f_{1/6}$ be the corresponding postcritically finite polynomial. The kneading sequence of $1/6$ is $\nu = \texttt{1}\overline{\texttt{10}}$, and its angled internal address is
$$
(1)_{1/3} \to (3)_{1/2} \to (5)_{1/2} \to (7)_{1/2} \to \cdots
$$
It is straightforward to verify by induction that the internal address consists of all odd integers in increasing order. Moreover, by (\ref{equation: denominators}), the only denominator $q \geq 3$ is at the entry $S_0 = 1$. Since the critical orbit falls into a trivial orbit portrait, with a period $2$ cycle of rays having external arguments $1/3$ and $2/3$ respectively, no extremity of the Hubbard tree is an attaching point of the preimage tree.

The only satellite orbit portrait with valence $v \geq 3$ corresponds to the $\alpha$ fixed point, with external arguments $\{1/7, 2/7, 4/7\}$. The preperiodic preimage $-\alpha$ has external arguments $\{1/14, 9/14, 11/14\}$, and is not a branch point since no forward image of $1/6$ is in the arc $(11/14, 1/14)$. This shows that the forbidden region is $(11/14, 1/14)$, as depicted in Figure \ref{fig: 1/6 forbidden region}.

The set of marked angles is
$$
S = \left\{ \dfrac{0}{1}, \dfrac{1}{14}, \dfrac{1}{7}, \dfrac{2}{7}, \dfrac{1}{2}, \dfrac{4}{7}, \dfrac{11}{14} \right\},
$$
and the forbidden states are $0/1, (0, 1/14), (11/14, 0)$. A representation of $\mathcal{H}_S(1/6)$ is given below:
% % https://q.uiver.app/#q=WzAsMTAsWzEsMywiMCJdLFsyLDMsIjAiXSxbMSwyLCIxIl0sWzIsMiwiMSJdLFszLDMsIjAiXSxbMywxLCIxIl0sWzIsMCwiMSJdLFsyLDEsIjAiXSxbMSwxLCIwIl0sWzAsMSwiMCJdLFswLDFdLFsyLDBdLFszLDJdLFsxLDRdLFszLDFdLFszLDQsIiIsMix7ImN1cnZlIjotMX1dLFs0LDMsIiIsMix7ImN1cnZlIjotMX1dLFs0LDVdLFs1LDZdLFs3LDZdLFs2LDhdLFs4LDddLFsyLDldLFs5LDhdLFszLDddLFsxLDJdXQ==
\[\begin{tikzcd}
	&& 1 & \\
	0 & 0 & 0 & 1 \\
	& 1 & 1 \\
	& 0 & 0 & 0
	\arrow[from=1-3, to=2-2]
	\arrow[from=2-1, to=2-2]
	\arrow[from=2-2, to=2-3]
	\arrow[from=2-3, to=1-3]
	\arrow[from=2-4, to=1-3]
	\arrow[from=3-2, to=2-1]
	\arrow[from=3-2, to=4-2]
	\arrow[from=3-3, to=2-3]
	\arrow[from=3-3, to=3-2]
	\arrow[from=3-3, to=4-3]
	\arrow[curve={height=-6pt}, from=3-3, to=4-4]
	\arrow[from=4-2, to=4-3]
	\arrow[from=4-3, to=3-2]
	\arrow[from=4-3, to=4-4]
	\arrow[from=4-4, to=2-4]
	\arrow[curve={height=-6pt}, from=4-4, to=3-3]
\end{tikzcd}\]
The first two rows of states are the angles states, and the bottom two rows are the arc states. The single periodic cycle among the angles states corresponds to the periodic orbit of rays landing at the $\alpha$-fixed point. Again, since no angle $\phi$ whose ray lands on $H_{1/6}$ is isolated, we may erase the angle states, and the resulting automaton below has the same regular language:
% % https://q.uiver.app/#q=WzAsNSxbMCwxLCIwIl0sWzEsMSwiMCJdLFswLDAsIjEiXSxbMSwwLCIxIl0sWzIsMSwiMCJdLFswLDFdLFsyLDBdLFszLDJdLFsxLDRdLFszLDFdLFszLDQsIiIsMix7ImN1cnZlIjotMX1dLFs0LDMsIiIsMix7ImN1cnZlIjotMX1dLFsxLDJdXQ==
\[\begin{tikzcd}
	1 & 1 & \\
	0 & 0 & 0
	\arrow[from=1-1, to=2-1]
	\arrow[from=1-2, to=1-1]
	\arrow[from=1-2, to=2-2]
	\arrow[curve={height=-6pt}, from=1-2, to=2-3]
	\arrow[from=2-1, to=2-2]
	\arrow[from=2-2, to=1-1]
	\arrow[from=2-2, to=2-3]
	\arrow[curve={height=-6pt}, from=2-3, to=1-2]
\end{tikzcd}\]
\end{exm}

\begin{exm} \label{example: rabbit tuned by basilica}
Let $\theta = 10/63 = .\overline{001010}$, which has kneading sequence $\nu = \texttt{11011*}$, upper kneading sequence $\mathcal{A}(\nu) = \texttt{110111}$, and internal address
$$
(1)_{1/3} \to (3)_{1/2} \to 6.
$$
The parameter $c_{10/63}$ is the center of a hyperbolic component of period $6$, whose polynomial corresponds to a tuning of the $1/3$-rabbit polynomial by the Basilica polynomial $z^2 - 1$. As in Example \ref{example: 1/6}, the $\alpha$ fixed point is the only branch point of $H_{3/16}$. However, the preimage point $-\alpha$ is not in $H_{3/16}$, and the forbidden region is the arc $F = (40/63, 5/63)$, whose attaching point is the unique attaching point which is the root of a Fatou component. A representation of $\mathcal{H}_S(10/63)$ for $f_{10/63}$ is given below:
% https://q.uiver.app/#q=WzAsOSxbMiwxLCIwIl0sWzIsMiwiMCJdLFsxLDIsIjEiXSxbMSwzLCIwIl0sWzIsMCwiMCJdLFszLDMsIjAiXSxbMCwzLCIwIl0sWzAsMCwiMSJdLFszLDAsIjEiXSxbMCwxXSxbMSwyXSxbMiwwXSxbMiwzXSxbNSwzXSxbMyw2XSxbNiw3XSxbNyw0XSxbNCw4XSxbMCw0XSxbOCw1XV0=
\[\begin{tikzcd}
	1 && 0 & 1 \\
	&& 0 \\
	& 1 & 0 \\
	0 & 0 && 0
	\arrow[from=1-1, to=1-3]
	\arrow[from=1-3, to=1-4]
	\arrow[from=1-4, to=4-4]
	\arrow[from=2-3, to=1-3]
	\arrow[from=2-3, to=3-3]
	\arrow[from=3-2, to=2-3]
	\arrow[from=3-2, to=4-2]
	\arrow[from=3-3, to=3-2]
	\arrow[from=4-1, to=1-1]
	\arrow[from=4-2, to=4-1]
	\arrow[from=4-4, to=4-2]
\end{tikzcd}\]

The inner cycle of period $3$ consits of the all the arc states and corresponds to the $\alpha$ fixed point, which has external rays $\{1/7, 2/7, 4/7\}$, and the outer cycle is the collection of angle states, representing the orbit of $10/63$. All angles whose rays land on $H_{10/63}$, except for the external arguments of $\alpha$, are in the grand orbit of $10/63$; the set of these landing points is discrete, and so each of these angles is isolated as seen landing on the Hubbard tree. Erasing the angle states does not give an automaton having the same regular language, contrary to the previous examples. This also illustrates the need to consider the marked angles as states separately, and not just the arc states as closed arcs.

\begin{figure}[ht]
    \centering
    \includegraphics[width=0.7\linewidth]{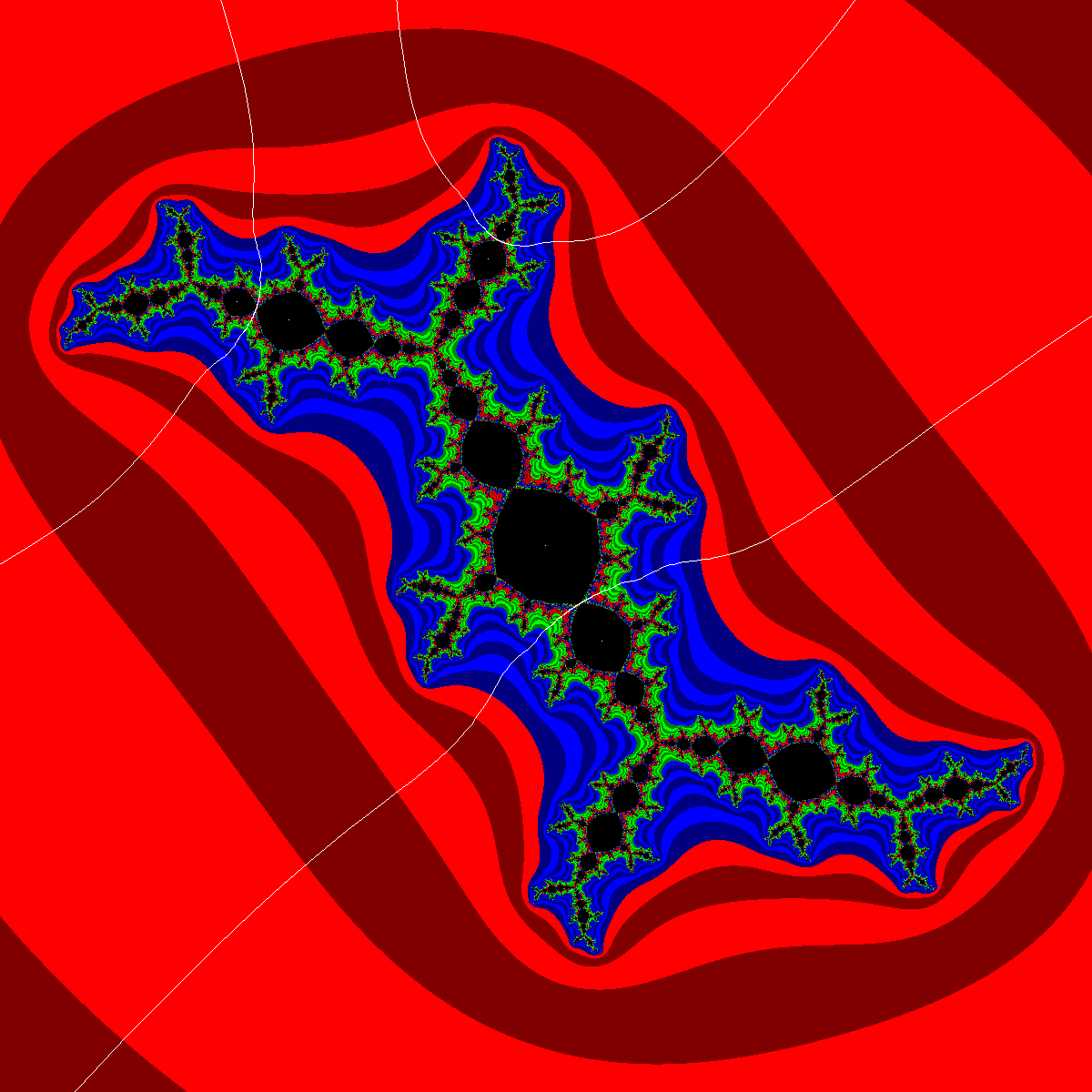}
    \caption{The filled Julia set of $f_{10/63}$, and the external rays landing at the roots of the periodic Fatou components. The intersection $J_{10/63} \cap H_{10/63}$ is discrete, except for the $\alpha$ fixed point.}
    \label{fig: rabbit tuned by basilica}
\end{figure}

\end{exm}

\begin{exm} \label{example: 43/127}
Let $\theta = 43/127 = .\overline{0101011}$. This is a periodic angle under doubling of period $7$, whose kneading sequence is $\nu = \overline{\texttt{101010*}}$. It has upper kneading sequence $\mathcal{A}(\nu(\theta)) = \overline{\texttt{1010100}}$, and angled internal address
$$
(1)_{1/2} \to (2)_{1/4} \to 7.
$$
In order to find the periodic branch points of $H_{43/127}$, we need only consider the entry $(2)_{1/4}$. By Proposition \ref{proposition: periodic approximation lands at root}, the angle $\varphi_2 = .\overline{01} = 1/3$ lands at the hyperbolic component associated to that entry, with partner angle $2/3 = .\overline{10}$. Substituting $t = 1/4$ in equation (\ref{equation: theta minus}), the angle $\theta_-$ landing at the $1/4$-sublimb of the main cardioid from the left is
$$
\theta_- = \left \lceil(2^4 - 1)\sum_{b = 1}^{10}\dfrac{\lceil b/4 \rceil - 1}{2^b}\right\rceil/ (2^4 - 1) = \dfrac{1}{15} = .\overline{0001},
$$
with partner angle $\theta_+ = 2/15 = .\overline{0010}$. By tuning, the angle $.\overline{01010110} = 86/255$ lands at the $1/4$-sublimb of the hyperbolic component associated to the entry $(2)_{1/4}$. Along with its partner angle, it generates the satellite orbit portrait
$$
A_1 = \left\{\dfrac{86}{255}, \dfrac{89}{255}, \dfrac{101}{255}, \dfrac{149}{255}\right\}, \quad A_2 = \left\{ \dfrac{172}{255}, \dfrac{178}{255}, \dfrac{202}{255}, \dfrac{43}{255} \right\}.
$$

By the iterated pullback procedure described for finding attaching points of $H$, we get the forbidden region
$$
F = \left( \dfrac{89}{510}, \dfrac{101}{510} \right) \sqcup \left( \dfrac{101}{510}, \dfrac{149}{510}\right) \sqcup \left( \dfrac{106}{127}, \dfrac{11}{127}\right),
$$
where the last angle interval corresponds to the periodic postcritical attaching point $c_6$. Not that $101/510$ is an isolated angle of the Hubbard tree, that is, it is not accumulated on either side by angles whose rays land on $H_{43/127}$, even though the landing point is not in the boundary of a Fatou component. Therefore erasing the angles states corresponding to it and its forward orbit would change the regular language. The Hubbard tree automaton $H_S(43/127)$ has 42 states, of which 17 are arc states and 25 are angle states. 

\begin{figure}[ht]
    \centering
    \includegraphics[width=0.7\linewidth]{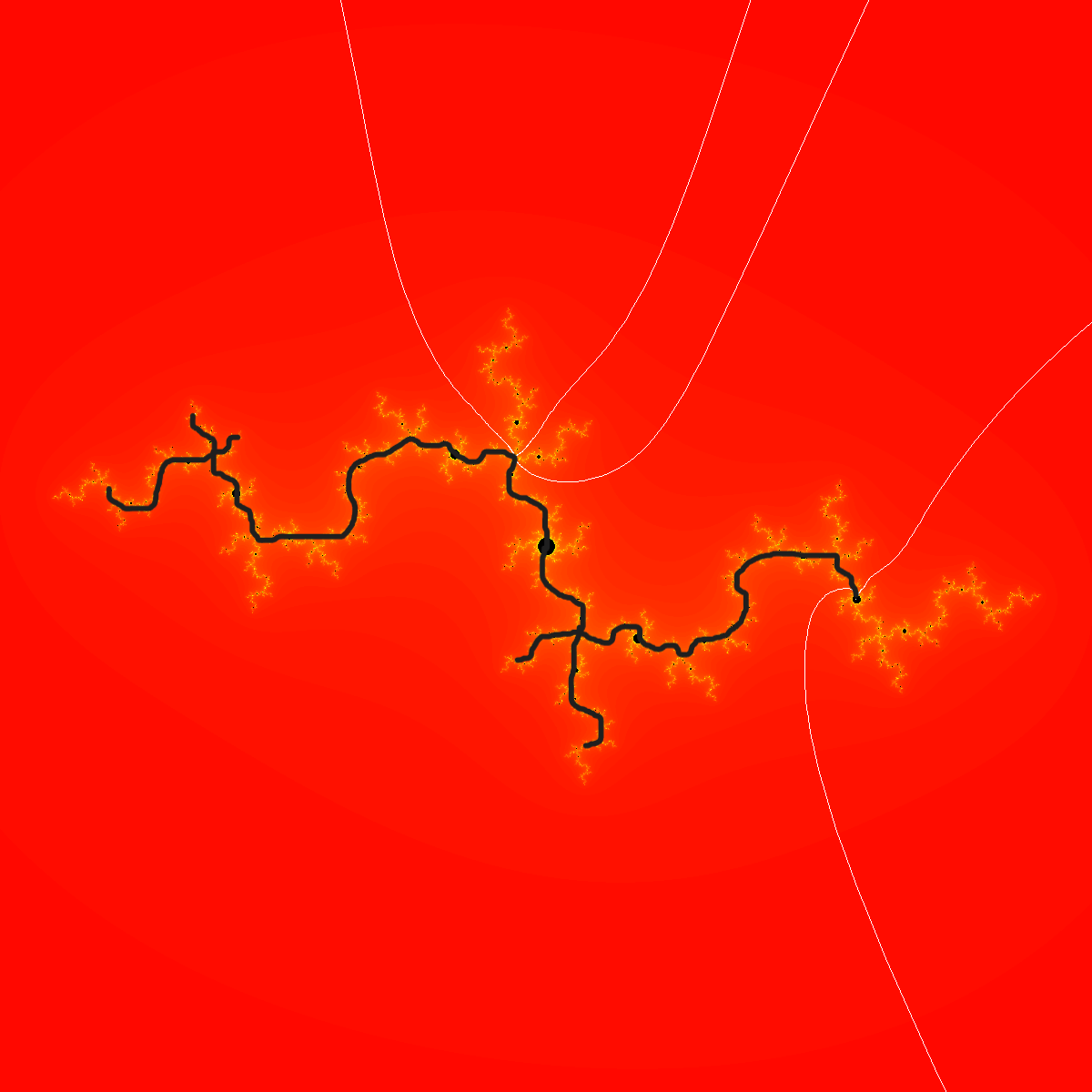}
    \caption{Sketch of the Hubbard tree $H_{43/127}$ and its forbidden region. The angle $101/510$ is isolated landing on the Hubbard tree.}
    \label{fig: 43over127}
\end{figure}
\end{exm}

%\begin{comment}
\begin{exm} \label{example: 3/16}
Let $\theta = 3/16 = .0011$, and consider the parameter $c_{3/16}$; this example will be relevant in Section \ref{section: ray connections between hubbard trees} for exhibiting a topological mating where the Hubbard trees are disjointly embedded. The kneading sequence of $\theta$ is $\nu = \texttt{1101}\overline{\texttt{0}}$, and its angled internal address is
$$
(1)_{1/3} \to (3)_{1/2} \to (5)_{1/2} \to (6)_{1/2} \to (7)_{1/2} \to \cdots
$$
Because $\theta$ is dyadic, the kneading sequence has periodic block $B = \overline{\texttt{0}}$, and  so the sequence of entries in the angled internal address will contain all integers greater than or equal to $3$. Again by (\ref{equation: denominators}), the rotation numbers $p_k/q_k$ are all equal to $1/2$ for $k \geq 1$. The critical orbit falls into a trivial orbit portrait, so no extremity of the Hubbard tree can be an attaching point of the preimage Hubbard tree. Because the valencies are $v \leq 3$ and $H_{3/16} \subseteq J_{3/16}$, no angle whose ray lands on $H_{3/16}$ is isolated, as was the case of Example \ref{example: 43/127}.

The only satellite orbit portrait with valence $v \geq 3$ corresponds to the $\alpha$ fixed point (for the entry $S_0 = 1$), with external arguments $\{1/7, 2/7, 4/7\}$. The preperiodic preimage $-\alpha$ has external arguments $\{1/14, 9/14, 11/14\}$, but it still is a branch point of $H_{3/16}$, since it contains elements of the orbit of $3/16$ in all of the angle intervals determined by adjacent external arguments.

The two preimages $\alpha^2_1$, $\alpha^2_2$ of $-\alpha$ have external arguments
$$
A^2_1 = \{2/28, 23/28, 25/28\}, \quad A^2_2 = \{ 9/28, 11/28, 15/28 \},
$$
and they are grouped this way due to the partition of $S^1$ given by $\{3/32, 19/32\}$. The point $\alpha^2_1$ is not a branch point of $H$, because only two of its angle intervals contain elements in the orbit of $\theta$. Since $(23/28, 25/28)$ does not, it is a forbidden arc. The point $\alpha^2_2$ is still a branch point, having ``lost none of its branches'' under pullback, so we consider its two preimages $\alpha^3_1$, $\alpha^3_2$, with external arguments
$$
A^3_1 = \{9/56, 11/56, 15/56\}, \quad A^3_2 = \{37/56, 39/56, 43/56\}. 
$$
The branches of $H$ lost at these pullbacks correspond to the forbidden arcs $(11/56, 15/56)$ and $(37/56, 39/56)$. The forbidden region is thus
$$
F = \left(\dfrac{11}{56}, \dfrac{15}{56} \right) \sqcup \left( \dfrac{37}{56}, \dfrac{39}{56} \right) \sqcup\left( \dfrac{23}{28}, \dfrac{25}{28} \right).
$$

From the previous observations, we may simplify a Hubbard automaton representing $\mathcal{H}(3/16)$ by erasing the angle states, since no angle landing on the Hubbard tree is isolated. A representation is given below, with $14$ states:
% https://q.uiver.app/#q=WzAsMTQsWzIsMSwiMCJdLFszLDEsIjAiXSxbNCwxLCIwIl0sWzQsMiwiMSJdLFszLDIsIjAiXSxbMiwyLCIwIl0sWzAsMiwiMSJdLFsyLDAsIjEiXSxbMSwxLCIxIl0sWzIsMywiMSJdLFszLDMsIjAiXSxbNCw0LCIxIl0sWzQsNSwiMSJdLFsyLDQsIjAiXSxbMCwxXSxbMSwyXSxbMSw0XSxbMywyXSxbMyw0XSxbMTAsM10sWzksMTBdLFs0LDVdLFs1LDldLFs1LDZdLFs2LDhdLFs2LDcsIiIsMSx7ImN1cnZlIjotM31dLFs3LDFdLFsyLDcsIiIsMSx7ImN1cnZlIjoxfV0sWzEyLDExXSxbMTIsMTIsIiIsMSx7InJhZGl1cyI6MSwiYW5nbGUiOjkwfV0sWzExLDNdLFswLDAsIiIsMSx7InJhZGl1cyI6MX1dLFs4LDBdLFs0LDEwXSxbMTMsMTJdLFs2LDEzXSxbMTMsMTFdXQ==
\[\begin{tikzcd}
	&& 1 && \\
	& 1 & 0 & 0 & 0 \\
	1 && 0 & 0 & 1 \\
	&& 1 & 0 \\
	&& 0 && 1 \\
	&&&& 1
	\arrow[from=1-3, to=2-4]
	\arrow[from=2-2, to=2-3]
	\arrow[from=2-3, to=2-3, loop, in=60, out=120, distance=5mm]
	\arrow[from=2-3, to=2-4]
	\arrow[from=2-4, to=2-5]
	\arrow[from=2-4, to=3-4]
	\arrow[curve={height=6pt}, from=2-5, to=1-3]
	\arrow[curve={height=-18pt}, from=3-1, to=1-3]
	\arrow[from=3-1, to=2-2]
	\arrow[from=3-1, to=5-3]
	\arrow[from=3-3, to=3-1]
	\arrow[from=3-3, to=4-3]
	\arrow[from=3-4, to=3-3]
	\arrow[from=3-4, to=4-4]
	\arrow[from=3-5, to=2-5]
	\arrow[from=3-5, to=3-4]
	\arrow[from=4-3, to=4-4]
	\arrow[from=4-4, to=3-5]
	\arrow[from=5-3, to=5-5]
	\arrow[from=5-3, to=6-5]
	\arrow[from=5-5, to=3-5]
	\arrow[from=6-5, to=5-5]
	\arrow[from=6-5, to=6-5, loop, in=330, out=30, distance=5mm]
\end{tikzcd}\]

Some periodic cycles for the automaton give the angles
$$
\begin{array}{rcl}
    .\overline{001} & = & 1/7, \\
    .\overline{00101} & = & 5/31, \\
    .\overline{000011} & = & 3/64,
\end{array}
$$
and an example of an irrational angle whose external ray lands on $H_{3/16}$ is
$$
.(0)^411(0)^511(0)^6(11)\ldots,
$$
where $(0)^k$ denotes a string of $k$ consecutive zeros. We can also find examples of angles whose external rays do not land on the Hubbard tree: for instance, $.\overline{00111} = 7/31$ is not recognized by the automaton, since every sequence of two consecutive zeros in the graph is followed by at most two consecutive ones.
\end{exm}

\section{Polynomial Mating and Automata} \label{section: mating}

We recall the basic constructions of matings of polynomials: see \cite{Buff2013, Jung2017, Lei1992}. Let $f$ and $g$ be two monic polynomials of degree $d \geq 2$ having connected Julia sets. Let $\CC_f$ and $\CC_g$ be the dynamical planes of $f$ and $g$ respectively. Due to the existence of Böttcher coordinates near infinity, we may attach circles at infinity $S^1_f$ and $S^1_g$ to the dynamical planes so that the polynomials extend continuously as the map $\phi \mapsto d \cdot \phi$ to them. We glue the two topological closed disks $\overline{\CC}_f$ and $\overline{\CC}_g$ obtained this way along their circles at infinity with the identification $(\phi, f) \sim (-\phi, g)$, that is, flipping the orientation of the $g$-plane. This results in a topological $2$-sphere $S^2_{f,g}$ whose equator is identified with the circles at infinity. 

\begin{definition}
The \textit{formal mating} $f \sqcup g : S^2_{f,g} \to S^2_{f,g}$ is the topological branched cover acting as $f$ on the $f$-hemishpere, as $g$ on the $g$-hemisphere, and as $\phi \mapsto d \cdot \phi$ on the equator.
\end{definition}

If $f$ and $g$ are PCF, the resulting object is a postcritically finite topological branched cover of the sphere, for which the ideas of Thurston theory apply \cite{Douady1993}. Suppose now that the Julia sets of $f$ and $g$ are locally connected, so that external rays land on the Julia sets continuously with respect to the angle.

\begin{definition}
A \textit{ray} $\R_{f \sqcup g}(\theta)$ on the formal mating $f \sqcup g :S^2_{f,g} \to S^2_{f,g}$ is the union $\R_{f}(\theta) \cup \{(\theta, f)\} \cup \R_g(-\theta)$, that is, the union of the external ray at angle $\theta$ in the $f$-plane, the external ray at angle $-\theta$ in the $g$-plane, and the point of intersection of their closures at the glued equator. A \textit{ray connection} is a topological embedded arc on $S^2_{f, g}$ consisting of a union of rays, and the \textit{ray equivalence classes} are the equivalence classes obtained by the transitive closure of the relation of two points belonging to the same ray. Each ray equivalence class is a union of rays, possibly infinite. A \textit{direct ray connection} between $X \subseteq J_c$ and $Y \subseteq J_g$ is a ray connection consisting of a single ray whose landing points on $f$-plane and $g$-plane are in $X$ and $Y$, respectively.
\end{definition}

One should think of ``indirect'' ray connections as those which connect $X$ to $Y$, but may pass through points not in either set along a ray connection having more rays. We will mainly be concerned with direct ray connections beween Hubbard trees. Note that $f \sqcup g$ preserves ray equivalence classes.

\begin{definition}
If the quotient of $S^2_{f,g}$ by collapsing ray equivalence classes to points happens to still be homeomorphic to a $2$-sphere, we say that the self-map $f \coprod g$ on the quotient is the \textit{topological mating} of $f$ and $g$. If $f \coprod g$ is topologically conjugate to a rational map $R: \hat{\CC} \to \hat{\CC}$ by a conjugacy which is holomorphic on the interiors of (the images in the quotient of) the filled Julia sets $K_f$ and $K_g$, we say that $R$ is a \textit{geometric mating} of $f$ and $g$.
\end{definition}

\begin{remark} \label{remark: quadratic PCF matings}
Tan Lei \cite{Lei1992} and Shishikura \cite{Shishikura2000} showed that the geometric mating of two PCF quadratic polynomials exists if and only if they are not in conjugate limbs of the Mandelbrot set. This uses Thurston's characterization of rational maps in terms of pullbacks of multicurves, and ideas of M. Rees. For most of this paper, we will be concerned only with formal matings of quadratic polynomials and the ray connections thereof.
\end{remark}

Let $c_1$ and $c_2$ be the parameters for $f_{1} \coloneqq f_{c_1}$ and $f_{2} \coloneqq f_{c_2}$ two quadratic PCF polynomials. We introduce the notion of mating Hubbard tree automata, which will allows us to find all the angles $\phi$ such that $\R_{f_1}(\phi)$ lands on $H_1$ in the dynamical plane of $c_{1}$ and $\R_{f_2}(-\phi)$ lands on $H_2$ in the dynamical plane of $c_{2}$; that is, $\R_{f_1 \sqcup f_2}(\phi)$ is a direct ray connection between the Hubbard trees.

\begin{definition}
A \textit{mating automaton} for $f_{1} \sqcup f_{2}$ is any automaton whose regular language is the set of direct ray connections between $H_{c_1}$ and $H_{c_2}$.
\end{definition}

We provide two ways of constructing a mating automaton. The first is based on the product automaton construction, through the following lemma:

\begin{lemma} \label{lemma: flip bits}
Let $\mathcal{H}(c)$ be a Hubbard tree automaton, and $\overline{\mathcal{H}(c)}$ be the automaton obtained by flipping the bit associated to each state, that is, a state that has a bit value $0$ becomes a states with bit value $1$, and vice-versa. Then $\overline{\mathcal{H}(c)}$ is a Hubbard tree automaton associated to the complex conjugate $\overline{c}$.
\end{lemma}

\begin{proof}
This is immediate from the fact that $f_{\overline{c}}$ is anti-conformally conjugate to $f_c$, and the landing point of $\R_{c}(\phi)$ corresponds to the landing point of $\R_{\overline{c}}(-\phi)$. The only ambiguity is with attributing bit values to the marked angles $0$ and $1/2$, but as observed in the proof of Proposition \ref{proposition: regular language lands on tree}, as long as we attribute distinct values to them, the regular language corresponds to the angles whose external rays land on the Hubbard tree.
\end{proof}

\begin{definition} \label{definition: product automaton}
Given two automata $\mathcal{A}_1 = (\mathcal{G}_1, \chi_1)$ and $\mathcal{A}_2 = (\mathcal{G}_2, \chi_2)$, seen as $\{0,1\}$-labeled directed graphs, the \textit{product automaton} $\mathcal{A}_1 \times \mathcal{A}_2$ is the automaton such that:
\begin{itemize}
    \item The set of states is the subset $\mathcal{S} \subseteq \mathcal{S}_1 \times \mathcal{S}_2$ of those pairs $(A_1, A_2)$ such that $\chi_1(A_1) = \chi_2(A_2)$;
    \item There is a transition $(A_1, A_2) \to (B_1, B_2)$ if and only if there are transitions $A_1 \to B_1$ and $A_2 \to B_2$.
\end{itemize}
\end{definition}

The following is immediate from the definitions.

\begin{prop} \label{proposition: product automaton}
The product automaton $\mathcal{H}(c_1) \times \overline{\mathcal{H}(c_2)}$ is a mating automaton for $f_{1} \sqcup f_{2}$.
\end{prop}

\begin{remark}
As we are interested only in the infinite allowable paths within the automata, we may again apply the previously described trimming procedure of iteratively erasing states that do not have any transitions from them, so as to obtain a simpler mating automaton having the same regular language. Letting $G$ be the transition matrix representing this automaton, where the entry $g_{ij}$ is $1$ if $j$ transitions to $i$ and $0$ otherwise, this procedure corresponds to erasing the $j$-th column and the $j$-th row whenever the $j$-th column consists of all zeros.
\end{remark}

While the product automaton construction is more general, applicable to any pair of automata, the second construction of a mating automaton takes advantage of the fact that the Hubbard tree automata come from two Markov partitions of $S^1$. The following is immediate from Proposition \ref{proposition: forbidden region and external arguments}:

\begin{lemma} \label{lemma: union automaton}
If $F_1$ and $F_2$ are the forbidden regions for $H_{c_1}$ and $H_{c_2}$ respectively, $\R_{f_1 \sqcup f_2}(\phi)$ is a direct ray connection between the Hubbard trees if and only if $\phi$ never enters $F_1 \cup (-F_2)$.
\end{lemma}

\begin{definition} \label{definition: union automaton}
Given sets of marked angles $S_1$ and $S_2$ for $f_1$ and $f_2$ respectively, let $S = S_1 \cup (-S_2)$ be an enlargement of $S_1$. The \textit{union automaton} $\mathcal{M}(c_1,c_2)$ of $c_1$ and $c_2$ is the automaton constructed from $S$ as in Section \ref{section: hubbard tree automaton}, but for which the forbbiden states are those contained in $F_1 \cup (-F_2)$. Its set of states $\mathcal{S}_{\mathcal{M}}$ after trimming are called the \textit{$\mathcal{M}$-states}.
\end{definition}

Lemma \ref{lemma: union automaton} guarantees that the union automaton $\mathcal{M}(c_1,c_2)$ is a mating automaton for $f_{c_1} \sqcup f_{c_2}$, and its states form a subset of the good states $\mathcal{S}_{\mathrm{good}}$ for $\mathcal{H}_S(c_1)$. When there is no risk of confusion, we may denote the union automaton simply by $\mathcal{M}$, or any mating automaton for $f_{1} \sqcup f_{2}$ by $\mathcal{M}$ when the specific construction is irrelevant.

This construction has some advantages: $\mathcal{M}$ can be seen as a sub-automaton of $\mathcal{H}_S(c_1)$, with some of the states and transitions erased, and so the notions of arc states and angle states still make sense for $\mathcal{M}$. Moreover, number of states of $\mathcal{M}$ is $O(|S_1| + |S_2|)$, instead of $O(|S_1| \cdot |S_2|)$. Note that embedded in this definition is the preferential viewpoint of taking $\mathcal{H}(c_1)$ as a reference, instead of $\mathcal{H}(c_2)$. The analogous construction of $\mathcal{M}$ as a subautomaton of $\mathcal{H}_{-S}(c_2)$ is possible.

\section{Laminations and Veins} \label{section: laminations}

In aiming to understand the combinatorics and dynamics of degree $d$ polynomials on $\CC$, Thurston introduced the notion of degree-$d$ invariant laminations. We write here some definitions and results which will be relevant to us in the quadratic case, and defer others to Thurston's original manuscript in \cite{Thurston2022} and the subsequent additions and appendix in \cite{Thurston2009}. Essentially, degree $d$ invariant laminations serve as a combinatorial model for locally connected Julia sets of polynomials.

\subsection{Quadratic invariant laminations} \label{subsection: quadratic invariant laminations}

A \textit{leaf} $l$ on the closed unit disk $\overline{\DD}$ is, formally, an unordered pair of angles on $S^1$, regarded as its \textit{endpoints}. It is identified with the geodesic between these ideal points on the boundary $S^1 = \partial \overline{\DD}$ with respect to the hyperbolic metric on $\overline{\DD}$. (In the Poincaré model, this is an arc of a circle meeting $S^1$ orthogonally, and in the Klein model, this is a segment of a straight line. The choice of which model to use is immaterial.) Its \textit{length} $|l|$ is the length on $S^1$ of the shortest arc determined by the endpoints, measured with proportion to the total length of $S^1$. The leaf is \textit{degenerate} if $|l| = 0$, or equivalently the endpoints coincide. A \textit{(geodesic) lamination} $L$ on $\overline{\DD}$ is a closed collection of leaves which are disjoint, except possibly at their common endpoints on $S^1$. A \textit{gap} of $L$ is the closure of a component of $\overline{\DD} \setminus \bigcup L$. 

\begin{definition} \label{definition: quadratic invariant lamination}
A \textit{degree $d$ invariant lamination} $L$, following Thurston, is a closed collection of disjoint leaves such that:
\begin{itemize}
    \item (Forward invariance) If $l = \overline{pq}$ is a leaf of $L$, then $f(l) \coloneqq \overline{f(p)f(q)}$ is a leaf of $L$, where $p, q \in S^1$ are the endpoints and $f: S^1 \to S^1$ is the map $f(z) = z^d$.

    \item (Backward invariance) If $l = \overline{pq}$ is a leaf of $L$ then there exists $d$ disjoint leaves $L$, each joining a preimage of $p$ to a preimage of $q$;

    \item (Gap invariance) For any gap $G$ of $L$, the convex hull of the image of $G \cap S^1$ is either a point on $S^1$, a leaf in $L$, or a gap of $L$. If it is a gap, the boundary of $G$ maps locally monotonically onto the boundary of the image gap with positive orientation.
\end{itemize}
\end{definition}

We will henceforth restrict ourselves to the quadratic case ($d = 2$). Given a quadratic invariant lamination $L$, its \textit{major leaves} are the leaves of greatest length. Either there is a unique one, and it is a diameter, or there are exactly two antipodal major leaves. The \textit{minor leaf} is the common image of the major leaves, being a degenerate leaf (a point on $S^1$) in the former case. Intuitively, the minor leaf represents the critical value in the quadratic polynomial that the lamination would represent. The lamination is \textit{clean} if no more than two leaves meet at a common endpoint on $S^1$.

\begin{definition}
A leaf $l$ of a quadratic invariant lamination $L$ \textit{characteristic} if it is the minor leaf, or it separates the minor leaf from all of the forward images of $l$.
\end{definition}

Note that a diameter can never be characteristic.

\begin{lemma} \label{lemma: characteristic leaves equivalent description}
A leaf $l$ is characteristic if and only if all of its forward iterates have length greater than or equal to $|l|$, and the minor leaf $m$ is contained in the shortest arc of $S^1$ determined by the endpoints of $l$.
\end{lemma}

\begin{proof}
The equivalence is immediate if $l$ is the minor leaf.
If $l$ satisfies the second condition in the statement, it cannot be a diameter: its image would be a degenerate leaf and have length $0$.

Suppose that $l$ is characteristic. There are two distinct ways of connecting the preimages of the endpoints of $l$ to form leaves that map to it; we claim that the preimages $l_1$ and $l_2$ of $l$ in the lamination must be such that their lengths are $(1 - |l|)/2$. Otherwise their lengths are $|l|/2$, and because the major leaves cannot cross $l_1$ and $l_2$, preservation of orientation under doubling would imply that $m$ is not contained in the shortest arc determined by the endpoints of $l$.

Suppose that some forward image $f^k(l)$ has length less then $|l|$, for $k \geq 1$. Its preimage $f^{k-1}(l)$ either has half the length, or the length must be $(1+|f^k(l)|)/2$. Because the length cannot keep halving under preimages as it must be $|l|$ eventually, there is some $k' < k$ such that the length is greater than $(1 + |l|)/2$, and therefore be contained in the region of $\overline{\DD}$ between $l_1$ and $l_2$ and connect the opposite arcs of $S^1$. But because the doubling map preserves orientation, the image of this leaf won't be separated from the minor leaf by $l$.
\end{proof}

In particular, the length of a characteristic leaf is always $|l| \leq 1/3$, and its preimages have length $|l_1| = |l_2| = (1 - |l|)/2 \geq 1/3$. Otherwise $f(l)$ would have length less than $1/3$, contradicting the equivalent condition for $l$ to be characteristic.

\begin{definition}
The \textit{characteristic lamination} of a quadratic invariant lamination $L$ is its set $\Char(L)$ of characteristic leaves.
\end{definition}

We call the set of characteristic leaves of a quadratic invariant lamination $L$ the \textit{characteristic lamination} of $L$, denoted by $\Char(L)$.

\begin{lemma} \label{characteristic leaves closed}
$\Char(L)$ is a closed sublamination of $L$.
\end{lemma}

\begin{proof}
Let $(l_n)_{n \geq 0}$ be a sequence of characteristic leaves converging to a leaf $l$. By hypothesis, for all $k$ and $n$ we have $|f^k(l_n)| \geq |l_n|$. The size of a leaf is continuous with respect to endpoints, and so for any fixed $k \geq 0$ we have $|f^k(l)| \geq |l|$. Moreover, because $m$ is contained in the minor arcs determined by all of the leaves $l_n$, it will also be contained in the minor arc determined bt $l$, since all characteristic leaves (and in the limit, also $l$) have length not greater than $1/3$.
\end{proof}

Let $m$ be a (possibly degenerate) leaf on $\overline{\DD}$ such that:
\begin{itemize}
    \item All of its forward images are pairwise unlinked;
    \item No forward image is shorter than $m$;
    \item If it is non-degenerate, the forward orbit of $m$ does not cross the two preimages of $m$ which have length greater than or equal to $1/3$.
\end{itemize}
From $m$, Thurston constructs a quadratic invariant lamination $I(m)$ having $m$ as its minor \cite{Thurston2009}. Under some suitable hypotheses on $m$, it may have the properties of being clean, unique, or minimal having $m$ as its minor; these depend on $m$ being degenerate or not, and on having preperiodic endpoints or not. Details will be recalled as convenient to our purposes. Let $\Char(m) \coloneqq \Char(I(m))$ be the characteristic lamination of $I(m)$. Note that every characteristic leaf $l$ of a lamination $L$ satisfies the required properties for it to generate a lamination $I(l)$, and is therefore the minor of some lamination and of QML.

\begin{comment}
\begin{lemma}
Let $l$ be a characteristic leaf for $I(m)$. If $m$ is non-denegerate, or if $m$ is degenerate and not periodic, then $I(l) \subseteq I(m)$ and $C(l) \subseteq C(m)$.
\end{lemma}

\begin{proof}
By Propositions II.6.5 and II.6.6 in \cite{S5}, every quadratic invariant lamination having $m$ as a minor leaf contains $I(m)$. Given $I(l)$, adjoin to it the leaf $m$, all of its forward images, and take the closure. This lamination is evidently contained in $I(m)$ and has $m$ as a minor, so by the minimality property also contains $I(m)$. In particular, $I(l) \subseteq I(m)$. 
\end{proof}
\end{comment}

The quadratic minor lamination has a natural tree-like structure, with an ordering $l \prec m$ of its leaves if $l$ separates $m$ from $1 \in S^1$. This gives another description of $C(m)$, and leads to a monotonicity property of characteristic laminations (see also \cite{Jung2014}):

\begin{prop} \label{proposition: characteristic coincides}
$\Char(m)$ coincides with the set of leaves in QML that precede $m$.
\end{prop}

\begin{proof}
If $l \preceq m$, Proposition II.6.10a in \cite{Thurston2009} guarantees that $l$ is a leaf of $I(m)$. As $l$ is the minor of some lamination, no forward image of $l$ is shorter than $l$, and so we only need to show that $m$ is contained in the minor arc of $l$. If not, then $1 \in S^1$ is contained in the minor arc of $l$, But this cannot occur for the minor of any quadratic invariant lamination; the preimages of length $\geq 1/3$ would cross $l$ itself. Therefore $l$ is characteristic for $I(m)$. Conversely, if $l \in \Char(m)$, it is the minor of some lamination, and therefore is a leaf of QML. Either $l = m$ or it separates $m$ from $1 \in S^1$, since neither contains $1$ in its minor arc, and so $l \preceq m$.
\end{proof}

\begin{cor} \label{corollary: monotonicity of characteristic leaves}
If $l$ and $m$ are leaves of QML, then $l \preceq m$ if and only if $\Char(l) \subseteq \Char(m)$.
\end{cor}

\begin{proof}
The forward implication is immediate from the fact that the leaves that precede $m$ in QML are linearly ordered. Conversely, if $\Char(l) \subseteq \Char(m)$, in particular $l \in I(m)$, and previous arguments guarantee that $l \preceq m$ in QML.
\end{proof}

We return to quadratic polynomials. Let $c \in M$ be a parameter for which the Julia set $J_c$ is locally connected. There is a canonical quadratic invariant lamination $L_c$ associated to $c$, which records which external rays land together on $J_c$, and serves as a combinatorial model for $J_c$. Two angles $\theta$ and $\theta'$ are related if their external rays land together, and this defines a \textit{laminational equivalence relation} \cite{Blokh2008}. From it, we construct the lamination $L_{c}$ whose leaves are the sides of the convex hulls of equivalence classes. The leaves represent pairs of adjacent external arguments of the same point on $J_c$, and gaps are either polygonal, corresponding to preperiodic or precritical points, or correspond to Fatou components. This lamination is always clean, and Proposition II.6.7 in \cite{Thurston2009} guarantees that it coincides with $I(m)$ where $m$ is its minor leaf. For $c$ a parameter such that $J_c$ is locally connected, we define $\Char(c) \coloneqq \Char(L_c)$.

The equivalent notion to a characteristic leaf in the polynomial setting is that of a characteristic point. The Lemma below serves to define them. Recall that $c$ is the critical value for $f_c$, and the interval notation denotes the unique regular arc in $K_c$ connecting the endpoints.

\begin{prop} \label{proposition: characteristic points}
Let $c$ be a parameter for which $J_c$ is locally connected. The following are equivalent for a non-precritical point $z \in J_{c}$:

\begin{itemize}
    \item[(1)] For every $k \geq 0$, $f^k(z) \notin (z, c]$;
    \item[(2)] For every $k \geq 0$, $z \in [f^k(z), c]$;

    \item[(3)] If $l$ is the leaf in the lamination $L_c$ corresponding to the ray pair $P(\theta, \theta')$ landing on $z$ and adjacent to $[z, c]$, then all forward images $f^k(l)$ have length greater than or equal to that of $l$, and the arc $(\theta, \theta')$ has length less than $1/2$;

    \item[(4)] If $l$ is the leaf in the lamination $L_c$ corresponding to the ray pair $P(\theta, \theta')$ landing on $z$ and adjacent to $[z, c]$, then $l$ is characteristic.
\end{itemize}
\end{prop}

\begin{proof}
By Lemma \ref{lemma: characteristic leaves equivalent description}, (3) and (4) are equivalent. Conditions (1) and (3) are equivalent under the laminational interpretation of the Julia sets, and likewise conditions (2) and (4).
\end{proof}

If $z$ satisfies any of the conditions above, we say that $z$ is \textit{characteristic}. If $z = c$ and it has only one external argument, we consider the leaf of conditions (3) and (4) to be the degenerate one; we emphasize that under this definition, the points $\beta$ and $c$ are always characteristic.

Every periodic orbit on $J_c$ with valence $v \geq 2$ has a unique characteristic point, but if an orbit is not periodic, even if it is contained in $H_c$, it may not contain any characteristic points. On the other hand, if an orbit is strictly preperiodic, it may contain more than one characteristic point. Since $\Char(L)$ is closed and the set of characteristic points is a closed quotient of the characteristic lamination, it is also closed, and is a subset of the arc $[\beta, c]$.

\subsection{Veins}

One of the ways of understanding the combinatorics of the Mandelbrot set is through the notion of veins. One may define them in purely combinatorial terms (see \cite{Douady1984}, where ``combinatorial veins'' are specific collections of postcritically finite parameters), but we instead choose to treat them as topological arcs within $M$, and to this end we compile some known results that relate to our purposes \cite{Branner1988, Jung2014, Schleicher1999, Riedl2000}.

An embedded arc within the Mandelbrot set is \textit{regular} if it passes through hyperbolic components only along internal rays, given their uniformizations by the multiplier of the corresponding attracting periodic cycle. The result below is Theorem 5.6 in \cite{Schleicher1999} combined with Proposition 4.8 in \cite{Jung2014}:

\begin{thm} \label{theorem: veins}
Given a postcritically finite parameter $c \in M$, there exists a regular arc $V_c$ in the Mandelbrot set connecting $1/4$ to it. Moreover:

\begin{enumerate}
    \item for every $c' \in V_c$:
    \begin{enumerate}
        \item $c'$ is at most finitely renormalizable or is in the image of $M \cap \RR$ by a tuning map;
        \item the Julia set $J_{c'}$ is locally connected;
        \item if $f_c'$ is not geometrically finite, then $c'$ is in the impression of exactly two parameter rays $\R_M(\theta)$ and $\R_M(\theta')$, and the angles form a dynamical ray pair $P_{c'}(\theta, \theta')$ landing together at the critical value $c'$ in $K_{c'}$;
        \item let $T_{c'}$ be the regular hull of the (possibly infinite) postcritical orbit or $f_{c'}$. Then $T_{c'}$ is a finite topological tree and compact;
    \end{enumerate}
    \item hyperbolic parameters are dense in $V_c$.
\end{enumerate}
\end{thm}

\begin{definition} \label{definition: vein}
The arc $V_c$ constructed in Theorem \ref{theorem: veins} is called the \textit{vein} of (or to) $c$.
\end{definition}

In fact, the theorem constructs veins for all parameters with trivial fiber, or on a tuning image of the real line; but the tree $T_c$ (which we call the \textit{generalized Hubbard tree} of $f_c$) may fail to be finite and compact. See for instance Example 4.3 in \cite{Jung2017}. When $c' \in V_c$ is not geometrically finite, the ray pair in item (1)(c) forms the minor leaf of its lamination $L_{c'}$.

The proof strategy for the existence of veins, first described by Kahn in the case of dyadic Misiurewicz parameters, starts by constructing a locally connected model $M^*$ for the Mandelbrot set, and a projection map $\pi_M : M \to M^*$ that is monotone, meaning that the preimage of connected sets is connected \cite[Theorem A.6]{Thurston2009}. One possibility is taking $M^* = \overline{\DD}/\mathrm{QML}$, the quotient of the unit disk by collapsing leaves of QML to points. Then, choosing a suitable arc $\varphi: [0,1] \to M^*$ in the model space, one construct a continuous lift $\psi: [0,1] \to M$, whereby continuity depends in an essential manner on Yoccoz's proof of triviality of fibers for finitely renormalizable parameters \cite{Hubbard1993}, and that only finitely many sublimbs having rotation number different from $1/2$ may appear in the internal address of a postcritically finite parameter (Propsotion \ref{proposition: denominator not 1/2}).

Results of Lyubich \cite{Lyubich1997}, Graczyk and Swiątek \cite{Graczyk1997}, and Levin and van Strien \cite{Levin1998} show local connectivity of the Julia sets and density of hyperbolicity along veins, when coupled with the fact that tuning preserves triviality of fibers in $M$ \cite[Corollary 3.2]{Schleicher1999} and triviality of fibers of the Julia sets \cite[Theorem 4.4]{Schleicher1999}. Moreover, density of hyperbolicity along veins implies the following:

\begin{prop} \label{proposition: uniqueness of veins}
For all postcritically finite parameters, there exists a unique regular arc within $M$ connecting $1/4$ to it up to reparametrization, and it is the vein $V_c$ as constructed in Theorem \ref{theorem: veins}.
\end{prop}

\begin{proof}
Let $\psi$ be a vein to $c$ satisfying the properties of Theorem \ref{theorem: veins}, and $\eta$ another regular arc connecting $1/4$ to $c$ within $M$. If their embedded images are not equal, the difference corresponds to Jordan curves in $M$: whenever the arcs ``diverge'' in their paths, they must meet again at some first point of intersection before reaching $c \in M$, and this forms a loop.

Since $M$ is full, the regions bounded by these curves are contained in the interior of $M$. They must correspond to hyperbolic components by density of hyperbolicity along $\psi$, and the assumption that $\psi$ and $\eta$ are regular arcs contradicts the assumption that their images are different.
\end{proof}

The same arguments show that the collection of all veins forms a topological tree in $M$, which is branched only at centers of hyperbolic components and at Misiurewicz points \cite[Lemma 5.5]{Schleicher1999}. Given any parameter $c$ in the tree of veins, we may also define its vein $V_c$ as the regular arc within the tree from $1/4$ to $c$. For two parameters $c, c'$ in the topological tree of veins, we say that $c$ \textit{precedes} $c'$, and denote this by $c \preceq c'$, if $c$ is contained in the vein to $c \in V_{c'}$. As a consequence of the previous results, $\preceq$ is a partial order.

We recall the classical form of the parameter-dynamical correspondence for periodic characteristic leaves, due to Lavaurs \cite{Lavaurs1986, Schleicher2000}:

\begin{lemma}
If $z \in J_c$ is a periodic characteristic point, then there exists a unique parameter $c_z$ on the combinatorial vein of $c$ such that if $P_c(\theta, \theta')$ is the characteristic ray pair landing on $z$, then $P_M(\theta, \theta')$ is a parameter ray pair landing at $c_z$. Conversely, if $c'$ is a Misiurewicz parameter or the root of a hyperbolic component and $P_M(\theta, \theta')$ is its characteristic parameter ray pair with respect to $c$, then there exists a preperiodic or periodic characteristic point $z_{c'} \in J_{c'}$, respectively, such that $P(\theta, \theta')$ is the characteristic ray pair landing on $z_{c'}$.
\end{lemma}

These characteristic periodic points on the Julia set are linearly ordered as a subset of the dynamical regular arc $[\beta, c]$, and this order is preserved in the combinatorial vein. Moreover, the next characteristic periodic point is the $\alpha$ fixed point, and there are no characteristic points (preperiodic or otherwise) in $(\beta, \alpha)$. This is because every biaccessible point has an iterate on the arc $[\alpha, c]$ \cite{Lei1992}.

The main relevance of veins to us is that the characteristic lamination of the Julia set lamination varies monotonically with respect to parameters along a vein, as a direct consequence of Corollary \ref{corollary: monotonicity of characteristic leaves} and the fact that $L_c = I(m)$ where $m$ is the minor leaf of the lamination.

\begin{lemma} \label{lemma: preperiodic counterpart on vein}
Let $c$ be a parameter in the tree of veins, and $z$ a characteristic preperiodic point on $J_c$. If $c_z$ is the parameter counterpart of $c$ as above, then $c_z \in V_c$.
\end{lemma}

\begin{proof}
The parameter ray pair $P_M(\theta, \theta')$ associated to $c_z$ must separate $c$ and $1/4$, since $c$ is in its wake by realizing the orbit portrait. Because $V_c$ is strictly contained in $M$, it must cross $P_M(\theta, \theta')$ exactly at the landing point $c_z \in M$.
\end{proof}

\begin{prop} \label{proposition: monotonicity along veins}
$c \preceq c'$ implies $\Char(c) \subseteq \Char(c')$.
\end{prop}

\begin{proof}
We will use the fact that gaps are dense in QML \cite[Theorem II.6.13]{Thurston2009}. Assume $c \neq c'$, and let $m$ and $m'$ be the minors of $L_c$ and $L_{c'}$, respectively. By Corollary \ref{corollary: monotonicity of characteristic leaves}, it's sufficient to show that $m \prec m'$ in QML. If that's not the case, then by the branch theorem for QML \cite[Corollary II.6.13b]{Thurston2009}, either $m' \prec m$, or they are separated by distinct sides of a gap corresponding to a hyperbolic component or a Misiurewicz parameter. 

Suppose $m' \prec m$. If these leaves are sides of the same gap, then the parameters either correspond to the same Misiurewicz point, or to a sublimb $W_{p/q}$ of a hyperbolic component $W$. In either case the assumptions that $c \in V_{c'}$ and $c \neq c'$ exclude the possibility that $m' \prec m$. If $m$ and $m'$ are not sides of the same gap, by density of gaps and their classification there must be a preperiodic leaf $l$ that separates $m'$ and $m$, and by Proposition \ref{proposition: characteristic coincides}, $l \in \Char(m)$ but $l \notin \Char(m')$. If $P_M(\theta, \theta')$ is the parameter ray pair corresponding to the leaf $l$, then $c$ is in the wake of $l$, but $c'$ isn't. Therefore the regular arc $(c,c']$ must cross the ray pair $P_M(\theta, \theta')$ to leave the corresponding wake, but it can only do so at $c_l \in M$ and already $c_l \in V_c$, a contradiction with the arc $V_{c'}$ being embedded.

In the latter case, we may assume that $m$ and $m'$ are not sides of the gap that separates them by the arguments above. Therefore there exist periodic leaves $l$ and $n$ such that $l \in \Char(m)$, $l \notin \Char(m')$, and $n \in \Char(m')$, $n \notin \Char(m)$. If $c_l$ and $c_n$ are the corresponding parameters in $M$, $c_l \in V_c \setminus V_{c'}$ and $c_n \in V_{c'} \setminus V_c$, a contradiction with the hypothesis that $V_c \subseteq V_{c'}$.
\end{proof}

A \textit{PCF vein} is the vein of a postcritically finite parameter $c$. A parameter $c$ is the \textit{tip of a (PCF) vein} if it is a Misiurewicz parameter and it has only one external argument $\theta$ (either in the parameter or dynamical plane). This is equivalent to the fact that, if $\theta = \theta_1$ and $\{\theta_m, \theta_{m+1}, \ldots, \theta_{m+k}\}$ is the periodic orbit that $\theta$ iterates into, it corresponds to either the zero portrait, or to a trivial portrait. The former case is that of Misiurewicz dyadic parameters, and an example of the latter is $c_{1/6}$, the landing point of the parameter ray $\R_M(1/6)$. Tips of PCF veins are maximal elements in the partial order described on the tree of veins.

The following Proposition complements Theorem \ref{theorem: veins} (1)(c).

\begin{prop} \label{proposition: irrational point on vein}
Suppose $t$ is the tip of a PCF vein and $l = (\phi, \phi')$ is a characteristic leaf in $L_t$ at an irrational angle. Then there exists a unique $c \in V_t$ such that the dynamical ray pair $P_c(\phi, \phi')$ lands at the critical value $c$ in the $c$-plane.
\end{prop}

\begin{proof}
Due to our hypotheses, $l$ is a leaf of QML. We use the fact that the leaves of QML having periodic endpoints are dense \cite[Corollary II.6.13a]{Thurston2009}, and so there is a sequence $l_k = (\phi_k, \phi_k')$ of leaves having periodic endpoints such that $l_k \to l$. In fact, since $l$ is not the side of a gap, we may assume that the leaves $l_k$ precede $l$, accumulating on it from ``below'' (See the discussion preceding Theorem II.6.13 in \cite{Thurston2009}). Because the leaves $l_k$ separate $t$ from $1$, they correspond to characteristic periodic points on $[\beta_t, t]$ for $f_t$, and also to parabolic parameters $c_k$ which contain the parameter $t$ in their wakes. In particular, the parameters $c_k$ are in the vein $V_t$. The sequence of $c_k$ is strictly increasing on $V_t$, and we let $c$ be the limit of the $c_k$.

The lamination $L_c$ contains all the leaves $l_k$, and because it closed, it also contains $l$. Hence $\R_c(\phi)$ and $\R_c(\phi')$ land together at some point $z \in [\alpha_c, c] \subset T_c$. If $z \neq c$, then $l$ is not the minor leaf $m$ of $L_c$. If there is no periodic leaf separating $l$ and $m$, then both leaves are sides of the same gap in QML, a contradiction with $l$ having irrational endpoints. Let $p = (\theta, \theta')$ be a leaf having periodic endpoints separating $l$ and $m$, so that the parameter ray pair $P_M(\theta, \theta')$ lands at a parabolic parameter $c^* \succ c$ on $V_t$. The leaf $m$ is characteristic for $L_c$, and by Proposition \ref{proposition: monotonicity along veins}, $m$ must be present in the lamination $L_{c^*}$. However, $p$ is the minor leaf for $L_{c^*}$ and precedes $m$ in QML, so $m$ cannot be in $\Char(c^*)$. Therefore $z = c$.

Suppose that there were another $c' \in V_t$ such that its minor leaf is $l$ and without loss of generality $c \prec c'$. By density of hyperbolicity on $V_t$ there must be a parabolic parameter $c^* \in V_t$ such that $c \prec c^* \prec c'$, having characteristic leaf $p$ with periodic endpoints. Then $l \prec p \prec l$, a contradiction.
\end{proof}

\section{Ray Connections Between Hubbard Trees} \label{section: ray connections between hubbard trees}

Given two postcritically finite parameters $c_1$ and $c_2$, the question of which ray connections $\R_{f_{c_1} \sqcup f_{c_2}}(\phi)$ land on $H_{c_1}$ and on $H_{c_2}$ in the formal mating is essentially that of finding the intersection of two Cantor sets, which are, in this case, defined by excluding iterated preimages of certain open arcs of $S^1$. This intersection may be empty, finite, countable, another Cantor set, or a Cantor set with countably many isolated points. In the language of Section \ref{section: mating}, this will correspond to the regular language of a mating automaton $\mathcal{M}(c_1, c_2)$.

\begin{exm} \label{example: 1/6 with 3/16}
Let $\theta_1 = 1/6$ and $\theta_2 = 3/16$, and consider the formal mating $f_{1/6} \sqcup f_{3/16}$. We show that there are no direct ray connections between the Hubbard trees. As a consequence of Theorem \ref{theorem: first ray connection} and monotonicity of the characteristic lamination, the same will hold true for all pairs of parameters along their respective veins.

To see this concretely, consider below the automaton $\mathcal{H}(1/6)$ from Example \ref{example: 1/6} on the left, and the automaton $\overline{\mathcal{H}(3/16)}$ from Example \ref{example: 3/16} on the right, but with its bits flipped.
% https://q.uiver.app/#q=WzAsMTksWzYsMCwiMSJdLFs3LDAsIjEiXSxbNywxLCIxIl0sWzYsMSwiMSJdLFs1LDAsIjAiXSxbNiwyLCIwIl0sWzcsMiwiMSJdLFs4LDMsIjAiXSxbMCwxLCIwIl0sWzEsMSwiMCJdLFsyLDEsIjAiXSxbMiwyLCIxIl0sWzEsMiwiMSJdLFs0LDMsIjEiXSxbNSwyLCIwIl0sWzgsMiwiMCJdLFs4LDEsIjEiXSxbOCwwLCIwIl0sWzQsMCwiMCJdLFswLDFdLFsxLDJdLFs1LDZdLFsyLDNdLFszLDVdLFswLDAsIiIsMSx7InJhZGl1cyI6MSwiYW5nbGUiOjE4MCwic3R5bGUiOnsidGFpbCI6eyJuYW1lIjoiYXJyb3doZWFkIn0sImhlYWQiOnsibmFtZSI6Im5vbmUifX19XSxbNCwwXSxbMiw2XSxbOCw5XSxbOSwxMF0sWzEwLDExLCIiLDAseyJjdXJ2ZSI6LTF9XSxbMTEsOV0sWzExLDEyXSxbMTIsOF0sWzksMTJdLFsxMSwxMCwiIiwwLHsiY3VydmUiOi0xfV0sWzEzLDcsIiIsMCx7ImN1cnZlIjoxfV0sWzEzLDE0XSxbMTQsN10sWzE0LDE0LCIiLDAseyJyYWRpdXMiOjF9XSxbNiwxNV0sWzE1LDJdLFsxNSwxNl0sWzcsMTVdLFsxNywxXSxbMTYsMTddLFsxOCwxM10sWzE4LDE3LCIiLDAseyJjdXJ2ZSI6LTR9XSxbMywxOF0sWzE4LDRdLFsxLDE2XV0=
\[\begin{tikzcd}
	&&&& 0 & 0 & 1 & 1 & 0 \\
	0 & 0 & 0 &&&& 1 & 1 & 1 \\
	& 1 & 1 &&& 0 & 0 & 1 & 0 \\
	&&&& 1 &&&& 0
	\arrow[from=1-5, to=1-6]
	\arrow[curve={height=-24pt}, from=1-5, to=1-9]
	\arrow[from=1-5, to=4-5]
	\arrow[from=1-6, to=1-7]
	\arrow[tail reversed, no head, from=1-7, to=1-7, loop, in=240, out=300, distance=5mm]
	\arrow[from=1-7, to=1-8]
	\arrow[from=1-8, to=2-8]
	\arrow[from=1-8, to=2-9]
	\arrow[from=1-9, to=1-8]
	\arrow[from=2-1, to=2-2]
	\arrow[from=2-2, to=2-3]
	\arrow[from=2-2, to=3-2]
	\arrow[curve={height=-6pt}, from=2-3, to=3-3]
	\arrow[from=2-7, to=1-5]
	\arrow[from=2-7, to=3-7]
	\arrow[from=2-8, to=2-7]
	\arrow[from=2-8, to=3-8]
	\arrow[from=2-9, to=1-9]
	\arrow[from=3-2, to=2-1]
	\arrow[from=3-3, to=2-2]
	\arrow[curve={height=-6pt}, from=3-3, to=2-3]
	\arrow[from=3-3, to=3-2]
	\arrow[from=3-6, to=3-6, loop, in=60, out=120, distance=5mm]
	\arrow[from=3-6, to=4-9]
	\arrow[from=3-7, to=3-8]
	\arrow[from=3-8, to=3-9]
	\arrow[from=3-9, to=2-8]
	\arrow[from=3-9, to=2-9]
	\arrow[from=4-5, to=3-6]
	\arrow[curve={height=6pt}, from=4-5, to=4-9]
	\arrow[from=4-9, to=3-9]
\end{tikzcd}\]

We must find all infinte paths that are simultaneously allowable for both automata. Note that there are no sequences of three consecutive 1's for $\mathcal{H}(1/6)$, so that we may erase some states and transitions for $\overline{\mathcal{H}(3/16)}$, along with all the states that must transition to them. These are indicated with slashes and dashed arrows.
% https://q.uiver.app/#q=WzAsMTksWzYsMCwiXFxjYW5jZWx7MX0iXSxbNywwLCIxIl0sWzcsMSwiMSJdLFs2LDEsIjEiXSxbNSwwLCJcXGNhbmNlbHswfSJdLFs2LDIsIjAiXSxbNywyLCIxIl0sWzgsMywiMCJdLFswLDEsIjAiXSxbMSwxLCIwIl0sWzIsMSwiMCJdLFsyLDIsIjEiXSxbMSwyLCIxIl0sWzQsMywiMSJdLFs1LDIsIjAiXSxbOCwyLCIwIl0sWzgsMSwiMSJdLFs4LDAsIjAiXSxbNCwwLCIwIl0sWzAsMSwiIiwwLHsic3R5bGUiOnsiYm9keSI6eyJuYW1lIjoiZGFzaGVkIn19fV0sWzEsMiwiIiwwLHsic3R5bGUiOnsiYm9keSI6eyJuYW1lIjoiZGFzaGVkIn19fV0sWzUsNl0sWzIsM10sWzMsNV0sWzAsMCwiIiwxLHsicmFkaXVzIjoxLCJhbmdsZSI6MTgwLCJzdHlsZSI6eyJ0YWlsIjp7Im5hbWUiOiJhcnJvd2hlYWQifSwiYm9keSI6eyJuYW1lIjoiZGFzaGVkIn0sImhlYWQiOnsibmFtZSI6Im5vbmUifX19XSxbNCwwLCIiLDAseyJzdHlsZSI6eyJib2R5Ijp7Im5hbWUiOiJkYXNoZWQifX19XSxbMiw2XSxbOCw5XSxbOSwxMF0sWzEwLDExLCIiLDAseyJjdXJ2ZSI6LTF9XSxbMTEsOV0sWzExLDEyXSxbMTIsOF0sWzksMTJdLFsxMSwxMCwiIiwwLHsiY3VydmUiOi0xfV0sWzEzLDcsIiIsMCx7ImN1cnZlIjoxfV0sWzEzLDE0XSxbMTQsN10sWzE0LDE0LCIiLDAseyJyYWRpdXMiOjF9XSxbNiwxNV0sWzE1LDJdLFsxNSwxNl0sWzcsMTVdLFsxNywxXSxbMTYsMTddLFsxOCwxM10sWzE4LDE3LCIiLDAseyJjdXJ2ZSI6LTR9XSxbMywxOF0sWzE4LDQsIiIsMSx7InN0eWxlIjp7ImJvZHkiOnsibmFtZSI6ImRhc2hlZCJ9fX1dLFsxLDE2XV0=
\[\begin{tikzcd}
	&&&& 0 & {\cancel{0}} & {\cancel{1}} & 1 & 0 \\
	0 & 0 & 0 &&&& 1 & 1 & 1 \\
	& 1 & 1 &&& 0 & 0 & 1 & 0 \\
	&&&& 1 &&&& 0
	\arrow[dashed, from=1-5, to=1-6]
	\arrow[curve={height=-24pt}, from=1-5, to=1-9]
	\arrow[from=1-5, to=4-5]
	\arrow[dashed, from=1-6, to=1-7]
	\arrow[dashed, tail reversed, no head, from=1-7, to=1-7, loop, in=240, out=300, distance=5mm]
	\arrow[dashed, from=1-7, to=1-8]
	\arrow[dashed, from=1-8, to=2-8]
	\arrow[from=1-8, to=2-9]
	\arrow[from=1-9, to=1-8]
	\arrow[from=2-1, to=2-2]
	\arrow[from=2-2, to=2-3]
	\arrow[from=2-2, to=3-2]
	\arrow[curve={height=-6pt}, from=2-3, to=3-3]
	\arrow[from=2-7, to=1-5]
	\arrow[from=2-7, to=3-7]
	\arrow[from=2-8, to=2-7]
	\arrow[from=2-8, to=3-8]
	\arrow[from=2-9, to=1-9]
	\arrow[from=3-2, to=2-1]
	\arrow[from=3-3, to=2-2]
	\arrow[curve={height=-6pt}, from=3-3, to=2-3]
	\arrow[from=3-3, to=3-2]
	\arrow[from=3-6, to=3-6, loop, in=60, out=120, distance=5mm]
	\arrow[from=3-6, to=4-9]
	\arrow[from=3-7, to=3-8]
	\arrow[from=3-8, to=3-9]
	\arrow[from=3-9, to=2-8]
	\arrow[from=3-9, to=2-9]
	\arrow[from=4-5, to=3-6]
	\arrow[curve={height=6pt}, from=4-5, to=4-9]
	\arrow[from=4-9, to=3-9]
\end{tikzcd}\]

The upper right part of the right automaton now represents the cycle $\overline{011}$, which is not allowable for $\mathcal{H}(1/6)$. This allows us to erase these states.
% https://q.uiver.app/#q=WzAsMTcsWzcsMCwiXFxjYW5jZWx7MX0iXSxbNywxLCIxIl0sWzYsMSwiMSJdLFs2LDIsIjAiXSxbNywyLCIxIl0sWzgsMywiMCJdLFswLDEsIjAiXSxbMSwxLCIwIl0sWzIsMSwiMCJdLFsyLDIsIjEiXSxbMSwyLCIxIl0sWzQsMywiMSJdLFs1LDIsIjAiXSxbOCwyLCIwIl0sWzgsMSwiXFxjYW5jZWx7MX0iXSxbOCwwLCJcXGNhbmNlbHswfSJdLFs0LDAsIjAiXSxbMyw0XSxbMSwyXSxbMiwzXSxbMSw0XSxbNiw3XSxbNyw4XSxbOCw5LCIiLDAseyJjdXJ2ZSI6LTF9XSxbOSw3XSxbOSwxMF0sWzEwLDZdLFs3LDEwXSxbOSw4LCIiLDAseyJjdXJ2ZSI6LTF9XSxbMTEsNSwiIiwwLHsiY3VydmUiOjF9XSxbMTEsMTJdLFsxMiw1XSxbMTIsMTIsIiIsMCx7InJhZGl1cyI6MX1dLFs0LDEzXSxbMTMsMV0sWzEzLDE0LCIiLDAseyJzdHlsZSI6eyJib2R5Ijp7Im5hbWUiOiJkYXNoZWQifX19XSxbNSwxM10sWzE1LDAsIiIsMCx7InN0eWxlIjp7ImJvZHkiOnsibmFtZSI6ImRhc2hlZCJ9fX1dLFsxNCwxNSwiIiwwLHsic3R5bGUiOnsiYm9keSI6eyJuYW1lIjoiZGFzaGVkIn19fV0sWzE2LDExXSxbMTYsMTUsIiIsMCx7ImN1cnZlIjotNCwic3R5bGUiOnsiYm9keSI6eyJuYW1lIjoiZGFzaGVkIn19fV0sWzIsMTZdLFswLDE0LCIiLDEseyJzdHlsZSI6eyJib2R5Ijp7Im5hbWUiOiJkYXNoZWQifX19XV0=
\[\begin{tikzcd}
	&&&& 0 &&& {\cancel{1}} & {\cancel{0}} \\
	0 & 0 & 0 &&&& 1 & 1 & {\cancel{1}} \\
	& 1 & 1 &&& 0 & 0 & 1 & 0 \\
	&&&& 1 &&&& 0
	\arrow[curve={height=-24pt}, dashed, from=1-5, to=1-9]
	\arrow[from=1-5, to=4-5]
	\arrow[dashed, from=1-8, to=2-9]
	\arrow[dashed, from=1-9, to=1-8]
	\arrow[from=2-1, to=2-2]
	\arrow[from=2-2, to=2-3]
	\arrow[from=2-2, to=3-2]
	\arrow[curve={height=-6pt}, from=2-3, to=3-3]
	\arrow[from=2-7, to=1-5]
	\arrow[from=2-7, to=3-7]
	\arrow[from=2-8, to=2-7]
	\arrow[from=2-8, to=3-8]
	\arrow[dashed, from=2-9, to=1-9]
	\arrow[from=3-2, to=2-1]
	\arrow[from=3-3, to=2-2]
	\arrow[curve={height=-6pt}, from=3-3, to=2-3]
	\arrow[from=3-3, to=3-2]
	\arrow[from=3-6, to=3-6, loop, in=60, out=120, distance=5mm]
	\arrow[from=3-6, to=4-9]
	\arrow[from=3-7, to=3-8]
	\arrow[from=3-8, to=3-9]
	\arrow[from=3-9, to=2-8]
	\arrow[dashed, from=3-9, to=2-9]
	\arrow[from=4-5, to=3-6]
	\arrow[curve={height=6pt}, from=4-5, to=4-9]
	\arrow[from=4-9, to=3-9]
\end{tikzcd}\]

Every instance of the bit sequence $11$ for $\mathcal{H}(1/6)$ has to be continued by $00$. This allows us to erase more states on the right automaton:
% https://q.uiver.app/#q=WzAsMTQsWzcsMSwiXFxjYW5jZWx7MX0iXSxbNiwxLCIxIl0sWzYsMiwiXFxjYW5jZWx7MH0iXSxbNywyLCJcXGNhbmNlbHsxfSJdLFs4LDMsIlxcY2FuY2VsezB9Il0sWzAsMSwiMCJdLFsxLDEsIjAiXSxbMiwxLCIwIl0sWzIsMiwiMSJdLFsxLDIsIjEiXSxbNCwzLCIxIl0sWzUsMiwiMCJdLFs4LDIsIlxcY2FuY2VsezB9Il0sWzQsMCwiMCJdLFsyLDMsIiIsMCx7InN0eWxlIjp7ImJvZHkiOnsibmFtZSI6ImRhc2hlZCJ9fX1dLFswLDEsIiIsMCx7InN0eWxlIjp7ImJvZHkiOnsibmFtZSI6ImRhc2hlZCJ9fX1dLFsxLDIsIiIsMCx7InN0eWxlIjp7ImJvZHkiOnsibmFtZSI6ImRhc2hlZCJ9fX1dLFswLDMsIiIsMCx7InN0eWxlIjp7ImJvZHkiOnsibmFtZSI6ImRhc2hlZCJ9fX1dLFs1LDZdLFs2LDddLFs3LDgsIiIsMCx7ImN1cnZlIjotMX1dLFs4LDZdLFs4LDldLFs5LDVdLFs2LDldLFs4LDcsIiIsMCx7ImN1cnZlIjotMX1dLFsxMCw0LCIiLDAseyJjdXJ2ZSI6MSwic3R5bGUiOnsiYm9keSI6eyJuYW1lIjoiZGFzaGVkIn19fV0sWzEwLDExXSxbMTEsNCwiIiwwLHsic3R5bGUiOnsiYm9keSI6eyJuYW1lIjoiZGFzaGVkIn19fV0sWzExLDExLCIiLDAseyJyYWRpdXMiOjF9XSxbMywxMiwiIiwwLHsic3R5bGUiOnsiYm9keSI6eyJuYW1lIjoiZGFzaGVkIn19fV0sWzEyLDAsIiIsMCx7InN0eWxlIjp7ImJvZHkiOnsibmFtZSI6ImRhc2hlZCJ9fX1dLFs0LDEyLCIiLDAseyJzdHlsZSI6eyJib2R5Ijp7Im5hbWUiOiJkYXNoZWQifX19XSxbMTMsMTBdLFsxLDEzXV0=
\[\begin{tikzcd}
	&&&& 0 &&&& \\
	0 & 0 & 0 &&&& 1 & {\cancel{1}} \\
	& 1 & 1 &&& 0 & {\cancel{0}} & {\cancel{1}} & {\cancel{0}} \\
	&&&& 1 &&&& {\cancel{0}}
	\arrow[from=1-5, to=4-5]
	\arrow[from=2-1, to=2-2]
	\arrow[from=2-2, to=2-3]
	\arrow[from=2-2, to=3-2]
	\arrow[curve={height=-6pt}, from=2-3, to=3-3]
	\arrow[from=2-7, to=1-5]
	\arrow[dashed, from=2-7, to=3-7]
	\arrow[dashed, from=2-8, to=2-7]
	\arrow[dashed, from=2-8, to=3-8]
	\arrow[from=3-2, to=2-1]
	\arrow[from=3-3, to=2-2]
	\arrow[curve={height=-6pt}, from=3-3, to=2-3]
	\arrow[from=3-3, to=3-2]
	\arrow[from=3-6, to=3-6, loop, in=60, out=120, distance=5mm]
	\arrow[dashed, from=3-6, to=4-9]
	\arrow[dashed, from=3-7, to=3-8]
	\arrow[dashed, from=3-8, to=3-9]
	\arrow[dashed, from=3-9, to=2-8]
	\arrow[from=4-5, to=3-6]
	\arrow[curve={height=6pt}, dashed, from=4-5, to=4-9]
	\arrow[dashed, from=4-9, to=3-9]
\end{tikzcd}\]
but then the only allowable infinite bit sequences for the right automaton are shifts of $.101\overline{0} = 5/8$, which are not allowed for $\mathcal{H}(1/6)$. This concludes that there are no direct ray connections between the Hubbard trees in the formal mating.

Any biaccessible point eventually iterates into the Hubbard tree (Lemma \ref{lemma: biaccessible}), so there are also no indirect ray connections between the two Hubbard trees; any one could be iterated to a direct ray connection.
\end{exm}

\begin{exm} \label{example: 1/4 with 1/6}
We now consider an example analogous to that of Appendix C in \cite{Milnor2004}. Computing the mating automaton for $f_{1/4} \sqcup f_{1/6}$ according to the product construction, we obtain the labeled directed graph below:
% https://q.uiver.app/#q=WzAsNixbMCwxLCIxIl0sWzEsMSwiMSJdLFsyLDEsIjAiXSxbMiwyLCIwIl0sWzEsMiwiMSJdLFswLDAsIjAiXSxbMCwxXSxbMSwyXSxbMyw0XSxbNCwwXSxbMiwzXSxbNSwwXSxbNSwxXV0=
\[\begin{tikzcd}
	0 && \\
	1 & 1 & 0 \\
	& 1 & 0
	\arrow[from=1-1, to=2-1]
	\arrow[from=1-1, to=2-2]
	\arrow[from=2-1, to=2-2]
	\arrow[from=2-2, to=2-3]
	\arrow[from=2-3, to=3-3]
	\arrow[from=3-2, to=2-1]
	\arrow[from=3-3, to=3-2]
\end{tikzcd}\]
It is clear from the automaton that there are only finitely many direct ray connections between the Hubbard trees. They are all all preperiodic, iterating into the periodic cycle
$$
\left\{ \dfrac{7}{31}, \dfrac{14}{31}, \dfrac{28}{31}, \dfrac{25}{31}, \dfrac{19}{31} \right\}
$$
with respect to the $c_{1/4}$-plane, or
$$
\left\{ \dfrac{24}{31}, \dfrac{17}{31}, \dfrac{3}{31}, \dfrac{6}{31}, \dfrac{12}{31} \right\} 
$$
with respect to the $c_{1/6}$-plane. The angle $7/31$ is characteristic for $f_{1/4}$, where the ray $\R_{c_{1/4}}(7/31)$ lands together with $\R_{c_{1/4}}(8/31)$ at a hyperbolic component of period $5$ on the vein of $c_{1/4}$. 

Now fix the parameter $c_{1/6}$, and move the parameter $c_{1/4}$ ``down'' along its vein. We do so until we pass the root of the hyperbolic component $W$ which $P_M(7/31, 8/31)$ lands at the root of. After passing it, because of the realizability of orbit portraits within wakes \cite{Milnor2000, Schleicher2000}, the dynamical rays $\R_c(7/31)$ and $\R_c(8/31)$ no longer land together, and they do not land at the Hubbard tree at all; if they did, they wouldn't land at an extremity since the angles are periodic, and so they'd land at biaccessible points. These would form orbit portraits containing these two angles, which would be preserved by moving the parameter back up along the vein; but $7/31$ and $8/31$ are already part of a different portrait in the wake of $W$.

Moving $c \in V_{1/4}$ past the root of $W$ ``destroys'' the ray connection: more precisely, for any parameter $c$ on the vein to $c_{1/4}$ but preceding this root, its mating with $c_{1/6}$ has no direct ray connections between the Hubbard trees, even for parameters that are not PCF. This is elaborated upon in the course of the proof of Theorem \ref{theorem: first ray connection}. Essentially, any direct ray connection between the parameters preceding them in the vein would be present for the mating between the tips $c_{1/4}$ and $c_{1/6}$.
\end{exm}

Destroying the direct ray connections between the Hubbard trees is possible because we can identify what is the ``first'' ray connection between $c_{1/4}$ and $c_{1/6}$ along the vein of $c_{1/4}$, and that it was periodic, corresponding to the root of some hyperbolic component. It may be the case that this first ray connection does not appear at parabolic parameters, as in the examples below.

\begin{exm} \label{example: tips of veins}
Consider the mating $f_{1/4} \sqcup f_{5/16}$, which, after a suitable trimming, has the product automaton below:
% https://q.uiver.app/#q=WzAsNCxbMCwwLCIxIl0sWzEsMCwiMCJdLFsxLDEsIjEiXSxbMCwxLCIwIl0sWzAsMV0sWzMsMl0sWzEsMSwiIiwwLHsicmFkaXVzIjoxLCJhbmdsZSI6OTB9XSxbMiwyLCIiLDAseyJyYWRpdXMiOjEsImFuZ2xlIjo5MH1dXQ==
\[\begin{tikzcd}
	1 & 0 \\
	0 & 1
	\arrow[from=1-1, to=1-2]
	\arrow[from=1-2, to=1-2, loop, in=330, out=30, distance=5mm]
	\arrow[from=2-1, to=2-2]
	\arrow[from=2-2, to=2-2, loop, in=330, out=30, distance=5mm]
\end{tikzcd}\]
The only direct ray connections between the Hubbard trees are those given by $\phi = 0$ and $\phi = 1/2$, that is, connecting the $\beta$ and $-\beta$ points. In particular, the Hubbard trees form a topological graph with a single cycle in the formal mating. The two direct ray connections are picked up exactly at the tips of the respective veins, which are Misiurewicz points, and any mating of $c_{1/4}$ with $c' \prec c_{5/16}$ would have no direct ray connections between the Hubbard trees, since $\beta$ and $-\beta$ are no longer points of the Hubbard tree $H_{c'}$.

The ray connections may not come from just the forward orbit of the external argument of a tip of a vein. In the mating $f_{1/4} \sqcup f_{11/32}$, the mating automaton from the product construction, after a suitable trimming of states, takes the form:
% https://q.uiver.app/#q=WzAsNyxbMCwxLCIwIl0sWzEsMSwiMSJdLFsyLDEsIjAiXSxbMywwLCIxIl0sWzQsMCwiMCJdLFszLDEsIjAiXSxbNCwxLCIxIl0sWzIsM10sWzMsNF0sWzIsNV0sWzQsNCwiIiwyLHsicmFkaXVzIjoxLCJhbmdsZSI6OTB9XSxbNiw2LCIiLDIseyJyYWRpdXMiOjEsImFuZ2xlIjo5MH1dLFs1LDZdLFswLDFdLFsxLDJdXQ==
\[\begin{tikzcd}
	&&& 1 & 0 \\
	0 & 1 & 0 & 0 & 1
	\arrow[from=1-4, to=1-5]
	\arrow[from=1-5, to=1-5, loop, in=330, out=30, distance=5mm]
	\arrow[from=2-1, to=2-2]
	\arrow[from=2-2, to=2-3]
	\arrow[from=2-3, to=1-4]
	\arrow[from=2-3, to=2-4]
	\arrow[from=2-4, to=2-5]
	\arrow[from=2-5, to=2-5, loop, in=330, out=30, distance=5mm]
\end{tikzcd}\]
which also has a finite collection of angles, explicitly
$$
\left\{ \dfrac{5}{16}, \dfrac{5}{8}, \dfrac{1}{4}, \dfrac{1}{2}, \dfrac{0}{1} \right\},
$$
and some of them are contained in the backward orbit of $1/4$.
\end{exm}

\begin{exm} \label{example: 1/4 with 7/8}
A mating automaton for $f_{1/4} \sqcup f_{7/8}$, after a suitable trimming that preserves the regular langauge, is given below:
% https://q.uiver.app/#q=WzAsMTMsWzEsMSwiMCJdLFsyLDEsIjAiXSxbMiwwLCIxIl0sWzEsMCwiMCJdLFswLDEsIjAiXSxbMCwyLCIwIl0sWzEsMiwiMSJdLFszLDEsIjAiXSxbMywwLCIxIl0sWzQsMCwiMCJdLFs1LDAsIjEiXSxbNSwxLCIwIl0sWzQsMSwiMCJdLFswLDFdLFsxLDJdLFsyLDNdLFszLDMsIiIsMCx7InJhZGl1cyI6MSwiYW5nbGUiOi00NX1dLFszLDBdLFswLDRdLFs0LDVdLFs1LDZdLFs2LDBdLFsxLDddLFs3LDhdLFs4LDldLFs5LDEwXSxbOCw4LCIiLDAseyJyYWRpdXMiOjEsImFuZ2xlIjotNDV9XSxbMTAsMTFdLFsxMSwxMl0sWzEyLDldXQ==
\[\begin{tikzcd}
	& 0 & 1 & 1 & 0 & 1 \\
	0 & 0 & 0 & 0 & 0 & 0 \\
	0 & 1
	\arrow[from=1-2, to=1-2, loop, in=105, out=165, distance=5mm]
	\arrow[from=1-2, to=2-2]
	\arrow[from=1-3, to=1-2]
	\arrow[from=1-4, to=1-4, loop, in=105, out=165, distance=5mm]
	\arrow[from=1-4, to=1-5]
	\arrow[from=1-5, to=1-6]
	\arrow[from=1-6, to=2-6]
	\arrow[from=2-1, to=3-1]
	\arrow[from=2-2, to=2-1]
	\arrow[from=2-2, to=2-3]
	\arrow[from=2-3, to=1-3]
	\arrow[from=2-3, to=2-4]
	\arrow[from=2-4, to=1-4]
	\arrow[from=2-5, to=1-5]
	\arrow[from=2-6, to=2-5]
	\arrow[from=3-1, to=3-2]
	\arrow[from=3-2, to=2-2]
\end{tikzcd}\]
it has entropy approximately $0.36354\ldots$, and so there are uncountably many direct ray connections between the Hubbard trees $H_{1/4}$ and $H_{7/8}$. Fixing the parameter $c_{7/8}$ and varying the parameter $c_{1/4}$ along its vein, it is not obvious at which parameter $c \in V_{c_{1/4}}$ the mating $f_{c} \sqcup f_{7/8}$ first acquires ray connections between the Hubbard trees. Theorem \ref{theorem: first ray connection} guarantees that this happens at some root of a hyperbolic component, and the proof will suggest an algorithm to find the rays landing at this parameter explicitly.
\end{exm}

\section{Characteristic Points and Characteristic States} \label{section: characteristic states automaton}

With the intent of proving Theorem \ref{theorem: first ray connection}, the pathological behavior to be excluded is that the first ray connection in a mating $f_{c_1} \sqcup f_{c_2}$ along a vein could be picked up at an irrational angle $\phi$. This would imply that there is a parameter $c \prec c_1$ and on the impression of $\R_M(\phi)$ which is accumulated on from above by parabolic parameters corresponding to periodic direct ray connections. To rule this phenomenon out, we need to better understand the relationship between characteristic points, direct ray connections between generalized trees, and mating automata for PCF parameters.

\subsection{Characteristic direct ray connections}

\begin{definition}
A direct ray connection $\R_{f_{c_1} \sqcup f_{c_2}}(\phi)$ is \textit{characteristic} for $f_{c_1}$ if the landing point of $\R_{c_1}(\phi)$ is characteristic. By a \textit{characteristic direct ray connection} we mean one that is characteristic for the first argument in the mating operation, when there is no risk of confusion.
\end{definition}

In the above definition, the landing point $z$ of $\R_{c_1}(\phi)$ is therefore in $\{\beta_{c_1}\} \cup [\alpha_{c_1}, c_1]$, since there are no characteristic points on $(\beta_{c_1}, \alpha_{c_1})$. Even though $\beta_{c_1}$ is a characteristic point for $f_{c_1}$, it is on the generalized Hubbard tree $T_{c_1}$ if and only if $c_1$ is a dyadic Misiurewicz point; this is because $\beta_{c_1}$ is has only one external argument, so it must be an extremity of $T_{c_1}$, and therefore in the postcritical set. Hence $c_1$ has only one external argument which is a dyadic rational, which is equivalent to the parameter $c_1$ being dyadic Misiurewicz. This fact, coupled with details about the parameter-dynamical correspondence of characteristic points, will justify our rationale for working with the arc $[\alpha_{c_1}, c_1]$ instead of $[\beta_{c_1}, c_1]$.

Note that $\R_{f_{c_1} \sqcup f_{c_2}}(\phi)$ being a characteristic direct ray connection does not mean that $\phi$ is an endpoint of a characteristic leaf for $L_{c_1}$. Indeed, if $z \in T_{c_1} \cap J_{c_1}$ is characteristic and has $3$ or more external arguments, exacly one leaf formed by its external arguments is characteristic if $z$ periodic, and exactly two if $z$ is strictly preperiodic.

\begin{lemma} \label{lemma: characteristic ray connections closed}
Let $c_1$ and $c_2$ be parameters on PCF veins, and $T_{c_1}$, $T_{c_2}$ be their generalized Hubbard trees. Then the set of characteristic landing points on $T_{c_1}$ of direct ray connections between $T_{c_1}$ and $T_{c_2}$ is closed.
\end{lemma}

\begin{proof}
As stated in Section \ref{subsection: quadratic invariant laminations}, the set of characteristic points of $f_{c_1}$ is closed. Suppose $\phi_n \to \phi$ is a sequence of angles characteristic such that $\R_{c_1}(\phi)$ lands on $z_n \in T_{c_1}$ and $\R_{c_2}(-\phi)$ lands on $w_n \in T_{c_2}$. By continuity and closedness of the generalized Hubbard trees, $\R_{c_1}(\phi)$ lands on $z = \lim z_n \in T_{c_1}$ and $\R_{c_2}(-\phi)$ lands on $w = \lim w_n \in T_{c_2}$. Therefore the set of angles of direct ray connections between the trees is closed, and the projection under the Carathéodory loop to $\partial K_{c_1}$ is closed. Therefore the set of characteristic landing points of direct ray connections is the intersection of two closed sets.
\end{proof}

Since closed subsets of $\RR$ always have minima and $\{\beta_{c_1}\}$ is discrete as a characteristic point on $[\beta_{c_1}, c_1]$, the following is immediate:

\begin{cor} \label{corollary: exists minimal characteristic}
Let $c_1$ and $c_2$ be parameters on PCF veins, and suppose that there exists a direct ray connection between $T_{c_1}$ and $T_{c_2}$ that is characteristic for $f_{c_1}$. Then there exists a characteristic landing point $z_0 \in T_{c_1} \cap [\alpha_{c_1}, c_1]$ of a direct ray connection between the generalized trees that is minimal with respect to the linear order on $[\alpha_{c_1}, c_1]$. More precisely, it precedes every other characteristic landing point of a direct ray connection between the trees which is on $[\alpha_{c_1}, c_1]$.
\end{cor}

Corollary \ref{corollary: exists minimal characteristic} goes a long way in proving Theorem \ref{theorem: first ray connection}; essentially the only missing pieces are showing that $z_0$ is preperiodic if it exists, and understanding what that means in the parameter plane. Lemma \ref{lemma: z0 strictly preperiodic} will show that if $z_0$ is strictly preperiodic, then the first two statements of case 2) in Theorem \ref{theorem: first ray connection} hold, and we will show in Section \ref{section: algorithm for first ray connection} that if $z_0 \in [\alpha, c_1]$ and it is periodic, then case 3) holds.

\begin{lemma} \label{lemma: exists characteristic}
Let $c_1$ and $c_2$ be parameters on PCF veins, and suppose that there exists a ray connection between the generalized Hubbard tree $T_{c_1}$ and $T_{c_2}$. Then there exists a direct ray connection between $T_{c_1}$ and $T_{c_2}$ that is characteristic for $f_{c_1}$, or $c_1$ is the tip of a PCF vein.
\end{lemma}

\begin{proof}
Analogously to the case of Hubbard trees, every biaccessible point in $J_c$ eventually iterates into $[\alpha, c_1] \subset T_{c_1}$. Therefore we may assume that there exists a direct ray connection between the generalized trees, since any indirect ray connection could be iterated to contain a direct one.

Let $\R_{f_{c_1} \sqcup f_{c_2}}(\phi)$ be a direct ray connection between the generalized trees, $z$ be the landing point of $\R_{c_1}(\phi)$, and $z_n = f^n(z)$. If the orbit of $z$ never enters the arc $[\alpha_{c_1}, c_1]$, it is because $z$ is not biaccessible, and it must be an extremity of the generalized tree. But the only way an extremity of the tree is not biaccessible is if $c_1$ has only one external argument, and so $c_1$ is the tip of a vein, and by hypothesis it is a PCF parameter.

If the orbit of $z$ passes through $[\alpha_{c_1}, c_1]$ but does not return to this arc infinitely often, we may just take the maximal such iterate on $[\alpha, c_1]$, which is by assumption characteristic and the landing point of a forward image of $\R_{f_{c_1} \sqcup f_{c_2}}(\phi)$.

Assume then that the orbit of $\phi$ returns infinitely often to $[\alpha, c_1]$, and let $z_* = \limsup_{n \to \infty} z_n$ in the arc $[\alpha_{c_1}, c_1]$, linearly ordered from $\alpha_{c_1}$ to $c_1$. We claim that $z_*$ is characteristic and the landing point of a direct ray connection between the generalized Hubbard trees. Taking a subsequence $\phi_{n_k}$ such that $\phi_{n_k} \to \phi_*$, $\R_{c_1}(\phi_*)$ lands at $z_*$ by continuity of the Carathéodory loop, and the same arguments from Lemma \ref{lemma: characteristic ray connections closed} guarantee that $z_*$ is the landing point of a direct ray connection between the generalized Hubbard trees. Therefore it is enough to show that $z_*$ is characteristic. If $z_*$ were not characteristic, then there would be some iterate $f^j(z_*)$ such that $f^j(z_*) \in (z_*, c_1]$, due to Proposition \ref{proposition: characteristic points}. As this is an open condition, by continuity there would be some $z_{n_k}$ for $k$ sufficiently large such that $f^j(z_{n_{k}}) = z_{n_k + j}$ is in $(z_*, c_1]$, a contradiction with the definition of $z_*$.
\end{proof}

The following result strengthens the previous lemma if we assume that the parameters are PCF:

\begin{prop} \label{proposition: exists periodic characteristic}
Let $c_1$ and $c_2$ be PCF parameters, and suppose that there exists a ray connection between the Hubbard trees $H_{c_1}$ and $H_{c_2}$. Then there exists a periodic characteristic direct ray connection between $H_{c_1}$ and $H_{c_2}$, or $c_1$ is the tip of a PCF vein, and all direct ray connections between the Hubbard trees are in the grand orbit of $\R_{f_1 \sqcup f_2}(\theta_1)$.
\end{prop}

\begin{proof}
As in the proof of Lemma \ref{lemma: exists characteristic}, we may assume there exists a direct ray connection, so that $\mathcal{L}(\mathcal{M}) \neq \emptyset$ for a union automaton $\mathcal{M}$ for $f_{c_1} \sqcup f_{c_2}$. By Lemma \ref{lemma: nonempty regular language implies existence of periodic path}, there exists a periodic path in $\mathcal{M}$, and so a periodic direct ray connection between the Hubbard trees. If the angle $\phi$ of this direct ray connection is part of a non-trivial orbit portrait, then one of the landing points of $\R_{c_1}(2^i\phi)$ is characteristic and periodic.

Suppose now that there are no characteristic periodic direct ray connections between the Hubbard trees, so that the angles associated to all periodic paths in $\mathcal{M}$ belong to trivial orbit portraits for $f_{c_1}$. The periodic points of these portraits must be in $H_{c_1}$ and are not biaccessible, so must be extremities of the Hubbard tree having only one external argument. In particular $c_1$ has only one external argument $\theta_1$, and so it is the tip of a PCF vein. Since all extremities of $H_{c_1}$ are postcritical points, all periodic paths in the automaton have associated angles in the periodic cycle $\{\theta_{m+1}, \ldots \theta_{m +n}\}$ that $\theta_1$ iterates into. There are at most three such periodic cycles of states, being the path on angle states $X_1 \to \cdots \to X_n \to X_1$, and the two paths on arc states that accumulate on the marked angles from both sides, $Y_1 \to \cdots \to Y_n \to Y_1$ and $Z_1 \to \cdots \to Z_n \to Z_1$, if the latter ones are recognized by $\mathcal{M}$.

As the same arguments below can be applied in a simpler form if either only one or none of the paths on arc states are recognized, we will assume that both are. Let $\mathcal{X} = \{X_i\}_i$, $\mathcal{Y} = \{Y_j\}_j$, and $\mathcal{Z} = \{Z_k\}_k$, and $\mathcal{U} = \mathcal{X} \cup \mathcal{Y} \cup \mathcal{Z}$. We claim that $\mathcal{U}$ is a cross section (Definition \ref{definition: cross section}). This is because if there were an infinite path that avoided $\mathcal{U}$, then there is an infinite periodic path that does so, a contradiction with the fact that every periodic path has associated periodic orbit of angles $\{\theta_{m+1}, \ldots \theta_{m +n}\}$, and these are only given by the $\mathcal{X}$, $\mathcal{Y}$ and $\mathcal{Z}$. We can be more precise: there are no transitions from $\mathcal{X}$ to states not in $\mathcal{X}$ and the only transitions are $X_i \to X_{i+1}$. If there were a path from $\mathcal{Y}$ to $\mathcal{Y}$ that is not just $Y_j \to Y_{j+1}$, this could be concatenated with the periodic cycle of $\mathcal{Y}$ to give a different periodic path in $\mathcal{M}$. The same applies to $\mathcal{Z}$.

There can be no paths from $\mathcal{Y}$ to $\mathcal{Z}$ and paths from $\mathcal{Z}$ to $\mathcal{Y}$ simultaneously, otherwise concatenating these paths (with suitable transitions of the periods of $\mathcal{Y}$ and $\mathcal{Z}$) we could produce yet another different periodic path in $\mathcal{M}$. Hence either $\mathcal{X} \cup \mathcal{Y}$ or $\mathcal{X} \cup \mathcal{Z}$ has no transitions out of these sets; assume without loss of generality that it is $\mathcal{X} \cup \mathcal{Y}$. The only missing possibility is that there is a path from $\mathcal{Y}$ to $\mathcal{X}$. There being one or not, every infinite path in $\mathcal{M}$ eventually falls into the periodic cycles of either $\mathcal{X}$, $\mathcal{Y}$ or $\mathcal{Z}$. This concludes that the angle of any direct ray connection between the Hubbard trees must be in the grand orbit of $\theta_1$.
\end{proof}

\begin{remark}
An example where there are direct ray connections between the Hubbard trees (even periodic), but none that are characteristic for $f_{c_1}$, is the mating $f_{1/6} \sqcup f_{1/3}$, where $c_{1/3}$ is the parameter of the Basilica polynomial. The set of angles forming ray connections is $\{1/3, 2/3\}$, but this corresponds to a trivial orbit portrait for $f_{1/6}$, whose points are extremities of $H_{1/6}$ which are not biaccessible. Additionally, the possibilities in Lemma \ref{lemma: exists characteristic} are not exclusive. If $c_1$ is a dyadic Misiurewicz parameter, the $\beta$ fixed point is the characteristic periodic point of the zero portrait, which by assumption is non-trivial.
\end{remark}

\begin{lemma} \label{lemma: z0 strictly preperiodic}
Let $c_1$ and $c_2$ be paramaters on PCF veins. If $z_0$ as defined in Corollary \ref{corollary: exists minimal characteristic} exists and is strictly preperiodic, then $z_0 = c_1$ and $c_1$ is the tip of a PCF vein. Moreover, if $c_2$ is PCF, then the only direct ray connections between the Hubbard trees are in the grand orbit of $\R_{f_{c_1} \sqcup f_{c_2}}(\theta_1)$, where $\theta_1$ is the external argument of $c_1$.
\end{lemma}

\begin{proof}
Suppose first that the periodic cycle that $z_0$ iterates into corresponds to a non-trivial orbit portrait. If $z_k$ is the characteristic point of said periodic cycle, it must precede $z_0$ on $T_{c_1} \cap [\beta, c_1]$, and still be the landing point of a direct ray connection between the generalized trees. The only possibility that does not contradict minimality of $z_0$ is $z_k = \beta$. Hence $z_0$ is an extremity of $T_{c_1}$, $z_0 = c_1$ and $c_1$ is a dyadic Misiurewicz parameter having external argument $\theta_1$.

As $\R_{c_2}(-\theta_1)$ must land on $T_{c_2}$ and $-\theta_1$ is dyadic, the $\beta$ fixed point $\beta_2$ for $f_{c_2}$ must be in the generalized tree, which implies that $\beta_2$ is in the postcritical orbit and $c_2$ is also a dyadic Misiurewicz parameter. In particular, both parameters are PCF. Apart from $\R_{f_{c_1 \sqcup f_{c_2}}}(0)$, there can be no other periodic direct ray connections in non-trivial orbit portraits for $f_{c_1}$, otherwise some forward iterate would land at a characteristic point for $f_{c_1}$ and contradict minimality of $z_0$. The same arguments in the proof of Proposition \ref{proposition: exists periodic characteristic} then that the periodic cycles of states in the union automaton $\mathcal{M}$ corresponding to the $0$ angle form a cross section, and that the rest of the statement follows.

If $z_0$ iterates into a trivial orbit portrait, then again it iterates to an extremity of the generalized tree. It wouldn't be an extremity of $T_{c_1}$ only if it were the critical point or precritical, but those are never characteristic since they iterate to $c_1$. So $z_0 = c_1$, $c_1$ is not biaccessible, and $c_1$ is the tip of a PCF vein. With the hypothesis that $c_2$ is PCF, again the arguments in the proof of \ref{proposition: exists periodic characteristic} show that all direct ray connections are in the grand orbit of $\R_{f_{c_1} \sqcup f_{c_2}}(\theta_1)$.
\end{proof}

\subsection{Characteristic states and words} \label{subsection: characteristic states and words}

Let $c$ be a PCF parameter, and $\mathcal{H}(c)$ a Hubbard tree automaton for $f_c$. We assume that the set of marked angles has ben sufficiently enlarged as in Section \ref{subsection: an involution on arc states}. 

There is a natural counterclockwise cyclic order on $\mathcal{S}$, corresponding to the counterclockwise arrangement of open arcs and marked angles on $S^1$. Given $A_0 \in \mathcal{S}$ as a fixed initial state, this cyclic order induces a linear lexicographic order on the set of words (infinite, or finite of the same length) $u$ such that $u_0 = A_0$. If $u$ and $v$ are such words such that $u_0 = v_0 = A_0$, we denote this relation by $u \preceq v$. This relation extends also to words of finite but different lengths, as long as one is not an extension of the other. This order reflects exactly the partitions of $A_0$ with respect to the subsets $J_{A_0A_{j_1}\ldots A_{j_i}}$ depending on the path length, arranged linearly along $A$.

As a consequence of the linear arrangement of the collections of nested arcs, and the fact that the associated angles are in the nested intersections of their closures, we obtain the following:

\begin{lemma} \label{lemma: lexicographic order on words transfers to angles}
If $u \preceq v$ are infinite words in $\mathcal{S}^{\NN}$ starting at the same state, and $\phi, \psi$ are the associated angles respectively, then $\phi \leq \psi$.
\end{lemma}

\begin{remark}
Note that the only case where $u \prec v$ but $\phi = \psi$ is if both paths accumulate on the same preimage of a marked angle, $u$ from below and $v$ from above.  
\end{remark}

As described in Section \ref{subsection: an involution on arc states}, there is a well defined involution $\iota : \mathcal{S}_{\mathrm{arc}} \cap \mathcal{S}_{\mathrm{good}} \to \mathcal{S}_{\mathrm{arc}} \cap \mathcal{S}_{\mathrm{good}}$, which extends naturally to the set of words in $\mathcal{S}_{\mathrm{arc}} \cap \mathcal{S}_{\mathrm{good}}$ and arc paths in $\mathcal{H}(c)$, as a consequence of Corollary \ref{corollary: involution preserves arc paths}. The involution reverses the lexicographic order:

\begin{lemma} \label{lemma: involution comparison}
Suppose that $u$ and $v$ are arc paths in $\mathcal{H}$ starting at the same state $A_0$ such that $u \preceq v$. Then $\iota(u) \succeq \iota(v)$. 
\end{lemma}

\begin{proof}
Consider the associated regular arcs on $H_c$ to the itinerary arcs $J_{A_0 A_{j_1} \ldots A_{j_i}}$. Since the external arguments of their endpoints cannot cross, and they are nested according to Corollary \ref{corollary: landing points of J arcs}, the linear arrangement of the $J_{A_0 A_{j_1} \ldots A_{j_i}}$ exactly corresponds to how the associated regular arcs are arranged on the regular arc associated to $A_0$, in a preferred order. Moreover, Proposition \ref{proposition: J arcs respect involution} guarantees that that order is reflected by the involution, since the reflected arcs have the same associated regular arcs on $H_c$. 
\end{proof}

The set of marked angles induces a partition of $H_{c}$ by their landing points. If $c$ is a Misiurewicz parameter, the arc $[\alpha, c]$ is therefore divided into a union of segments $\bigcup_{i = 0}^N [x_i, x_{i+1}]$, with $\alpha = x_0 < x_1 < \cdots < x_{N+1} = c$, such that the common endpoints are the landing points of marked angles, and there are no marked angles whose rays land on $(x_i, x_{i+1})$. Each segment $[x_i, x_{i+1}]$ has two arc states $A_i$ and $B_i$ associated to it, corresponding to the two sides of the edge. More precisely, each segment determines four angles in counterclockwise order: $\{\phi_i^1, \phi^2_i, \phi^3_i, \phi^4_i\}$, such that $\phi^1_i$ and $\phi^4_i$ are the external arguments of $x_i$ adjacent to $[x_i, x_{i+1}]$, $\phi^2_i$ and $\phi^3_i$ are the external arguments of $x_{i+1}$ adjacent to $[x_i, x_{i+1}]$, and $A_i = (\phi^1_i, \phi^2_i)$, $B_i = (\phi_i^3, \phi_i^4)$. (The angles $\phi_i^2$ and $\phi_i^3$ may coincide if $x_{i+1} = c$ is not biaccessible.) Some of these arc states may not be good, occuring if and only if $(x_i, x_{i+1})$ is contained in the Fatou set (see the proof of Lemma \ref{lemma: A' is good arc state}). We may still define an involution on these states, and we have $\iota(A_i) = B_i$. If instead $f_c$ is critically periodic, we consider the partition only of $[\alpha, c^*]$, where $c^* = x_{N+1}$ is the root if the Fatou component containing the critical value. The same description as above holds.

\begin{definition} \label{definition: characteristic states}
The \textit{characteristic states} of the Hubbard tree automaton $\mathcal{H}(c)$ is the collection of states $\{A_0, B_0, \ldots, A_{N}, B_{N}\}$ as constructed above.
\end{definition}

Even though the definition of the characteristic states depends on the choice of set of marked angles (and requires a preliminary enlargement), the results that follow do not depend on the choice of enlargement.

We will require some condition on paths in $\mathcal{H}(c)$ to guarantee that they correspond to characteristic points for $f_c$. The following is a way to do so for paths passing through only arc states.

\begin{definition} \label{definition: arc characteristic word}
Let $w$ be an arc path in $\mathcal{H}(c)$ and such that $w_0 \in \{A_i, B_i\}_i$. Let $\mathcal{D}_i$ be the set of states $D$ such that $A_i < D < B_i$. The arc path $w$ is \textit{characteristic} if:

\begin{itemize}
    \item For all $k > 0$, $w_k \notin \mathcal{D}_i$;

    \item If $w_0 = A_i$, and $k > 0$ is such that $w_k = A_i$ or $B_i$, then $w \succeq \sigma^k(w)$ or $\iota(w) \preceq \sigma^k(w)$, respectively;

    \item If $w_0 = B_i$, and $k > 0$ is such that $w_k = A_i$ or $B_i$, then $\iota(w) \succeq \sigma^k(w)$ or $w \preceq \sigma^k(w)$, respectively.
\end{itemize}    
\end{definition}

As a consequence of Corollary \ref{corollary: landing points of J arcs}, Proposition \ref{proposition: J arcs respect involution}, and Lemmas \ref{lemma: lexicographic order on words transfers to angles} and \ref{lemma: involution comparison}, the last two conditions on $w = Z_0Z_1\ldots$ being a characteristic arc path in $\mathcal{H}(c)$, with $Z_0 = A_i$ or $B_i$, can be restated in terms of associated regular arcs as follows. Let $[u^k_n, v^k_n]$ be the regular arc on $H_c$ associated to the itinerary arc $J_{Z_k Z_{k+1} \ldots Z_{k +n}}$. Then the last two conditions translate to:
\begin{center}
If $Z_k = Z_0$, then $v^k_n \leq v^0_n$ (or equivalently, $u^k_n \leq u^0_n$) for all $n \geq 0$.
\end{center}

\begin{prop} \label{proposition: characteristic arc words and characteristic points}
If $w$ is a characteristic arc path and $\phi$ is its associated angle, then the landing point of $\R_c(\phi)$ is characteristic. Conversely, if $\R_{c}(\phi)$ lands at a characteristic point $z$ in $[\alpha, c]$, there are three possibilities:
\begin{enumerate}[a)]
    \item there exists a characteristic arc path $w$ whose associated bit sequence is $\phi$;
    \item $z$ is periodic, the landing point of a marked angle, and is part of a satellite orbit portrait;
    \item $\phi$ is an isolated angle landing on $H_c$ and an iterated preimage of a marked angle.
\end{enumerate}
\end{prop}

\begin{proof}
Let $w = Z_{0}Z_{1}\ldots$ be a characteristic arc path (starting at the characteristic state $Z_{0} = A_i$ or $B_i$), $\phi$ its associated angle, and $z$ the landing point of $\R_c(\phi)$. If $z$ were not characteristic, then for some iterate $z_k$ of $z$ we would have $(z, c]$ by Proposition \ref{proposition: characteristic points}. If $x_{i+1} < z_k$ in the linear order on $[\alpha, c]$ from $\alpha$ to $c$, then necessarily $\phi \in (\phi_{i}^2, \phi_{i}^3)$, and $w$ passes through some state in $\mathcal{D}$, which is not possible. Hence $z_n \in (z, x_{i+1}]$, so that $2^k\phi$ is in $A_i$ or $B_i$. Moreover, since $z \neq z_k$, and $\bigcap_{n \geq 0}[u^k_n, v^k_n] = \{z_k\}$ for all $k \geq 0$, there must be disjoint regular arcs $[u^0_n, v^0_n]$ and $[u^k_n, u^k_n]$ such that $v^0_n < u^k_n < v^k_n$ in the linear order on $[\alpha, c]$, a contradiction with $w$ being characteristic.

Now suppose $\R_c(\phi)$ lands at a characteristic point $z \in [x_i, x_{i+1}] \subseteq [\alpha, c]$. First assume that there is some infinite arc path $w = Z_0Z_1\ldots$ in $\mathcal{H}(c)$ that accumulates on $\phi$, and suppose that it is not characteristic. If for some $i \geq k$ we have $Z_k \in \mathcal{D}_i$, the only possibility for $z$ to be characteristic is if $z = z_k = x_{i+1}$, and $w$ accumulates on $w$ from ``below'' in the orientation of $[\alpha, c]$ from $\alpha$ to $c$. But this means that the return map $f^k$ on $z$ maps a point on $[x_i, x_{i+1})$ ``infinitesimally close'' to $x_{i+1}$ to a point on another branch of the Hubbard tree at $x_{i+1}$. This implies that the first return map does not fix the rays landing at $x_{i+1}$, and its orbit portrait is satellite. Similarly, if for some $k \geq 1$ we had $Z_k = z_0$, $v^k_n > v^0_n$ for some $n$, it is necessary that $u^k_n = v^0_n$, otherwise the regular arcs $[u_n^0, v^0_n]$, $[u^k_n, v^k_n]$ are disjoint and $z < z_k$. Moreover, $z = u^k_n = z_k$, so that is it is again periodic. Because an external argument of $z$ is an iterated preiamge of a marked angle, at it itself is periodic, then some (and therefore all) external arguments of $z$ are marked angles, and $z = x_{i+1}$ by the arrangement of $[u_n^0, v^0_n]$, $[u^k_n, v^k_n]$. By the same arguments as before, it has a satellite orbit portrait.

Finally, if there is no infinite arc path $w$ accumulating on $\phi$, then it must be an isolated angle landing on $H_c$, and it iterates to a marked angle.
\end{proof}

\section{An Algorithm for the First Ray Connection Between Hubbard Trees} \label{section: algorithm for first ray connection}

Let $c_1$ and $c_2$ be PCF parameters, and $\mathcal{M}$ a union automaton for $f_{c_1} \sqcup f_{c_2}$. We consider a sufficient enlargement $S'$ of the set of marked angles $S = S_1 \cup (-S_2)$ as in the previous section, and assume this choice throughout. Recall that $\mathcal{M}$ is seen as a subautomaton of $\mathcal{H}(c_1)$, with set of states $\mathcal{S}_{\mathcal{M}} \subseteq \mathcal{S}_{\mathrm{good}}$. Moreover, the involution $\iota : \mathcal{S}_{\mathrm{arc}} \cap \mathcal{S}_{\mathrm{good}}$ on good arc states does not necessarily leave $\mathcal{S}_{\mathcal{M}}$ invariant, and if $w$ is an arc path in $\mathcal{M}$, it is not necessarily the case that $\iota(w)$ also is.

As remarked previously, our goal to prove Theorem \ref{theorem: first ray connection} is to show that $z_0 \in [\alpha, c]$ as defined in Corollary \ref{corollary: exists minimal characteristic} is preperiodic, if it exists. If case 2 of Theorem \ref{theorem: first ray connection} does not hold, Lemma \ref{lemma: z0 strictly preperiodic} guarantees that $z_0$ will be periodic, and will fall into case 3 of the Theorem. The proof of the essential Lemma \ref{lemma: z0 preperiodic} suggests a algorithm that finds $z_0$ and its external arguments, based on the automaton $\mathcal{M}$.

\begin{lemma} \label{lemma: z0 preperiodic}
Let $c_1$ and $c_2$ be PCF parameters. Suppose that there exists a ray connection between $H_{c_1}$ and $H_{c_2}$, and that $c_1$ is not the tip of a vein. Then $z_0$ as in Corollary \ref{corollary: exists minimal characteristic} is preperiodic.
\end{lemma}

\begin{proof}
Consider the partition $\bigcup_{i=0}^{N}[x_i, x_{i+1}]$ of the arc $[\alpha_c, c_1]$ as given in Section \ref{section: characteristic states automaton} if $c_1$ is Misiurewicz, or of the arc $[\alpha_{c_1}, c_1^*]$ if $f_{c_1}$ is critically periodic. Starting from $i = 0$, let $i$ be the current index in the algorithm we now describe. We assume that there exists no characteristic landing point of a direct ray connection between $H_{c_1}$ and $H_{c_2}$ on $[\alpha, x_{i})$.

Let $\mathcal{M}_i$ be the subautomaton of $\mathcal{M}$ obtained by erasing all the states in $\mathcal{S}_{\mathcal{M}} \cap \mathcal{D}_i$ and trimming, where $\mathcal{D}_i$ is as defined in Definition \ref{definition: arc characteristic word}. If $\mathcal{M}_i$ is the empty automaton, meaning that $\mathcal{L}(\mathcal{M}_i) = \emptyset$, then there are no characteristic landing points of direct ray connections between the Hubbard trees on $[x_i, x_{i+1})$. If there were, the angle $\phi$ of such a direct ray connection would be a infinte path in $\mathcal{M}_i$ not passing through the states in $\mathcal{D}_i$, by virtue of the landing point being characteristic and distinct from $x_{i+1}$. If $\mathcal{L}(\mathcal{M}_N) = \emptyset$, so that $\mathcal{L}(\mathcal{M}_i) = \emptyset$ for all $0 \leq i \leq N$, then the only possible characteristic landing point of a direct ray connection is $x_{N + 1} = c_1$ or $c_1^*$, so $z_0 = c_1$ or $c_1^*$ is preperiodic and we are done. Suppose then that $\mathcal{L}(\mathcal{M}_N) \neq \emptyset$, and let $i$ be the first index for which $\mathcal{L}(\mathcal{M}_i) \neq \emptyset$.

Consider the point $x_i$, which is preperiodic, being the landing point of a ray at a preperiodic marked angle. If it characteristic and the landing point of a direct ray connection between $H_{c_1}$ and $H_{c_2}$, then $z_0 = x_i$ is preperiodic and we are done.

Suppose that $x_i$ is not a characteristic landing point of a direct ray connection between the trees. We claim that the set of states $\{A_i, B_i\}$ is a cross section of $\mathcal{M}_i$. If there were an infinite path in $\mathcal{M}_i$ that avoided $\{A_i, B_i\}$, there is a periodic such path, say $w$, with associated angle $\phi$. Since $c_1$ is not the tip of a vein, the landing point of $\R_{c_1}(\phi)$ is biaccessible, being part of a non-trivial orbit portrait. Some element of its orbit is characteristic, and we may assume that $\R_c(\phi)$ lands at this characteristic point. Hence either $w$ must pass through a characteristic state $A_j, B_j$ for $j < i$, or $\R_{c_1}(\phi)$ lands at $x_i$, since $\R_{c_1}(\phi)$ cannot land at $\beta_{c_1} \notin H_{c_1}$. But both cases are not possible by our assumptions; the former would imply that there is a characteristic landing point of a direct ray connection between the trees in $[\alpha, x_i]$, and we assumed $x_i$ is not such a point. As a consequence of $\{A_i, B_i\}$ being a cross section, all states in $\mathcal{M}_i$ are arc states.

Analogously to the discussion following Proposition \ref{proposition: cross section shape of automaton}, the automaton $\mathcal{M}_i$ has the following shape:
% https://q.uiver.app/#q=WzAsMyxbMCwxLCJBX2kiXSxbMiwxLCJCX2kiXSxbMSwwLCJcXG1hdGhjYWx7VH0iXSxbMCwwLCIiLDEseyJvZmZzZXQiOjEsImFuZ2xlIjotOTB9XSxbMCwwLCIiLDEseyJyYWRpdXMiOjEsImFuZ2xlIjotOTB9XSxbMiwwLCIiLDEseyJvZmZzZXQiOjEsImN1cnZlIjoxfV0sWzIsMCwiIiwxLHsib2Zmc2V0IjotMSwiY3VydmUiOjF9XSxbMiwxLCIiLDEseyJvZmZzZXQiOjEsImN1cnZlIjotMX1dLFsyLDEsIiIsMSx7Im9mZnNldCI6LTEsImN1cnZlIjotMX1dLFsxLDEsIiIsMSx7InJhZGl1cyI6MSwiYW5nbGUiOjkwfV0sWzEsMSwiIiwxLHsib2Zmc2V0IjoxLCJhbmdsZSI6OTB9XSxbMCwxLCIiLDEseyJjdXJ2ZSI6LTF9XSxbMCwxLCIiLDEseyJvZmZzZXQiOi0yLCJjdXJ2ZSI6LTF9XSxbMSwwLCIiLDEseyJvZmZzZXQiOi0yLCJjdXJ2ZSI6LTF9XSxbMSwwLCIiLDEseyJjdXJ2ZSI6LTF9XV0=
\[\begin{tikzcd}
	& {\mathcal{T}} & \\
	{A_i} && {B_i}
	\arrow[shift right, curve={height=6pt}, from=1-2, to=2-1]
	\arrow[shift left, curve={height=6pt}, from=1-2, to=2-1]
	\arrow[shift right, curve={height=-6pt}, from=1-2, to=2-3]
	\arrow[shift left, curve={height=-6pt}, from=1-2, to=2-3]
	\arrow[shift right, from=2-1, to=2-1, loop, in=145, out=215, distance=10mm]
	\arrow[from=2-1, to=2-1, loop, in=150, out=210, distance=5mm]
	\arrow[curve={height=-6pt}, from=2-1, to=2-3]
	\arrow[shift left=2, curve={height=-6pt}, from=2-1, to=2-3]
	\arrow[shift left=2, curve={height=-6pt}, from=2-3, to=2-1]
	\arrow[curve={height=-6pt}, from=2-3, to=2-1]
	\arrow[from=2-3, to=2-3, loop, in=330, out=30, distance=5mm]
	\arrow[shift right, from=2-3, to=2-3, loop, in=325, out=35, distance=10mm]
\end{tikzcd}\]
where $\mathcal{T} = \mathcal{M}_{i} \setminus \{A_i, B_i\}$, and the arrows describe finite paths as in Proposition \ref{proposition: cross section shape of automaton}. We will use this to construct a minimal characteristic arc path in $\mathcal{M}_i$ starting at $A_i$, and a maximal characteristic arc path starting at $B_i$, when either paths exist. As a consequence of the construction, they will be preperiodic, corresponding to a preperiodic angle, and the landing point of the ray will be $z_0$.

Invoking the notation of Proposition \ref{proposition: cross section shape of automaton}, for $X, Y \in \{A_i, B_i\}$, let $\mathcal{P}_{XY}$ be the finite set of all finite paths starting at $X$ and ending at $Y$, and not passing through $\{A_i, B_i\}$ execept for the start and end states. Since no path in any of the four collections defined is an extension of any other, these sets, along with the unions $\mathcal{P}_{A_iA_i} \cup \mathcal{P}_{A_i B_i}$ and $\mathcal{P}_{B_i A_i} \cup \mathcal{P}_{B_i B_i}$, are linearly ordered with respect to the lexicographic order on arc paths described in Section \ref{subsection: characteristic states and words}. We also write $\iota(\mathcal{P}_{XY}) \coloneq \{\iota(p) : p \in \mathcal{P}_{XY}\}$, since states are arc states.

A path $p \in \mathcal{P}_{A_iA_i}$ can always be continued to form an infinite characteristic arc path, simply by repeating this path from $A_i$ to $A_i$; this is analogous for $q \in \mathcal{P}_{B_i B_i}$. However, a path $p \in \mathcal{P}_{A_i B_i}$ only admits a continuation into an infinite characteristic arc path if

\begin{equation} \label{characteristic continuation}
\iota(p) \preceq \max(\mathcal{P}_{B_i A_i} \cup \mathcal{P}_{B_i B_i}).
\end{equation}
Indeed, $\iota(p)$ is not an extension or restriction of a path in $\mathcal{P}_{B_i A_i} \cup \mathcal{P}_{B_i B_i}$, and so is comparable with them. If $\iota(p) \succ \max(\mathcal{P}_{B_i A_i} \cup \mathcal{P}_{B_i B_i})$, then any continuation of $p$ starting at $B_i$ cannot be characteristic, and if (\ref{characteristic continuation}) holds, by taking $q \in \mathcal{P}_{B_i A_i} \cup \mathcal{P}_{B_i B_i}$ such that $\iota(p) \preceq q$, either $q \in \mathcal{P}_{B_i B_i}$ and the path $p\overline{q} = pqqq\ldots$ is characteristic, or $q \in \mathcal{P}_{B_i A_i}$ and the path $\overline{pq} = pqpqpq\ldots$ is characteristic. (Here we engage in the abuse of notation of concatenating paths by contracting the end and start states between consecutive ones.)

Let
$$
\begin{array}{rcl}
    p_* & \coloneqq & \min (\mathcal{P}_{A_i A_i} \cup \mathcal{P}_{A_i B_i}), \\
    p_0 & \coloneqq & \min \{ \mathcal{P} \in \mathcal{P}_{A_i A_i} \cup \mathcal{P}_{A_i B_i} \ | \ p \in \mathcal{P}_{A_i A_i}, \text{ or } p \in \mathcal{P}_{A_i B_i} \text { and (\ref{characteristic continuation}) holds} \},
\end{array}
$$
if the latter set is nonempty. We divide into two cases:

\textit{Case 1: $p_0 \in \mathcal{P}_{A_i A_i}$.} Then $p_0 = \min \mathcal{P}_{A_i A_i}$, and the path $\overline{p_0}$ is characteristic and minimal with said property starting from $A_i$. Indeed, if $u$ is an infinite path starting at $A_i$ and $u \prec \overline{p_0}$, then for some iterate where $u_k = A_i$, the path taken by $u$ at this iterate differs from $p_0$ and is smaller in the lexicographic order, which contradicts minimality of $p_0$.

\textit{Case 2: $p_0 \in \mathcal{P}_{A_i B_i}$.} Let
$$
q_0 \coloneqq \min \{ q \in \mathcal{P}_{B_i A_i} \cup \mathcal{P}_{B_i B_i} \ | \ \iota(p_0) \preceq q \},
$$
which exists by assumption. If $q_0 \in \mathcal{P}_{B_i B_i}$, then $p_0 \overline{q_0}$ is characteristic and minimal having this property starting at $A_i$, by analogous arguments as in Case 1: any other choice of path at an iterate $u_k = B_i$ would either result in a larger path in the lexicographic order, or violate the word being characteristic. 

Suppose then that $q_0 \in \mathcal{P}_{B_i A_i}$. Any state $p \prec p_0$ must be in $\mathcal{P}_{A_i B_i}$, and in particular $p_* \in \mathcal{P}_{A_i B_i}$. If $\iota(p_0) \prec q_0$, we claim that the word $p_0q_0\overline{p_*q_0}$ is characteristic and minimal with said property. any smaller choice of path $q \prec q_0$ in the path would prevent it from being characteristic whenever the current state is $B_i$, since by minimality of $q_0$ we would have $\iota(p_0) \succ q$; and the choice of $p_* \preceq p_0$ is possible and minimizing for each state $A_i$ after the first occurence because of the inequality $\iota(p_0) \prec q_0$, which already guarantees that the path is characteristic when passing through $B_i$ without having to check more states.

On the other hand, if $\iota(p_0) = q_0$, then $\overline{p_0q_0}$ is minimal characteristic. As above, the choice of path $q_0$ is optimal when the state is $B_i$. If the current state is $A_i$, we claim that the choice of path $p_0$ is optimal for minimizing the path while still being characteristic. If $p \prec p_0$, then $p \in P_{A_i B_i}$, and the path
$$
p_0q_0p_0q_0\ldots p_0q_0p
$$
is such that, at the last occurence of $q_0$ shown, we have
$$
\iota(p_0q_0) = q_0p_0 \succ q_0p
$$
due to Lemma \ref{lemma: involution comparison}, which prevents the path from being characteristic. In all cases, if a characteristic path exists, the minimal characteristic path is preperiodic, with preperiodic associated angle $\phi_{A_i}$ and landing point $z_{A_i}$ of $\R_{c_1}(\phi_{A_i})$.

By analogous arguments for paths in $\mathcal{M}_i$ starting at $B_i$, but with the lexicographic order reversed, we also find the minimal characteristic path starting at $B_i$, being preperiodic if any exist, and landing point $z_{B_i}$. If both exist, let $z_* = \min(z_{A_i}, z_{B_i})$. As a consequence of Lemma \ref{lemma: lexicographic order on words transfers to angles}, Proposition \ref{proposition: characteristic arc words and characteristic points}, and minimality of $z_0$, we necessarily have $z_0 = z_*$, and we are done.
\end{proof}

\begin{proof}[Proof of Theorem \ref{theorem: first ray connection}]
Suppose that there exists a ray connection between the Hubbard trees $H_{c_1}$ and $H_{c_2}$, and assume first that $c_1$ is not the tip of a PCF. Lemmas \ref{lemma: z0 strictly preperiodic} and \ref{lemma: z0 preperiodic} guarantee that the minimal characteristic landing point $z_0$ of a direct ray connection between $H_{c_1}$ and $H_{c_2}$ is periodic. We claim that the parameter counterpart $c_0$ of $z_0$ satifies the conditions for case 3 in the statement of the Theorem. If not, suppose that for some parameter $c \prec c_0$, possibly not postcritically finite, there were a direct ray connection between the trees $T_c$ and $H_{c_2}$. Lemma \ref{lemma: exists characteristic} implies the existence of a characteristic such ray connection, say at angle $\psi$, which by monotonicity of the characteristic lamination is also present for $f_{c_1} \sqcup f_{c_2}$. But the landing point on $H_{c_1}$ of $\R_{c_1}(\psi)$ corresponds to a characteristic leaf that precedes the one corresponding to $z_0$ in $[\alpha, c_1]$, a contradiction with the minimality of $z_0$.

Now suppose that $c_1$ is the tip of a vein, and that there exists a direct ray connection between the Hubbard trees. If there exists a periodic characteristic one at angle $\psi$, we consider the mating $f_{c_{\psi}} \sqcup f_{c_2}$ and apply the same arguments above, invoking again monotonicity of the characteristic lamination and the parameter-dynamical correspondence of characteristic points. If not, by Proposition \ref{proposition: exists periodic characteristic} we fall into case 2 of the theorem: if there were any ray connections between $T_c$ and $H_{c_2}$ for $c \prec c_1$, Lemma \ref{lemma: exists characteristic} implies some characteristic direct ray connection landing at a biaccessible point in $H_{c}$ will be preserved for $f_{c_1} \sqcup f_{c_2}$. The direct ray connection would be in the grand orbit of $\R_{f_{c_1} \sqcup f_{c_2}}(\theta_1)$, where $\theta_1$ is the unique external argument of $c_1$, and therefore the landing point is biaccessible and precritical. But a precritical point is never characteristic, iterating to $c_1$.
\end{proof}

\begin{remark} \label{remark: infinitely many first ray connections}
We return to the mating $f_{1/4} \sqcup f_{7/8}$ of Example \ref{example: 1/4 with 7/8}. Within the mating automaton, we are able to identify the direct ray connection at angle $.\overline{0100} = 4/15$, corresponding to the unique hyperbolic component of period $4$ in the vein to $1/4$. By considering the mating $f_{4/15} \sqcup f_{7/8}$, it has a mating automaton of the form below, after a suitable trimming preserving the regular language:
% https://q.uiver.app/#q=WzAsMTMsWzIsMCwiMCJdLFszLDAsIjAiXSxbMywxLCIwIl0sWzUsMCwiMCJdLFs0LDEsIjEiXSxbMiwxLCIxIl0sWzEsMSwiMCJdLFswLDEsIjEiXSxbMCwwLCIwIl0sWzEsMCwiMSJdLFs1LDEsIjAiXSxbNiwxLCIwIl0sWzYsMCwiMSJdLFswLDFdLFsxLDJdLFsyLDRdLFsyLDVdLFs1LDBdLFs1LDZdLFs2LDddLFs3LDhdLFs4LDldLFs5LDZdLFs0LDEwXSxbMTAsMTFdLFsxMSwxMl0sWzEyLDNdLFszLDEwXV0=
\[\begin{tikzcd}
	0 & 1 & 0 & 0 && 0 & 1 \\
	1 & 0 & 1 & 0 & 1 & 0 & 0
	\arrow[from=1-1, to=1-2]
	\arrow[from=1-2, to=2-2]
	\arrow[from=1-3, to=1-4]
	\arrow[from=1-4, to=2-4]
	\arrow[from=1-6, to=2-6]
	\arrow[from=1-7, to=1-6]
	\arrow[from=2-1, to=1-1]
	\arrow[from=2-2, to=2-1]
	\arrow[from=2-3, to=1-3]
	\arrow[from=2-3, to=2-2]
	\arrow[from=2-4, to=2-3]
	\arrow[from=2-4, to=2-5]
	\arrow[from=2-5, to=2-6]
	\arrow[from=2-6, to=2-7]
	\arrow[from=2-7, to=1-7]
\end{tikzcd}\]
which shows that the only direct ray connections are preperiodic, and eventually fall into the cycle of $4/15$. Hence the parameter $c_1'$ as in case 3 of Theorem \ref{theorem: first ray connection} is $c_{4/15}$, and all ray connections between the Hubbard trees are destroyed once we move back from this parameter along the vein $V_{1/4}$.
\end{remark}

\begin{proof}[Proof of Corollary \ref{corollary: existence of topological mating}]
Let $c_1' \in V_{c_1}$ and $c_2' \in V_{c_2}$. By monotonicity of the characteristic lamination and Lemma \ref{lemma: exists characteristic}, there are no ray connections between the trees $T_{c_1'}$ and $T_{c_2'}$. Proposition 4.1 in \cite{Jung2017} asserts that the length of ray connections between $J_{c_1'}$ and $J_{c_2'}$ is uniformly bounded; iterating any ray connection passing through biaccessible points enough times makes it pass through the generalized Hubbard trees. Moreover, there are no cyclic ray connections. This can be seen from the Proposition itself or as an application of Proposition 4.4 in the same paper, where one would exist if and only if there is a direct ray connection between the critical values $c_1$ and $c_2$. By a result of Epstein (\cite[Proposition 2.5]{Jung2017}), the topological mating exists.
\end{proof}

\section{The Disjointness Locus} \label{section: disjointness locus}

Let $V_1$ and $V_2$ be two PCF veins of $M$, with tips $t_1$ and $t_2$ respectively; we will denote $V_i = [0, t_i]$ as a regular arc in $M$, oriented from the center $c = 0$ of the main cardioid to $t_i$. Recall that the \textit{disjointness locus} $D$ of $V_1$ and $V_2$, defined in the introduction, is the set of parameter pairs for which there are no ray connections between the generalized trees in the formal mating. Also let $E = (V_1 \times V_2) \setminus D$ be its complement (that is, the set of parameter pairs for which there exists a ray connection between the generalized trees).

If there is a ray connection between the Hubbard trees $H_{t_1}$ and $H_{t_2}$ whose angle (with respect to $f_{t_1}$) is not in the grand orbit of the external argument of $t_1$, Theorem \ref{theorem: first ray connection} guarantees that there exists a parabolic parameter $p_1 \in (0, t_1]$ such that $E \cap (V_1 \times \{t_2\}) = [p_1, t_1] \times \{t_2\}$, and if $P_M(\theta_1^-, \theta_1^+)$ is the ray pair landing at $p_1$, either $\R_{t_2}(-\theta_1^-)$ or $\R_{t_2}(-\theta^+)$ lands at the Hubbard tree $H_{t_2}$. If the ray does not land at an extremity of $H_{t_2}$, it lands at a biaccessible periodic point, and the characteristic point of this orbit determines a wake in the parameter plane for which the direct ray connection between the trees is preserved. More precisely, there exists $p_2^* \in V_2$ such that $[p_1, t_1] \times [p_2^*, t_2] \in E$. This is because for any $c_1 \in [p_1, t_1]$ and $c_2 \in [p_2^*, t_2]$, either $\R_{f_{c_1} \sqcup f_{c_2}}(\theta_1^-)$ or $\R_{f_{c_1} \sqcup f_{c_2}}(\theta_1^+)$ is a direct ray connection between the generalized trees $T_{c_1}$ and $T_{c_2}$, as $c_1$ is in the wake determined by (the hyperbolic component of) $p_1$ and $c_2$ is in the wake determined by $p_2^*$, so the orbit portraits are realized and land on the trees.

Analogously, if $p_2$ is the (parabolic) first parameter on $V_2$ such that $(t_1, p_2) \in E$, and the corresponding direct ray connection does not land at an extremity of $H_{t_1}$, then there exists a parabolic parameter $p_1^*$ such that $[p_1^*, t_1] \times [p_2, t_2] \in E$. Note that necessarily $p_1 \preceq p_1^*$, $p_2 \preceq p_2^*$.

\begin{exm} \label{example: varying both paramters, 1/4 with 1/6}

Consider the mating $f_{1/4} \sqcup f_{1/6}$ as in Example \ref{example: 1/4 with 1/6}. In this case, $p_1 = p_1^*$, being the landing point of the parameter ray pair $P_M(7/31, 8/31)$, and $p_2 = p_2^*$ is the landing point of $P_M(5/31, 6/31)$. Moreover, $E = [p_1, t_1] \times [p_2, t_2]$.
\end{exm}

\begin{exm}
Consider now the mating $f_{1/4} \sqcup f_{7/8}$ as in Example \ref{example: 1/4 with 7/8}. According to Remark \ref{remark: infinitely many first ray connections}, $p_1$ is the landing point of the parameter ray pair $P_M(3/15, 4/15)$, while $p_2^*$ is the landing point of $P_M(13/15, 14/15)$. A mating automaton for $f_{1/4} \sqcup f_{13/15}$ is given below:
% https://q.uiver.app/#q=WzAsNCxbMCwwLCIwIl0sWzEsMCwiMCJdLFsxLDEsIjEiXSxbMCwxLCIwIl0sWzAsMV0sWzEsMl0sWzIsM10sWzMsMF1d
\[\begin{tikzcd}
	0 & 0 \\
	0 & 1
	\arrow[from=1-1, to=1-2]
	\arrow[from=1-2, to=2-2]
	\arrow[from=2-1, to=1-1]
	\arrow[from=2-2, to=2-1]
\end{tikzcd}\]
showing that $p_2 = p_2^*$, $p_1 = p_1^*$, and $E = [p_1, t_1] \times [p_2, t_2]$ as well, even though there are uncoutably many periodic ray connections between $H_{t_1}$ and $H_{t_2}$.
\end{exm}

\begin{exm} \label{example: 1/4 with 3/8}
Consider the mating automaton for $f_{1/4} \sqcup f_{3/8}$ below:
% https://q.uiver.app/#q=WzAsMTEsWzAsMSwiMCJdLFsxLDEsIjEiXSxbMiwxLCIwIl0sWzMsMSwiMCJdLFs0LDEsIjEiXSxbNCwyLCIwIl0sWzMsMiwiMSJdLFsyLDIsIjAiXSxbNCwwLCIxIl0sWzUsMCwiMSJdLFs1LDEsIjAiXSxbMCwxXSxbMSwyXSxbMiwzXSxbMyw0XSxbNCw1XSxbNSw2XSxbNiw3XSxbNywyXSxbNyw3LCIiLDAseyJyYWRpdXMiOjEsImFuZ2xlIjotMTM1fV0sWzMsOF0sWzgsNF0sWzksOF0sWzEwLDldLFs1LDEwXSxbOSw5LCIiLDAseyJyYWRpdXMiOjEsImFuZ2xlIjo0NX1dXQ==
\[\begin{tikzcd}
	&&&& 1 & 1 \\
	0 & 1 & 0 & 0 & 1 & 0 \\
	&& 0 & 1 & 0
	\arrow[from=1-5, to=2-5]
	\arrow[from=1-6, to=1-5]
	\arrow[from=1-6, to=1-6, loop, in=15, out=75, distance=5mm]
	\arrow[from=2-1, to=2-2]
	\arrow[from=2-2, to=2-3]
	\arrow[from=2-3, to=2-4]
	\arrow[from=2-4, to=1-5]
	\arrow[from=2-4, to=2-5]
	\arrow[from=2-5, to=3-5]
	\arrow[from=2-6, to=1-6]
	\arrow[from=3-3, to=2-3]
	\arrow[from=3-3, to=3-3, loop, in=195, out=255, distance=5mm]
	\arrow[from=3-4, to=3-3]
	\arrow[from=3-5, to=2-6]
	\arrow[from=3-5, to=3-4]
\end{tikzcd}\]

By inspection, the angle $.\overline{001010} = 10/63$ is recognized, and necessarily $p_1 = c_{10/63}$. This is because the only other parabolic parameter preceding it in the vein to $c_{1/4}$ is $c_{1/7}$, and neither $1/7$ or $2/7$ forms a ray connection between the trees. Hence $p_2^* = c_{23/63}$, as $c_{3/8}$ is in the wake determined by the parameter ray pair $P_M(23/63, 24/63)$.

On the other hand, the angle $.\overline{10011} = 19/31$ is recognized, so that $\R_{c_{3/8}}(-19/31) = \R_{c_{3/8}}(12/31)$ lands on $H_{c_{3/8}}$ and $c_{3/8}$ is in the wake of $P_M(11/31, 12/31)$. The mating $f_{1/4} \sqcup f_{12/31}$ has the automaton below:
% https://q.uiver.app/#q=WzAsMTAsWzAsMCwiMCJdLFsxLDAsIjEiXSxbMiwwLCIwIl0sWzMsMCwiMCJdLFszLDEsIjAiXSxbNCwwLCIxIl0sWzQsMSwiMSJdLFs1LDEsIjAiXSxbNiwxLCIwIl0sWzUsMCwiMSJdLFswLDFdLFsxLDJdLFsyLDNdLFszLDVdLFs1LDZdLFs0LDZdLFs2LDddLFs3LDhdLFs4LDldLFs5LDVdXQ==
\[\begin{tikzcd}
	0 & 1 & 0 & 0 & 1 & 1 & \\
	&&& 0 & 1 & 0 & 0
	\arrow[from=1-1, to=1-2]
	\arrow[from=1-2, to=1-3]
	\arrow[from=1-3, to=1-4]
	\arrow[from=1-4, to=1-5]
	\arrow[from=1-5, to=2-5]
	\arrow[from=1-6, to=1-5]
	\arrow[from=2-4, to=2-5]
	\arrow[from=2-5, to=2-6]
	\arrow[from=2-6, to=2-7]
	\arrow[from=2-7, to=1-6]
\end{tikzcd}\]
which shows that all ray connections are in the grand orbit of $19/31 = -12/31$. Moving back along the vein of $c_{3/8}$ before $c_{12/31}$ would destroy all ray connections, and so $p_2 = c_{12/31} \prec c_{23/63}$. Moreover, $p_1^* = c_{7/31} \succ c_{10/63} = p_1$. See Figure \ref{fig: vein square}.

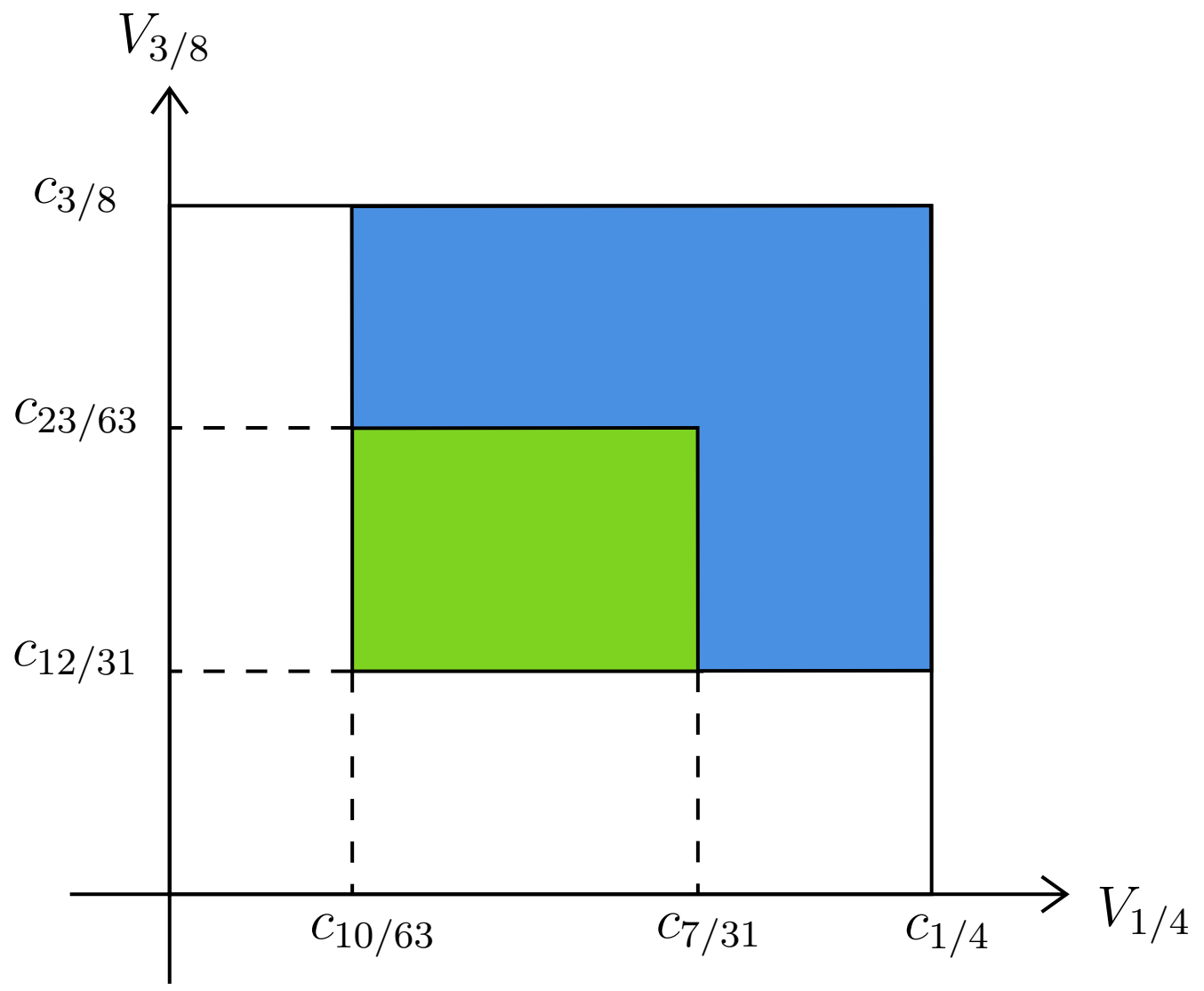
\begin{figure}[ht]

\centering

\tikzset{every picture/.style={line width=0.75pt}} %set default line width to 0.75pt        

\begin{tikzpicture}[x=0.75pt,y=0.75pt,yscale=-1,xscale=1]
%uncomment if require: \path (0,309); %set diagram left start at 0, and has height of 309

%Shape: Axis 2D [id:dp5040057837075325] 
\draw  (164.89,264.75) -- (446.13,264.75)(193.01,37.52) -- (193.01,290) (439.13,259.75) -- (446.13,264.75) -- (439.13,269.75) (188.01,44.52) -- (193.01,37.52) -- (198.01,44.52)  ;
%Shape: Rectangle [id:dp662204819167594] 
\draw   (193.01,70.54) -- (408.05,70.54) -- (408.05,264.75) -- (193.01,264.75) -- cycle ;
%Shape: Rectangle [id:dp5962221206941413] 
\draw  [fill={rgb, 255:red, 126; green, 211; blue, 33 }  ,fill opacity=1 ] (244.56,133.19) -- (343.2,133.19) -- (343.2,201.85) -- (244.56,201.85) -- cycle ;
%Straight Lines [id:da8531923287782177] 
\draw  [dash pattern={on 4.5pt off 4.5pt}]  (244.56,201.85) -- (244.56,264.75) ;
%Shape: L Shape [id:dp8496731361066926] 
\draw  [fill={rgb, 255:red, 74; green, 144; blue, 226 }  ,fill opacity=1 ] (407.98,201.68) -- (342.1,201.75) -- (342.03,133.19) -- (244.49,133.29) -- (244.42,70.71) -- (407.84,70.54) -- cycle ;
%Straight Lines [id:da8670395969608706] 
\draw  [dash pattern={on 4.5pt off 4.5pt}]  (342.1,201.75) -- (342.1,243) -- (342.1,264.65) ;
%Straight Lines [id:da09624623592264414] 
\draw  [dash pattern={on 4.5pt off 4.5pt}]  (244.56,201.85) -- (193.01,201.85) ;
%Straight Lines [id:da9596121131976305] 
\draw  [dash pattern={on 4.5pt off 4.5pt}]  (244.56,133.19) -- (193.01,133.19) ;

% Text Node
\draw (453.11,261.15) node [anchor=north west][inner sep=0.75pt]    {$V_{1/4}$};
% Text Node
\draw (176.5,14.76) node [anchor=north west][inner sep=0.75pt]    {$V_{3/8}$};
% Text Node
\draw (152.86,61.31) node [anchor=north west][inner sep=0.75pt]    {$c_{3/8}$};
% Text Node
\draw (398.9,269.04) node [anchor=north west][inner sep=0.75pt]    {$c_{1/4}$};
% Text Node
\draw (147.16,191.88) node [anchor=north west][inner sep=0.75pt]    {$c_{12/31}$};
% Text Node
\draw (147.36,123.57) node [anchor=north west][inner sep=0.75pt]    {$c_{23/63}$};
% Text Node
\draw (231.08,268.93) node [anchor=north west][inner sep=0.75pt]    {$c_{10/63}$};
% Text Node
\draw (328.78,268.7) node [anchor=north west][inner sep=0.75pt]    {$c_{7/31}$};

\end{tikzpicture}
\caption{Depiction of parameter pairs in $V_{1/4} \times V_{3/8}$. The open white region is in the disjointness locus $D$, while the blue region is in the complement $E$. At first glance, it is unclear how the green rectangle is split into $D$ and $E$.}
\label{fig: vein square}
\end{figure}
\end{exm}

Note that the parameter $c_1$ in the lemma below need not be PCF.

\begin{lemma} \label{lemma: closed arc in veins}
Let $V_1$ and $V_2$ be PCF veins, and $c_1 \in V_1$. Then the set of parameters $c_2 \in V_2$ such that there is a ray connection between the generalized Hubbard trees $T_{c_1}$ and $T_{c_2}$ in the formal mating $f_{c_1} \sqcup f_{c_2}$ is empty, or a closed arc containing the tip of $V_2$. 
\end{lemma}

\begin{proof}
Let $S \subseteq V_2$ be the set of parameters in the statement. We claim that if $c_2 \in S$ and $c_2 \prec c$, then $c \in S$. If $c_2$ is the tip of a vein, there is nothing to prove. If not, Lemma \ref{lemma: exists characteristic} implies the existence of a direct ray connection between $T_{c_1}$ and $T_{c_2}$ that is characteristic for $c_2$. Monotonicity of the characteristic lamination guarantees that this ray connection is also present in $f_{c_1} \sqcup f_{c}$, and so $c \in S$.

Therefore $S$ is empty or an arc containing $t_2$, and it remains to prove that it is closed. If it is empty or contains just the tip of $V_2$, there is nothing to prove, so assume there exists $c_2 \in S$ not the tip. Lemma \ref{lemma: exists characteristic} and Corollary \ref{corollary: exists minimal characteristic} imply the existence of a minimal characteristic landing point $z_0 \in [\alpha_2, c_2]$ of a direct ray connection between $T_{c_1}$ and $T_{c_2}$. If $z_0$ is preperiodic, by Lemma \ref{lemma: z0 strictly preperiodic} it must be periodic and Theorem \ref{theorem: first ray connection} proves the statement, so assume $z_0$ is not preperiodic.

Let $\phi$ be the irrational angle of this ray connection (with respect to the dynamical plane of $c_2$). Proposition \ref{proposition: irrational point on vein} asserts that there is a unique $c_0 \in V_t$ such that $\R_{c_0}(\phi)$ lands at the critical value $c_0$ in the dynamical plane of $c_0$, and the minor leaf $l_{c_0}$ is of the form $(\phi, \phi')$. In particular, there is a ray connection between $T_{c_1}$ and $T_{c_0}$, so $c_0 \in S$. We claim that $c_0 = \min S$.

Suppose that $c \prec c_0$, and that there is a ray connection between the trees $T_{c_1}$ and $T_c$. Again by Lemma \ref{lemma: exists characteristic}, there is a direct ray connection that is characteristic for $c$, whose angle is an endpoint of some characteristic leaf $l_c$. Note that $l_c \prec l_{c_0}$, and this characteristic direct ray connection is preserved for $T_{c_1}$ and $T_{c_2}$ by monotonicity of the characteristic lamination. But this implies that there is a characteristic landing point of a direct ray connection $z \prec z_0$ in the lamination $L_{c_2}$, a contradiction. Hence $c \notin S$, $c_0 = \min S$, and $S$ is closed.
\end{proof}

\begin{proof}[Proof of Theorem \ref{theorem: disjointness locus open}]
We will show that the set $E \subseteq V_1 \times V_2$ is closed. We may assume that $E$ contains a point $(x,y) \in [0, t_1) \times [0, t_2)$, otherwise the result follows immediately from Lemma \ref{lemma: closed arc in veins}.

Let $(x_n, y_n) \in E$ be a sequence of parameter pairs such that $(x_n, y_n) \to (x, y)$. Lemma \ref{lemma: closed arc in veins} implies that each horizontal segment $E \cap (V_1 \times \{c_2\})$ and each vertical segment $E \cap (\{c_1\} \times V_2)$ is (identified with) a closed arc in $V_1$ (respectively $V_2)$ containing the tip of the vein. Let $x' \succ x$, and $y_{x'} = \min (C \cap (\{x'\} \times V_2)$. Suppose by way of contradiction that $y_{x'} \succ y$. Let $n$ be large enough so that $x_n \prec x'$, $y_n \prec y_{x'}$. Since $(x_n, y_n) \in E$, we have $[x_n, t_1] \times \{y_n\} \subset E$. But then $(x', y_n) \in E$, contradicting minimality of $y_{x'}$. Therefore $y_{x'} \preceq y$, and in particular $(x', y) \in E$. Since this holds for all $x' \succ x$, we have $(x, v_1] \times \{y\} \in E$. By again applying Lemma \ref{lemma: closed arc in veins}, we must have $[x,v_1] \times \{y\} \subset C$, and in particular $(x,y) \in E$.
\end{proof}

\section{Proof of Theorem \ref{theorem: tuned rabbits}} \label{section: proof of tuned rabbits}

We recall some information about tunings and renormalizations (See \cite{Douady1986, McMullen1994, Milnor2000}). If $c \in M_{p/q}$ has locally connected Julia set, then $c = r_{p/q} * c^*$, where $r_{p/q}$ is the parameter of the critically periodic $p/q$-rabbit polynomial and $c^* \in M$ also has connected Julia set. The polynomial $f_c$ is (simply) renormalizable with period $q$, and the quadratic-like map $f_c^q$ near the critical point is hybrid-conjugate to $f_{c^*}$.

The filled Julia set $K_{c}$ contains a bouquet of $q$ homeomorphic copies $K^0, \ldots, K^{q-1}$ of the filled Julia set $K_{c^*}$, all intersecting at their respective $\beta$ fixed points. The common point of intersection is the $\alpha$ fixed point $\alpha_c$ of $K_c$, on which $p/q$ rays land and are permutted transitively among themselves with rotation number $p/q$ under the dynamics of $f_c$. These homeomorphic copies $K^{i}$ are called the \textit{small Julia sets} of $K_c$, and are mapped one to the next with combinatorial rotation number $p/q$ around $\alpha_c$. The map $K^0 \to K^1$ is a branched double cover, and $K^i \to K^{i+1}$ is a conformal isomorphism for $i \neq 0$. 

The main idea for the proof of Therem \ref{theorem: tuned rabbits} is that a direct ray connection between small Julia sets would have to follow their rotational dynamics, and this would imply that some of the rays cross over each other if the rotations don't match.

\begin{lemma} \label{lemma: indirect implies direct rabbit}
Let $c_1 \in M_{p_1/q_1}$ and $c_2 \in M_{p_2/q_2}$ have locally connected Julia sets. If there exists a ray connection between the small Julia sets in the formal mating, then there exists a direct such ray connection.
\end{lemma}

\begin{proof}
Any indirect ray connection between the small Julia sets must pass through a biaccessible point, say in $K_{c_1}$. Every biaccessible point eventually iterates into the Hubbard tree $H_{c_1}$, and the Hubbard tree for tunings of $r_{p_1/q_1}$ is a bouquet of $q_1$ homeomorphic copies of the \textit{extended} Hubbard trees $\widehat{H}_{c^*}$, which is the regular hull of postcritical set along with the $\beta$ and $-\beta$ points of $f_{c^*}$. Each copy is the regular hull of $P_{c_1} \cap K^i_{c_1}$ along with $\alpha_{c_1}$ (which corresponds to the $\beta$ fixed point for $f_{c^*}$) and a $q_1$-th preiamge of $\alpha_{c_1}$ (which corresponds to the $-\beta$ point for $f_{c^*}$). Each copy is contained in the respective small Julia set; therefore by iterating the indirect ray connection sufficiently many times we find a direct one.
\end{proof}

\begin{proof}[Proof of Theorem \ref{theorem: tuned rabbits}]
In the formal mating $f_{c_1} \sqcup f_{c_2}: S^2_{f_{c_1}, f_{c_2}} \to S^2_{f_{c_1}, f_{c_2}}$, we blow up the Carathéodory boundaries of the filled Julia sets $K_{c_1}$ and $K_{c_2}$ to circles. We obtain a closed cylinder $C$ for which the top boundary $S^1_{\mathrm{top}}$ corresponds angles of rays landing on $K_{c_1}$, the bottom boundary $S^1_{\mathrm{bot}}$ corresponds to angles of rays landing on $K_{c_2}$, and direct ray connections between the filled Julia sets are disjoint vertical paths. Alternatively, we take two copies of $S^1$, with opposite orientations, and connect them by an open cylinder. This space has a closed quotient onto $S^2_{f_{c_1}, f_{c_2}} \setminus ( \mathrm{int}K_{c_1} \cup \mathrm{int}K_{c_2})$ which is a homeomorphism on the open cylinder and the quotient induced by the Carathéodory loops on the boundary. The forward orbit of $\phi$ on $S^1_{\mathrm{top}}$ is connected by disjoint paths to the forward orbit of $-\phi$ on $S^1_{\mathrm{bot}}$, which are the preimages of direct ray connections under the quotient. By a slight abuse of notation, we call these paths rays, and say that one such path is mapped to another under the dynamics if the corresponding direct ray connections are mapped accordingly.

We in fact work with a slight modification of this cylinder, still denoted by $C$, obtained by blowing up each external argument of $\alpha_{c_1}$ on $S^1_{\mathrm{top}}$ and each external argument of $\alpha_{c_2}$ on $S^1_{\mathrm{bot}}$ to a small arc, and blowing up each corresponding ray on the open cylinder to a rectangle foliated by vertical paths. Strictly speaking, this makes us lose the dynamics on $C$, but the combinatorics of mappings coming from the orbit of rays is preserved.

To arrive at a contradiction, by Lemma \ref{lemma: indirect implies direct rabbit} it is sufficient to assume that there exists a direct ray connection $\R_{f_{c_1} \sqcup f_{c_2}}(\phi)$ in the formal mating between a small Julia set of $K_{c_1}$ and a small Julia set of $K_{c_2}$. If neither endpoint of the ray connection is the $\alpha$ fixed point or the preimage $-\alpha$ point of the respective polynomial, then the dynamics of $\phi$ under doubling must follow the dynamics of mappings between the small Julia sets. Therefore its forward orbit can be partitioned into $q_1$ unliked sets $U_0, \ldots, U_{q_1 - 1}$ of the circle, permuted transitively with rotation number $p_1/q_1$, corresponding to the angles whose external rays land at each small filled Julia set $K^i_{c_1}$. The analogous statements hold for $c_2$, with $q_2$ unlinked sets $V_0, \ldots, V_{q_2 - 1} \subseteq S^1$ partitioning the forward orbit of $-\phi$. In fact, we may find these unlinked sets even in the cases excluded above, since its external arguments have the same dynamics.

On the circle $S^1$, it could be the case that the closures $U_{i}$ intersect (at one of the external arguments for $\alpha_{c_1}$). However, the construction of $C$ by blowing up points to arcs guarantees that we can find disjoint closed arcs $I_0, \ldots, I_{q_1 - 1}$ on $S^1_{\mathrm{top}}$ and disjoint closed arcs $J_0, \ldots, J_{q_2 - 1}$ on $S^1_{\mathrm{bot}}$ such that $U_i \subseteq I_{i}$ and $V_{j} \subseteq J_j$, under the natural identifications of the blown up $S^1$ with $S^1_{\mathrm{top}}$ and $S^1_{\mathrm{bot}}$. Let $E(i,j)$ be the set of rays with endpoints on $U_{i}$ and $V_{j}$. Then $E(i,j) \neq \emptyset$ if and only if there exists a direct ray connection in the forward orbit of $\R_{f_{c_1} \sqcup f_{c_2}}(\phi)$ between $K^i_{c_1}$ and $K^j_{c_2}$.

Let $C'$ be the closed quotient of $C$ obtained by collapsing each of the closed arcs $I_i$ and $J_j$ to points $u_i$ and $v_j$, respectively. The quotient is still homeomorphic to a cylinder, and all of the paths in $E(i,j)$ now have a common endpoint on each boundary and are homotopic modulo endpoints. We denote this homotopy class by $\mathcal{E}(i,j)$. The homotopy classes $\mathcal{E}(i,j)$ for which $E(i,j) \neq \emptyset$ are disjoint, that is, they have disjoint representatives on the quotient cylinder, since they are represented by non-crossing direct ray connections between the small Julia sets. We say that $\mathcal{E}(i,j)$ \textit{is or is not present}, corresponding to whether $E(i,j) \neq \emptyset$ or not.

Due to the dynamics of the orbit of $\R_{f_{c_1} \sqcup f_{c_2}}(\phi)$, rays in $E(i,j)$ are mapped to rays in $E(i+1,j+1)$, modulo $(q_1, q_2)$. This implies that if $\mathcal{E}(i,j)$ is present, then $\mathcal{E}(i+1,j+1)$ also is. Moreover, since small Julia sets are permutted transitively, for every $i$ there will be some $j$ such that $\mathcal{E}(i,j)$ is present, and for every $j$ there will be some $i$ such that $\mathcal{E}(i,j)$ is present. Let $\mathcal{E}(i, \bullet)$ be the collection of all present homotopy classes having an endpoint on $u_i$, and we analogously define $\mathcal{E}( \bullet, j)$. Assume without loss of generality that $q_1 \leq q_2$.

We claim that $q_2$ must be a multiple of $q_1$. If not, and if $\mathcal{E}(i,j)$ is a present homotopy class with endpoint in $v_j$, then $\mathcal{E}(i + q_2, j)$ is also present and distinct from $\mathcal{E}(i,j)$. Hence $|\mathcal{E}(\bullet, j)| \geq 2$, and this holds for all $j = 0, \ldots, q_2 - 1$. We describe a cell structure for the quotient cylinder $C'$. The $0$-cells are the points $u_i$ and $v_j$; the $1$-cells are the circle arcs connecting the $0$-cells on the boundary and, for each present homotopy class $\mathcal{E}(i,j)$, one representative of it; and the $2$-cells are the faces that fill in the open cylinder with boundaries on the $1$-cells. There are
$$
V = q_1 + q_2
$$
$0$-cells, and
$$
E \geq q_1 + 3q_2
$$
$1$-cells, accounting for the $q_1 + q_2$ arcs on the boundary circles, and at least $2q_2$ present homotopy classes with endpoints on the points $v_j$, due our assumption that $q_2$ is not a multiple of $q_1$. Since the Euler characteristic of $C'$ is $0$, we have
$$
F = E - V \geq 2q_2
$$
$2$-cells. However, each $2$-cell has at least one boundary on a circle arc, since the present homotopy classes are distinct, and so $F \leq q_1 + q_2 \leq 2q_2$. Equality is attained only if $q_1 = q_2$, a contradiction with the assumption that $q_2$ is not a multiple of $q_1$.

Therefore $q_1 | q_2$. With respect to iterates in multiples of $q_1$, if $\mathcal{E}(i,j)$ is present, then
$$
\mathcal{E}(i,j + 1q_1), \  \ldots, \ \mathcal{E}(i,j + \dfrac{q_2}{q_1}q_1) = \mathcal{E}
(i,j)$$
are all present and distinct. Hence $|\mathcal{E}(i, \bullet)|$ is a multiple of $q_2/q_1$ for each $i$. In fact, the number $|E(i, \bullet)|$ is independent of $i$, since for every $\mathcal{E}(i,j)$ present, $\mathcal{E}(i+1, j+1)$ is also present. Therefore the sum $\sum_{i = 0}^{q_1 - 1}|E(i, \bullet)|$ is a multiple of $q_2/q_1$. If the sum is not $q_2/q_1$, it is at least $2q_2/q_1$. Hence, with the same cell structure for $C'$, $E \geq q_1 \cdot 2q_2/q_1 = 2q_2$, which again results in a contradiction. Hence for all $i$ we have $|E(i, \bullet)| = q_2 / q_1$.

For a given $i$, let $W'_i$ be the collection of bottom endpoints of homotopy classes in $\mathcal{E}(i, \bullet)$. The sets $W'_i$ have $q_2/q_1$ elements, and we claim that they are pairwise unlinked on $S^1_{\mathrm{bot}}$ (more precisely, on its quotient in $C'$). This is a consequence of the fact that if $v_{j_1}$ and $v_{j_2}$ are two points in $W'_i$, exactly one of the two open arcs determined by them on $S^1_{\mathrm{bot}}$ has the property that any point $v_{t}$ within it must also be in $W'_i$, that is, connect to $u_i$. If it didn't, and connected to another point $u_{s}$, the homotopy class $\mathcal{E}(s,t)$ would intersect either $\mathcal{E}(i, j_1)$ or $\mathcal{E}(i, j_2)$. See figure \ref{fig: unlinkedness}.

\begin{figure}[ht]
    \centering
\tikzset{every picture/.style={line width=0.75pt}} %set default line width to 0.75pt        

\begin{tikzpicture}[x=0.75pt,y=0.75pt,yscale=-1,xscale=1]
%uncomment if require: \path (0,265); %set diagram left start at 0, and has height of 265

%Shape: Can [id:dp24511313417958325] 
\draw   (478.46,66.68) -- (478.46,189.6) .. controls (478.46,204.15) and (398.89,215.94) .. (300.73,215.94) .. controls (202.57,215.94) and (123,204.15) .. (123,189.6) -- (123,66.68) .. controls (123,52.13) and (202.57,40.34) .. (300.73,40.34) .. controls (398.89,40.34) and (478.46,52.13) .. (478.46,66.68) .. controls (478.46,81.23) and (398.89,93.02) .. (300.73,93.02) .. controls (202.57,93.02) and (123,81.23) .. (123,66.68) ;
%Shape: Ellipse [id:dp7906964337376162] 
\draw  [fill={rgb, 255:red, 0; green, 0; blue, 0 }  ,fill opacity=1 ] (256.92,92.3) .. controls (256.92,90.81) and (258.64,89.6) .. (260.77,89.6) .. controls (262.89,89.6) and (264.61,90.81) .. (264.61,92.3) .. controls (264.61,93.79) and (262.89,95) .. (260.77,95) .. controls (258.64,95) and (256.92,93.79) .. (256.92,92.3) -- cycle ;
%Shape: Ellipse [id:dp3947218045470786] 
\draw  [fill={rgb, 255:red, 0; green, 0; blue, 0 }  ,fill opacity=1 ] (300.73,215.94) .. controls (300.73,214.45) and (302.45,213.24) .. (304.58,213.24) .. controls (306.7,213.24) and (308.42,214.45) .. (308.42,215.94) .. controls (308.42,217.43) and (306.7,218.64) .. (304.58,218.64) .. controls (302.45,218.64) and (300.73,217.43) .. (300.73,215.94) -- cycle ;
%Shape: Ellipse [id:dp8218565939652884] 
\draw  [fill={rgb, 255:red, 0; green, 0; blue, 0 }  ,fill opacity=1 ] (383.72,212.3) .. controls (383.72,210.81) and (385.44,209.6) .. (387.57,209.6) .. controls (389.69,209.6) and (391.41,210.81) .. (391.41,212.3) .. controls (391.41,213.79) and (389.69,215) .. (387.57,215) .. controls (385.44,215) and (383.72,213.79) .. (383.72,212.3) -- cycle ;
%Shape: Ellipse [id:dp56515318892197] 
\draw  [fill={rgb, 255:red, 0; green, 0; blue, 0 }  ,fill opacity=1 ] (218.45,213.3) .. controls (218.45,211.81) and (220.18,210.6) .. (222.3,210.6) .. controls (224.43,210.6) and (226.15,211.81) .. (226.15,213.3) .. controls (226.15,214.79) and (224.43,216) .. (222.3,216) .. controls (220.18,216) and (218.45,214.79) .. (218.45,213.3) -- cycle ;
%Shape: Ellipse [id:dp6701246576558068] 
\draw  [fill={rgb, 255:red, 0; green, 0; blue, 0 }  ,fill opacity=1 ] (396.54,89.3) .. controls (396.54,87.81) and (398.26,86.6) .. (400.39,86.6) .. controls (402.51,86.6) and (404.23,87.81) .. (404.23,89.3) .. controls (404.23,90.79) and (402.51,92) .. (400.39,92) .. controls (398.26,92) and (396.54,90.79) .. (396.54,89.3) -- cycle ;
%Curve Lines [id:da9514734883875805] 
\draw    (226.15,213.3) .. controls (283.14,183.3) and (203.78,122.3) .. (260.77,92.3) ;
%Curve Lines [id:da946965455434914] 
\draw    (387.57,212.3) .. controls (337.42,163.6) and (267.61,153.6) .. (260.77,92.3) ;
%Curve Lines [id:da15321733385416547] 
\draw  [dash pattern={on 4.5pt off 4.5pt}]  (304.58,215.94) .. controls (354.5,195.6) and (351.5,139.6) .. (400.39,89.3) ;
%Curve Lines [id:da18185812375487043] 
\draw  [dash pattern={on 0.84pt off 2.51pt}]  (123,189.6) .. controls (136.53,152.52) and (482.73,157.52) .. (478.46,189.6) ;

% Text Node
\draw (384.21,65.4) node [anchor=north west][inner sep=0.75pt]    {$u_{s}$};
% Text Node
\draw (253.11,65.4) node [anchor=north west][inner sep=0.75pt]    {$u_{i}$};
% Text Node
\draw (194.72,218.4) node [anchor=north west][inner sep=0.75pt]    {$v_{j_{1}}$};
% Text Node
\draw (396.15,215.7) node [anchor=north west][inner sep=0.75pt]    {$v_{j_{2}}$};
% Text Node
\draw (295.72,223.4) node [anchor=north west][inner sep=0.75pt]    {$v_{t}$};
% Text Node
\draw (496,62.4) node [anchor=north west][inner sep=0.75pt]    {$C'$};
% Text Node
\draw (173,135.4) node [anchor=north west][inner sep=0.75pt]    {$\mathcal{E}( i,j_{1})$};
% Text Node
\draw (377,176.4) node [anchor=north west][inner sep=0.75pt]    {$\mathcal{E}( i,j_{2})$};
% Text Node
\draw (391,120.4) node [anchor=north west][inner sep=0.75pt]    {$\mathcal{E}( s,t)$};
\end{tikzpicture}

    \caption{The sets $W'_i$ must be pairwise unlinked on $S^1_{\mathrm{bot}}$.}
    \label{fig: unlinkedness}
\end{figure}
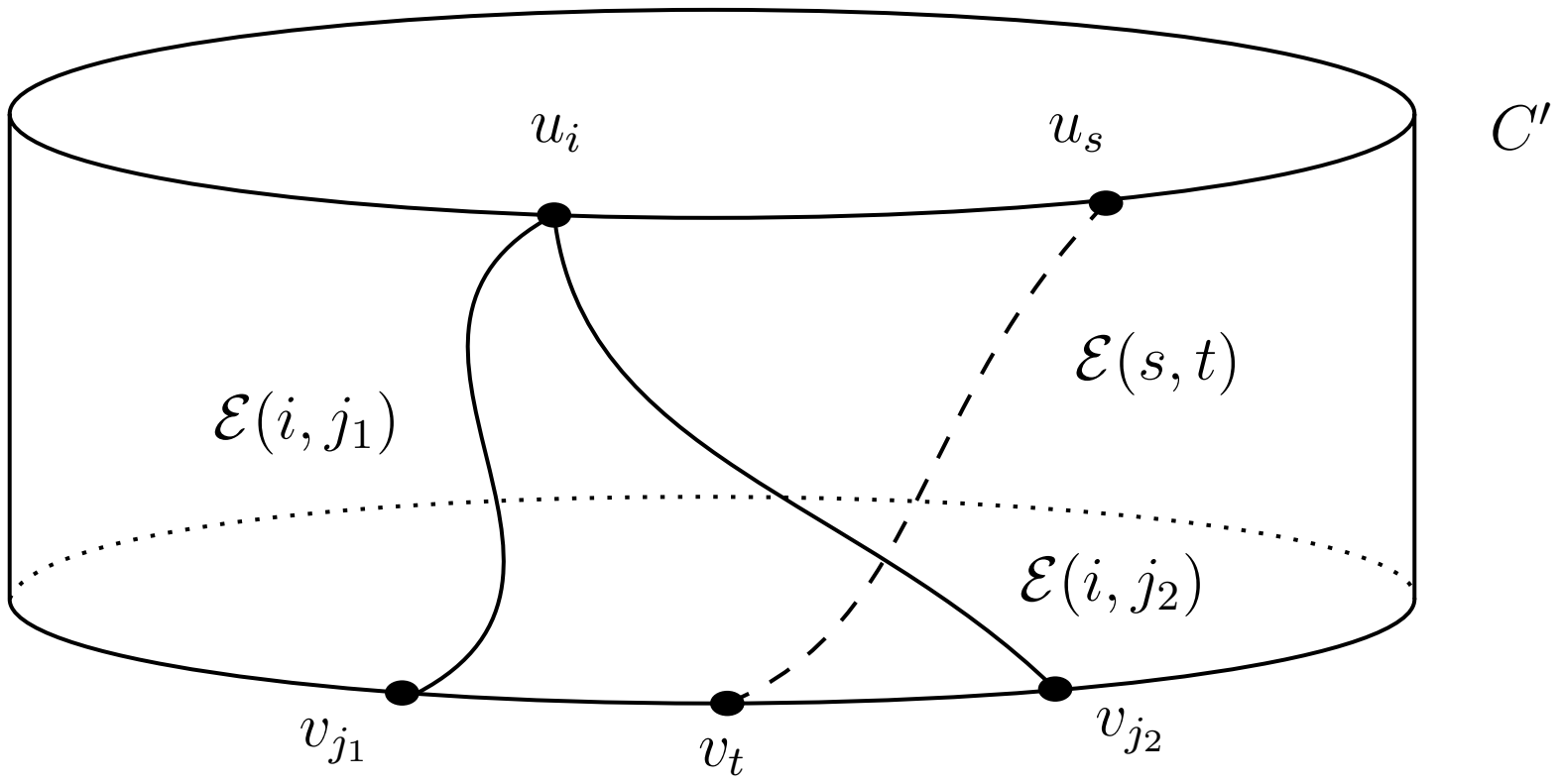

Suppose $W_i = \{v^i_{j_1}, \ldots, v^i_{j_l}\}$, where $l = q_2/q_1$, and let $W_i = V^i_{j_1} \cup \ldots \cup V^i_{j_l}$, where $V^i_{j_k}$ is the set among the $V_j \subset J_j$ corresponding to the endpoint $v^i_{j_k} \in W_i$. The sets $W_i$ are pairwise unlinked and permutted transitively on the unquotiented cylinder $C$, preserving the circular ordering, with rotation number $-p_1/q_1$. This is because the corresponding sets $W'_i$ are linked to the points $u_i$ on $S^1_{\mathrm{top}}$ with opposite orientation, and the dynamics of rays with endpoints on the sets $U_i$ has that combinatorial rotation number.

Each $W_i$ is a union of exactly $l$ of the sets $V_j$. We renumber the indices to match the circular ordering on $S^1_{\mathrm{bot}}$:
\begin{align*}
W_0 & = V_0 \cup \ldots \cup V_{l-1}, \\
W_1 & = V_l \cup \ldots \cup V_{2l - 1}, \\
\vdots & \\
W_{q_1 - 1} & = V_{(q_1 - 1)l} \cup \ldots \cup V_{q_2}.
\end{align*}
The dynamics of doubling on these sets is such that $W_i$ maps to $W_{i - p_1}$ modulo $q_1$, and $V_{j}$ maps to $V_{j + p_2}$ modulo $q_2$. Since the circular order of these sets is preserved, $W_0 = V_0 \cup \ldots \cup V_{l-1}$ maps to
$$
W_{-p_1} = V_{-p_1l} \cup \ldots \cup V_{-p_1l + l -1} = V_{p_2} \cup \ldots \cup V_{p_2 + l - 1}
$$
with matching order as well, and so $p_2 = -p_1l = -p_1q_2/q_1$ modulo $q_2$. This implies that $-p_1/q_1 = p_2/q_2$, a contradiction with the assumption that the combinatorial rotation numbers are not opposite.
\end{proof}

\section{Questions and Future Directions} \label{section: future}

We believe that the algorithm for computing direct ray connections between Hubbard trees will reveal several more interesting phenomena in matings, with a considerable amount of experimental data being available. Below are some questions and future directions of research based on the ideas explored in this paper.

\begin{question}
Is it possible to describe concisely the set of pairs of parameters $(c_1, c_2) \in M \times M$ such that there are no ray connections between the generalized trees in the formal mating, not just with respect to individual veins? What kinds of combinatorial and geometric properties will it have? What about finitely many direct ray connections, or countably many?  Theorems \ref{theorem: first ray connection}, \ref{theorem: disjointness locus open} and \ref{theorem: tuned rabbits} provide only very partial answers to this.
\end{question}

\begin{question}
Remark \ref{remark: infinitely many first ray connections} exhibits an example of a mating where ``infinitely many'' first ray connections are picked up at a parabolic parameter, despite all of them being iterated preimages of the same cycle of rays. On the other hand, Example \ref{example: tips of veins} shows matings for which there are only finitely many such ray connections. Is this always the case for tips of veins? Examining the proof of Proposition \ref{proposition: exists periodic characteristic}, there can be infinitely many direct ray connections in the grand orbit of $\R_{f_{c_1} \sqcup f_{c_2}}(\theta_1)$ only if there is a path from $\mathcal{Y}$ or $\mathcal{Z}$ to $\mathcal{X}$, so that there are dire ray connections accumulating on $\R_{f_{c_1} \sqcup f_{c_2}}(\theta_1)$ itself.
\end{question}

\begin{question}
The proof of Theorem \ref{theorem: first ray connection} implicitly establishes an upper bound on the period of the first periodic ray connection in case 3, based on the topology of the Hubbard trees (or, more fundamentally, on the angles $\theta_1$ and $\theta_2$). Can this be upper bound be sharpened in general?
\end{question}

\begin{question}
How can the mating automata be augmented to also detect long or indirect ray connections between the Hubbard trees? A possibility for this is taking products of copies of the Hubbard tree automata, ``twisted'' by the involution on arc states described in Section \ref{subsection: an involution on arc states}. We first look for arc paths $w$ in $\mathcal{M}$ such that $\iota(w)$ is also in $\mathcal{M}$. If $\psi$ is the angle associated to $\iota(w)$, we then look for an arc path $u$ in $\mathcal{H}(c_2)$ whose associated angle is $-\psi$ and see if $\iota(u)$ is in $\mathcal{M}$, where the involution is now with respect to arc states in $\mathcal{H}(c_2)$. This idea seems suitable to understand all long ray connections in a given mating, but it is unclear what this can say about families of matings.
\end{question}

\begin{question}
A natural quantity to attempt to understand in an invariant set under doubling is the entropy. In a given mating, we may compute the entropy of the set of direct ray connections between the Hubbard trees, and in the examples calculated it is lower than the minimum of the core entropies of the polynomials. Varying one of the parameters along a vein, when is the entropy first positive? Can we understand the locus of parameters pairs for which the entropy is zero, or positive?
\end{question}

\begin{question} \label{question: varying both parameters}
In the notation of Section \ref{section: disjointness locus} and in light of Example \ref{example: 1/4 with 3/8}, what can happen in general if $p_1 \prec p_1^*$ and $p_2 \prec p_2^*$? What more can be said about the disjointness locus for a pair of veins? Are the minima of the closed arcs in Lemma \ref{lemma: closed arc in veins} necessarily parabolic parameters?
\end{question}

\begin{question}
In the situations where there are ``few'' ray connections between the Hubbard trees in a geometric mating, what is the geometric and combinatorial structure of the embedded invariant graph containing the postcritical set? Is it canonical in any sense? How are these invariant graphs related to one another in the phenomenon of shared matings?
\end{question}

\begin{question}
A possible adaptation of the Hubbard tree automaton construction is one that also recognizes those angles landing at the boundaries of Fatou components. How can this automaton for the ``disked trees'' be employed to understand disjointness of immediate basins of attracting in quadratic rational maps coming from matings, in analogy with the previous questions?
\end{question}

\begin{question}
To what extent do the ideas contained in this paper extend to higher degree polynomials and their matings? For a PCF polynomial, the description of the forbidden region is the same, and adapting the definition of the Hubbard tree automaton is not complicated. On the other hand, to computationally find the forbidden region, internal addresses are only defined for unicritical polynomials $z^d + c$ of degree $d$ \cite{Schleicher2004}. Likewise, monotonicity of the characteristic lamination for degree $d$ invariant laminations would seem to require unicritically, along with known results about PCF matings.
\end{question}

\nocite{*}

\printbibliography

\end{document}